\documentclass[11pt,pdftex]{article}
\usepackage{latexsym,amssymb}
\usepackage{mathtools}
\usepackage{psfrag}
\usepackage{amsmath}
\usepackage{subfig}
\usepackage[abbrev]{amsrefs}
\usepackage{ulem}

\usepackage{tikz}
\usetikzlibrary{intersections, calc, angles, arrows.meta}
\usetikzlibrary{patterns}

\newtheorem{Theorem}{\bf Theorem}[section]
\newtheorem{Lemma}{\bf Lemma}[section]
\newtheorem{Proposition}{\bf Proposition}[section]
\newtheorem{Corollary}{\bf Corollary}[section]
\newtheorem{Remark}{\bf Remark}[section]
\newtheorem{Example}{\bf Example}[section]
\newtheorem{Definition}{\bf Definition}[section]

\newenvironment{theorem}{\begin{Theorem}$\!\!\!$}{\end{Theorem}}
\newenvironment{lemma}{\begin{Lemma}$\!\!\!$}{\end{Lemma}}
\newenvironment{proposition}{\begin{Proposition}$\!\!\!$}{\end{Proposition}}

\newenvironment{definition}{\begin{Definition}$\!\!\!$}{\end{Definition}}

\def\XXint#1#2#3{{\setbox0=\hbox{$#1{#2#3}{\int}$}
\vcenter{\hbox{$#2#3$}}\kern-.5\wd0}}

\numberwithin{equation}{section}

\allowdisplaybreaks[4]

\begin{document}

\title{Existence of solutions for a semilinear parabolic system\\ 
with singular initial data}
\author{Yohei Fujishima, Kazuhiro Ishige, and Tatsuki Kawakami}
\date{}
\maketitle
\begin{abstract}
Let $(u,v)$ be a solution to the Cauchy problem for a semilinear parabolic system
$$
\mbox{(P)}
\qquad
\left\{
\begin{array}{ll}
\partial_t u=D_1\Delta u+v^p\quad & \quad\mbox{in}\quad{\mathbb R}^N\times(0,T),\vspace{3pt}\\
\partial_t v=D_2\Delta v+u^q\quad &  \quad\mbox{in}\quad{\mathbb R}^N\times(0,T),\vspace{3pt}\\
(u(\cdot,0),v(\cdot,0))=(\mu,\nu) & \quad\mbox{in}\quad{\mathbb R}^N,
\end{array}
\right.
$$
where $N\ge 1$, $T>0$, $D_1>0$, $D_2>0$, $0<p\le q$ with $pq>1$, and 
$(\mu,\nu)$ is a pair of
nonnegative Radon measures or locally integrable nonnegative functions in ${\mathbb R}^N$.
In this paper we establish sharp sufficient conditions on the initial data for the existence of solutions to problem~(P) 
using uniformly local Morrey spaces and uniformly local weak Zygmund type spaces. 
\end{abstract}
\vspace{50pt}
\noindent Addresses:

\smallskip
\noindent Y.\,F.:  Department of Mathematical and Systems Engineering, Faculty of Engineering,\\
 Shizuoka University, 3-5-1 Johoku, Hamamatsu, Shizuoka 432-8561, Japan\\
\noindent
E-mail: {\tt fujishima@shizuoka.ac.jp}\\

\smallskip
\noindent 
K. I.: Graduate School of Mathematical Sciences, The University of Tokyo,\\ 
3-8-1 Komaba, Meguro-ku, Tokyo 153-8914, Japan\\
\noindent 
E-mail: {\tt ishige@ms.u-tokyo.ac.jp}\\

\smallskip
\noindent 
{T. K.}: Applied Mathematics and Informatics Course,\\ 
Faculty of Advanced Science and Technology, Ryukoku University,\\
1-5 Yokotani, Seta Oe-cho, Otsu, Shiga 520-2194, Japan\\
\noindent 
E-mail: {\tt kawakami@math.ryukoku.ac.jp}\\

\newpage
\section{Introduction}
We consider the Cauchy problem for a semilinear parabolic system
\begin{equation}
\tag{P}
\label{eq:P}
\left\{
\begin{array}{ll}
\partial_t u=D_1\Delta u+v^p & \quad\mbox{in}\quad{\mathbb R}^N\times(0,T),\vspace{3pt}\\
\partial_t v=D_2\Delta v+u^q & \quad\mbox{in}\quad{\mathbb R}^N\times(0,T),\vspace{3pt}\\
(u(\cdot,0),v(\cdot,0))=(\mu,\nu) & \quad\mbox{in}\quad{\mathbb R}^N,
\end{array}
\right.
\end{equation}
where $N\ge 1$, $T>0$, $D_1>0$, $D_2>0$, $0<p\le q$ with $pq>1$, 
and $(\mu,\nu)$ is a pair of nonnegative Radon measures or locally integrable nonnegative functions in ${\mathbb R}^N$.
Parabolic system~\eqref{eq:P} is the Cauchy problem for one of the simplest parabolic systems
and it is an example of reaction-diffusion systems describing heat propagation
in a two component combustible mixture.
Problem~\eqref{eq:P} has been studied extensively in many papers from various points of view. 
See e.g., \cites{AHV, EH01, EH02, FG, FI00, FI01, FI02, FIM, IKS,QS, MST, S} and references therein  
(see also \cite{QSBook}*{Chapter~32}). 
In this paper we establish sharp sufficient conditions on initial data for the existence of solutions to problem~\eqref{eq:P}.
\vspace{3pt}

We formulate the definition of solutions to problem~\eqref{eq:P}. 
Denote by ${\mathcal M}$ (resp.~${\mathcal L}$) 
the set of nonnegative Radon measures (resp.~locally integrable functions) in ${\mathbb R}^N$. 
We often identify $d\mu=\mu (x)\,dx$ in ${\mathcal M}$ for $\mu\in{\mathcal L}$.
For any $\mu\in{\mathcal M}$, let
$$
[S(t)\mu](x)\coloneqq\int_{{\mathbb R}^N} G(x-y,t)\, d\mu(y),
\quad\mbox{where}\quad
G(x,t)=(4\pi t)^{-\frac{N}{2}}\exp\left(-\frac{|x|^2}{4t}\right).
$$
\begin{definition}
\label{Definition:1.1}
Let $\mu$, $\nu\in{\mathcal M}$ and $T\in(0,\infty]$. 
Let $u$ and $v$ be nonnegative measurable and almost everywhere finite functions in ${\mathbb R}^N\times(0,T)$. 
We say that $(u,v)$ is a solution to problem~\eqref{eq:P} in ${\mathbb R}^N\times(0,T)$ if
$(u,v)$ satisfies 
\begin{equation}
\label{eq:1.1}
\begin{split}
  u(x,t) & =[S(D_1t)\mu](x)+\int_0^t [S(D_1(t-s))v(s)^p](x)\,ds,\\
  v(x,t) & =[S(D_2t)\nu](x)+\int_0^t [S(D_2(t-s))u(s)^q](x)\,ds,
\end{split}
\end{equation}
for almost all $(x,t)\in{\mathbb R}^N\times(0,T)$. 
If $(u,v)$ satisfies \eqref{eq:1.1} with $``="$ replaced by $``\ge"$, 
we say that $(u,v)$ is a supersolution to problem~\eqref{eq:P} in ${\mathbb R}^N\times(0,T)$. 
\end{definition}
For the existence of solutions to problem~\eqref{eq:P}, 
the following results have already been proved in \cites{EH02, IKS, QS} for the case of $D_1=D_2$.
\begin{itemize}
  \item[{\rm (1)}]
  Let $p\ge1$ and $r_1,r_2\in(1,\infty)$.
  Assume
  $$
  \max\{P(r_1,r_2),Q(r_1,r_2)\}\le 2,
  $$
  where
  $$
  P(r_1,r_2)\coloneqq N\left(\frac{p}{r_2}-\frac{1}{r_1}\right),
  \qquad
  Q(r_1,r_2)\coloneqq N\left(\frac{q}{r_1}-\frac{1}{r_2}\right).
 $$
  Then problem~(P) possesses a solution in ${\mathbb R}^N\times(0,T)$ for some $T>0$ if $(\mu,\nu)\in L^{r_1,\infty}\times L^{r_2,\infty}$.
  \item[{\rm (2)}]
  Assume that $\max\{P,Q\}>2$. Then there exists $(\mu,\nu)\in L^{r_1}\times L^{r_2}$
  such that problem~(P) possess no local-in-time solutions.
  \item[{\rm (3)}]
  Assume that 
  \begin{equation}
  \label{eq:1.2}
  \frac{q+1}{pq-1}<\frac{N}{2}
  \end{equation}
  and both $\|\mu\|_{L^{r_1^*,\infty}}$ and $\|\nu\|_{L^{r_2^*,\infty}}$ are small enough,
  where
  \begin{equation*}
  r_1^*\coloneqq\frac{N}{2}\frac{pq-1}{p+1},
  \qquad
  r_2^*\coloneqq\frac{N}{2}\frac{pq-1}{q+1}.
  \end{equation*}
  Then problem~(P) possesses a global-in-time solution. 
  On the other hand,
  if $(p,q)$ does not satisfy \eqref{eq:1.2}, then
  problem~(P) possesses no global-in-time non-trivial solutions.
\end{itemize}
Subsequently, in \cite{FI01} the first and the second authors of this paper divided problem~\eqref{eq:P} into the following six cases: 
\begin{equation*}
\begin{split}
 & {\rm (A)}\quad\frac{q+1}{pq-1}<\frac{N}{2};\\
 & {\rm (B)}\quad \frac{q+1}{pq-1}=\frac{N}{2}\quad\mbox{and}\quad p<q ;
  \qquad\qquad\,\,
	{\rm (C)}\quad \frac{q+1}{pq-1}=\frac{N}{2}\quad\mbox{and}\quad p=q ;\\
 & {\rm (D)}\quad \frac{q+1}{pq-1}>\frac{N}{2}\quad\mbox{and}\quad q>1+\frac{2}{N};
	\qquad
	{\rm (E)}\quad \frac{q+1}{pq-1}>\frac{N}{2}\quad\mbox{and}\quad q=1+\frac{2}{N};\\
 & {\rm (F)}\quad \frac{q+1}{pq-1}>\frac{N}{2}\quad\mbox{and}\quad q<1+\frac{2}{N},
\end{split}
\end{equation*}
\begin{figure}[htbp]
	\centering
	\subfloat{
		\begin{tikzpicture}[samples = 100]
			\draw[-{Latex}] (-0.5,0) -- (5.15,0) node[right] {\small $p$};
			\draw[-{Latex}] (0,-0.5) -- (0,5.15) node[above] {\small $q$};
			\path (0,0) node[above left] {\small O};
			\fill [ opacity = 0.2 ] plot [ domain = 0.2:1 ] ({\x}, {1/\x}) -- plot [ domain = 1:5 ] ({\x}, {\x});
			\draw plot[ domain = -0.4:5 ] ({\x}, {\x});
			\path (4,4) node[below right] {\small $p=q$};
      \draw [dashed] (3/5,5/3) -- (0,5/3) node[left] {\small $1+\frac{2}{N}$};
			\draw[color = white, thick] plot[ domain = 0.2:1 ] ({\x}, {1/\x});
			\draw[dashed] plot[ domain = 0.2:5 ] ({\x}, {1/\x});
			\draw[color = white, thick] plot[ domain = 1:5/3 ] ({\x}, {(5/3) / (\x - 2/3)});
			\draw[dashed] plot[ domain = 1:5/3 ] ({\x}, {(5/3) / (\x - 2/3)});
			\draw plot[ domain = 5/3:5 ] ({\x}, {(5/3) / (\x - 2/3)});
			\path (4.2,1.2) node {\small $\dfrac{q+1}{pq-1}=\dfrac{N}{2}$};
      \path (-0.8, 3.5) node {\small $pq=1$};
      \path (-0.8, 3.2) edge [-{Latex}, bend right] (1/2.7,2.7);
      \draw [thick, color = white] (3/5,5/3) -- (5/3,5/3);
			\draw [dashed] (3/5,5/3) -- (5/3,5/3);
      \draw [dashed] (5/3,5/3) -- (5/3,0) node[below] {\small $1+\frac{2}{N}$};
			\fill [color = white] (5/3, 5/3) circle (2pt);
			\draw (5/3, 5/3) circle (2pt);
			\fill [ color = white ] (1,1) circle (2pt);
			\draw (1,1) circle (2pt);
			\path (0.85,2.8) node {\small (D)};
			\path (2.3,3.3) node[above] {\small (A)};
			\path (1.03,1.1) node[above] {\small (F)};
		\end{tikzpicture}
	}
	\quad
	\subfloat{
		\begin{tikzpicture}[samples = 100]
			\draw[-{Latex}] (-0.5,0) -- (5.15,0) node[right] {\small $p$};
			\draw[-{Latex}] (0,-0.5) -- (0,5.15) node[above] {\small $q$};
			\path (0,0) node[above left] {\small O};
			\draw plot[ domain = -0.4:5 ] ({\x}, {\x});
			\path (4,4) node[below right] {\small $p=q$};
			\draw[dashed] plot[ domain = 0.2:5 ] ({\x}, {1/\x});
			\draw plot[ domain = 1:5 ] ({\x}, {(5/3) / (\x - 2/3)});
			\path (4.2,1.2) node {\small $\dfrac{q+1}{pq-1}=\dfrac{N}{2}$};
      \path (-0.8, 3.5) node {\small $pq=1$};
      \path (-0.8, 3.2) edge [-{Latex}, bend right] (1/2.7,2.7);
      \draw [dashed] (5/3,5/3) -- (5/3,0) node[below] {\small $1+\frac{2}{N}$};
      \draw [dashed] (3/5,5/3) -- (0,5/3) node[left] {\small $1+\frac{2}{N}$};
			\fill (5/3, 5/3) circle (2pt);
			\fill [ color = white ] (1,1) circle (2pt);
			\draw (1,1) circle (2pt);
			\draw (3/5,5/3) -- (5/3,5/3);
			\path (2.8,1.9) edge [-{Latex}, bend left = 30] (1.74,5/3);
			\path (2.84,1.8) node[above] {\small (C)};
			\path (0.8,2.3) edge [-{Latex}, bend right = 15] (1,1.7);
			\path (0.8,2.2) node[above] {\small (E)};
			\path (2,3.5) edge [-{Latex}, bend left] (1.25,3);
			\path (2,3.4) node[above] {\small (B)};
		\end{tikzpicture}
	}
	\caption{}\label{figure:1}
\end{figure}

\noindent
and obtained 
necessary conditions for the existence of solutions to problem~\eqref{eq:P}.
Subsequently, in \cite{FI02} they studied sufficient conditions for the existence of solutions to problem~\eqref{eq:P}, 
and identified the optimal singularity of the initial data for the existence of solutions to problem~\eqref{eq:P} (see \cite{FI02}*{Theorem~1.2}). 
\begin{proposition}
  \label{Proposition:1.1}
  Let $N\ge 1$ and $0<p\le q$ with $pq>1$.
  \begin{itemize}
    \item[{\rm (a)}]
    Consider case {\rm (A)}. Let
    \begin{align*}
      & \mu(x)=c_{a,1} |x|^{-\frac{2(p+1)}{pq-1}}\chi_{B(0,1)}(x)\quad\mbox{in}\quad{\mathbb R}^N,\\
      & \nu(x)=c_{a,2} |x|^{-\frac{2(q+1)}{pq-1}}\chi_{B(0,1)}(x)\quad\mbox{in}\quad{\mathbb R}^N,
    \end{align*}
    where $c_{a,1}$, $c_{a,2}>0$.
    Problem~{\rm (P)} possesses no local-in-time solutions
    if either $c_{a,1}$ or $c_{a,2}$ is large enough.
    On the other hand,
    problem~{\rm (P)} possesses a global-in-time solution
    if both of $c_{a,1}$ and $c_{a,2}$ are small enough.
    \item[{\rm (b)}]
    Consider case {\rm (B)}. Let
    \begin{align*}
      & \mu(x)=c_{b,1} |x|^{-\frac{2(p+1)}{pq-1}}\left[\log\left(e+\frac{1}{|x|}\right)\right]^{-\frac{p}{pq-1}}\chi_{B(0,1)}(x)
      \quad\mbox{in}\quad{\mathbb R}^N,\\
      & \nu(x)=c_{b,2} |x|^{-N}\left[\log\left(e+\frac{1}{|x|}\right)\right]^{-\frac{1}{pq-1}-1}\chi_{B(0,1)}(x)
     \quad\,\,\,\,\mbox{in}\quad{\mathbb R}^N,
    \end{align*}
    where $c_{b,1}$, $c_{b,2}>0$.
    Problem~{\rm (P)} possesses no local-in-time solutions
    if either $c_{b,1}$ or $c_{b,2}$ is large enough.
    On the other hand,
    problem~{\rm (P)} possesses a local-in-time solution
    if both of $c_{b,1}$ and $c_{b,2}$ are small enough.
    \item[{\rm (c)}]
    Consider case {\rm (C)}. Let
    \begin{align*}
      & \mu(x)=c_{c,1}|x|^{-N}\left[\log\left(e+\frac{1}{|x|}\right)\right]^{-\frac{N}{2}-1}\chi_{B(0,1)}(x)\quad\mbox{in}\quad{\mathbb R}^N,\\
      & \nu(x)=c_{c,2}|x|^{-N}\left[\log\left(e+\frac{1}{|x|}\right)\right]^{-\frac{N}{2}-1}\chi_{B(0,1)}(x)\quad\mbox{in}\quad{\mathbb R}^N,
    \end{align*}
    where $c_{c,1}$, $c_{c,2}>0$.
    Problem~{\rm (P)} possesses no local-in-time solutions
    if either $c_{c,1}$ or $c_{c,2}$ is large enough.
    On the other hand,
    problem~{\rm (P)} possesses a local-in-time solution
    if both of $c_{c,1}$ and $c_{c,2}$ are small enough.
    \item[{\rm (d)}]
    Consider case {\rm (D)}.
    Let
    $$
    \mu(x)=|x|^{-\frac{N+2}{q}}h_1(|x|)\chi_{B(0,1)}(x)\quad\mbox{in}\quad{\mathbb R}^N,
    $$
    where $h_1$ is a positive increasing function in $(0,1]$
    such that $h_1(1)<\infty$
    and $r^{-\epsilon}h_1(r)$ is decreasing in $r$ for some $\epsilon>0$.
    Let $\nu\in{\mathcal M}$. 
    Problem~{\rm (P)} possesses no local-in-time solutions if either
    $$
    \int_0^1 h_1(\tau)^q\tau^{-1}\,d\tau=\infty
    \quad\mbox{or}\,\,\,
    \sup_{x\in{\mathbb R}^N}\nu(B(x,1))=\infty.
    $$
    On the other hand,
    problem~{\rm (P)} possesses a local-in-time solution if
    $$
    \int_0^1 h_1(\tau)^q\tau^{-1}\,d\tau<\infty
    \quad\mbox{and}\,\,\,
    \sup_{x\in{\mathbb R}^N}\nu(B(x,1))<\infty.
    $$
    \item[{\rm (e)}]
    Consider case {\rm (E)}.
    Let
    $$
    \mu(x)=|x|^{-N}h_2(|x|)\chi_{B(0,1)}(x)\quad\mbox{in}\quad{\mathbb R}^N,
    $$
    where $h_2$ is a positive increasing function in $(0,1]$ such that $h_2(1)<\infty$.
    Let $\nu\in{\mathcal M}$. 
    Problem~{\rm (P)} possesses no local-in-time solutions if either
    $$
    \int_0^1\left[\int_0^r h_2(\tau)\tau^{-1}\,d\tau\right]^qr^{-1}\,dr=\infty
    \quad\mbox{or}\,\,\,
    \sup_{x\in{\mathbb R}^N}\nu(B(x,1))=\infty.
    $$
    On the other hand,
    problem~{\rm (P)} possesses a local-in-time solution if
    $$
    \int_0^1\left[\int_0^r h_2(\tau)\tau^{-1}\,d\tau\right]^qr^{-1}\,dr<\infty
    \quad\mbox{and}\,\,\,
    \sup_{x\in{\mathbb R}^N}\nu(B(x,1))<\infty.
    $$
     \item[{\rm (f)}]
    Consider case {\rm (F)}.
    Let $\mu$, $\nu\in{\mathcal M}$. 
    Problem~{\rm (P)} possesses no local-in-time solutions if either
    $$
    \sup_{x\in{\mathbb R}^N}\mu(B(x,1))=\infty
    \quad\mbox{or}\,\,\,
    \sup_{x\in{\mathbb R}^N}\nu(B(x,1))=\infty.
    $$
    On the other hand, problem~{\rm (P)} possesses a local-in-time solution if
    $$
    \sup_{x\in{\mathbb R}^N}\mu(B(x,1))<\infty
    \quad\mbox{and}\,\,\,
    \sup_{x\in{\mathbb R}^N}\nu(B(x,1))<\infty.
    $$
  \end{itemize}
\end{proposition}
Proposition~\ref{Proposition:1.1} with cases~(A), (C), and (F) can be regarded as a generalization of 
\cite{HI01}*{Corollary~1.2}~(ii), (i), and \cite{HI01}*{Theorem~1.3}, respectively, for the scalar semilinear parabolic equation 
$\partial_t w=\Delta w+w^p$, where $p>1$. (See also \cites{BP, FHIL}.) 
On the one hand, optimal singularities of the initial data in Proposition~\ref{Proposition:1.1} with cases~(B), (D), and (E) are peculiar to the parabolic system.

In this paper, taking into the account of Proposition~\ref{Proposition:1.1}, 
we obtain sharp sufficient conditions on the initial data for the existence of solutions to problem~(P) 
in the framework of Banach spaces. 
In cases~(A) and (F), we develop the arguments in \cite{FI02}*{Section~3} and \cite{IKS} to 
obtain our sharp sufficient conditions using uniformly local Morrey spaces 
(see Theorems~\ref{Theorem:1.1} and \ref{Theorem:1.2}). 

For the other cases (B)--(E), we develop the arguments in \cites{IIK} 
to introduce new uniformly local weak Zygmund type spaces. 
In \cite{IIK} the second and the third authors of this paper and Ioku introduced 
a uniformly local weak Zygmund type space ${\mathfrak L}_{{\rm ul}}^{r,\infty}(\log {\mathfrak L})^\alpha$, 
where $1\le r\le\infty$ and $0\le\alpha<\infty$, 
to obtain sharp sufficient conditions for the existence of solutions to the Cauchy problem 
for the critical fractional semilinear heat equation
$$
\partial_t u+(-\Delta)^{\frac{\theta}{2}}u=|u|^{\frac{\theta}{N}}u\quad\mbox{in}\quad{\mathbb R}^N\times(0,T),
\quad
u(\cdot,0)=\mu\quad\mbox{in}\quad{\mathbb R}^N,
$$
where $\theta\in(0,2]$.  
For the proof, they established sharp decay estimates of the fractional heat semigroup in ${\mathfrak L}_{{\rm ul}}^{r,\infty}(\log {\mathfrak L})^\alpha$.  
In this paper, applying the arguments in \cites{FI02, IIK}, we obtain sharp sufficient conditions for the existence of solutions to problem~\eqref{eq:P} in case~(C)
(see Theorem~\ref{Theorem:1.4}). 

For cases~(B), (D), and (E), 
in addition to ${\mathfrak L}_{{\rm ul}}^{r,\infty}(\log {\mathfrak L})^\alpha$, 
we treat somewhat standard uniformly local weak Zygmund type space $L_{{\rm ul}}^{r,\infty}(\log L)^\alpha$ 
and we also introduce more general uniformly local weak Zygmund type spaces $L_{{\rm ul}}^{r,\infty}\Phi(L)^\alpha$ and ${\mathfrak L}_{{\rm ul}}^{r,\infty}\Phi({\mathfrak L})^\alpha$. 
Then we establish sharp decay estimates of the heat semigroup  in these uniformly local weak Zygmund type spaces (see Proposition~\ref{Proposition:3.1}). 
Furthermore, we develop the arguments in \cites{FI02, IIK} to get uniform estimates of approximate solutions in suitable uniformly local weak Zygmund type spaces, 
and obtain sharp sufficient conditions for the existence of solutions in cases~(B), (D), and~(E).
\medskip

We introduce some notation. 
For any measurable set $E$ in ${\mathbb R}^N$, we denote by $\chi_E$ (resp.~$|E|$) the characteristic function of $E$ 
(resp.~the $N$-dimensional Lebesgue measure of $E$). 
For any $x\in{\mathbb R}^N$ and $R>0$, let $B(x,R)\coloneqq\{y\in{\mathbb R}^N\,:\,|x-y|<R\}$. 
Set $\omega_N\coloneqq|B(0,1)|$. 
For any $r\in[1,\infty]$, we denote by $\|\cdot\|_{L^r}$ the usual norm of $L^r$. 
For any $\mu\in{\mathcal M}$, we say that 
$\mu\in{\mathcal M}_{{\rm ul}}$ if 
$$
\|\mu\|_{{\mathcal M}_{{\rm ul}}}\coloneqq\sup_{x\in{\mathbb R}^N}\mu(B(x,1))<\infty. 
$$
Similarly, for any $f\in{\mathcal L}$ and $r\in[1,\infty]$, we say that $f\in L^r_{{\rm ul}}$ if 
$$
\|f\|_{L^r_{{\rm ul}}}\coloneqq\sup_{x\in{\mathbb R}^N}\|f\chi_{B(x,1)}\|_{L^r}<\infty.
$$
For any measurable function $f$ in ${\mathbb R}^N$, 
we denote by $\mu_f$ the distribution function of~$f$, that is, 
$$
\mu_f(\lambda)\coloneqq\left|\{x\,:\,|\,f(x)|>\lambda\}\right|,\quad \lambda > 0.
$$
We define the non-increasing rearrangement $f^*$ of $f$ by 
$$
f^{*}(s)\coloneqq\inf\{\lambda>0\,:\,\mu_f(\lambda)\le s\},\quad s\in[0,\infty).
$$
Here we adopt the convention $\inf\emptyset=\infty$. 
Then $f^*$ is non-increasing and right continuous in $[0,\infty)$, and it has the following properties (see \cite{Grafakos}*{Proposition~1.4.5}):
\begin{equation}
\label{eq:1.3}
\begin{split}
 & (kf)^*=|k|f^*,\qquad (|f|^q)^*=(f^*)^q,\qquad \|f^*\|_{L^r((0,\infty))}=\|f\|_{L^r},
\end{split}
\end{equation}
where $q\in(0,\infty)$, $k\in{\mathbb R}$, and $r\in[1,\infty]$.
For any $r\in[1,\infty]$, we define the weak $L^r$ space by
$$
L^{r,\infty}\coloneqq\left\{f\in {\mathcal L}\,:\,\|f\|_{L^{r,\infty}}\coloneqq\sup_{s>0}\left\{s^{\frac{1}{r}}f^*(s)\right\}<\infty\right\}. 
$$
Then $L^{\infty,\infty}=L^\infty$ and $L^r\subsetneq L^{r,\infty}$ if $1<r<\infty$.

Next, we introduce uniformly local Morrey spaces. 
For any $r\in[1,\infty]$, $\alpha\in[1,r]$, and $R\in(0,\infty]$, let 
\begin{equation}
\label{eq:1.4}
\|f\|_{M(r,\alpha;R)}\coloneqq\sup_{x\in{\mathbb R}^N}\sup_{\sigma\in(0,R)}\left\{|B(x,\sigma)|^{\frac{1}{r}-\frac{1}{\alpha}}
\|f\|_{L^\alpha(B(x,\sigma))}\right\},
\quad f\in {\mathcal L}.
\end{equation}
We write $\|f\|_{M(r,\alpha)}\coloneqq\|f\|_{M(r,\alpha;1)}$ for simplicity. 
We define the uniformly local Morrey space $M(r,\alpha)$ by 
$$
M(r,\alpha)\coloneqq\left\{f\in {\mathcal L}\,:\,\|f\|_{M(r,\alpha)}<\infty\right\}.
$$
Then $M(r,\alpha)$ is a Banach space equipped with the norm~$\|\,\cdot\,\|_{M(r,\alpha)}$. 
Notice that $M(\infty,\alpha)=L^\infty$ and
$$
\|\,\cdot\,\|_{M(\infty,\alpha; R)}=\|\cdot\|_{L^\infty}
$$ 
for $\alpha\in[1,\infty]$ and $R\in(0,\infty]$.
\vspace{3pt}

Now we state our main results in cases~(A) and (F). 
\begin{theorem}
\label{Theorem:1.1} 
Consider case~{\rm (A)}. 
Let 
\begin{equation}
\label{eq:1.5}
r_1^*\coloneqq\frac{N}{2}\frac{pq-1}{p+1},
\quad
r_2^*\coloneqq\frac{N}{2}\frac{pq-1}{q+1},
\quad 
\alpha_A\coloneqq\frac{q+1}{p+1}\beta_A,
\quad
1<\beta_A<\frac{q(p+1)}{q+1},
\quad
\beta_A\le r_2^*.
\end{equation}
Then there exists $\delta_A>0$ such that, if a pair $(\mu,\nu)\in{\mathcal L}\times{\mathcal L}$ satisfies
\begin{equation}
\label{eq:1.6}
\|\mu\|_{M(r_1^*,\alpha_A;T^{\frac{1}{2}})}^{\alpha_A}+\|\nu\|_{M(r_2^*,\beta_A;T^{\frac{1}{2}})}^{\beta_A}\le\delta_A
\end{equation}
for some $T\in(0,\infty]$, 
then there exists a solution $(u,v)$ to problem~\eqref{eq:P} in ${\mathbb R}^N\times(0,T)$ 
such that 
\begin{align}
\label{eq:1.7}
 &  \sup_{t\in(0,T)}\|u(t)\|_{M(r_1^*,\alpha_A;T^{\frac{1}{2}})}+\sup_{t\in(0,T)}\left\{t^{\frac{N}{2r_1^*}}\|u(t)\|_{L^\infty}\right\}<\infty,\\
\label{eq:1.8}
 & \sup_{t\in(0,T)}\|v(t)\|_{M(r_2^*,\beta_A;T^{\frac{1}{2}})}+\sup_{t\in(0,T)}\left\{t^{\frac{N}{2r_2^*}}\|v(t)\|_{L^\infty}\right\}<\infty,\\
\label{eq:1.9}
 &
 \lim_{t\to +0}\|u(t)-S_1(D_1t)\mu\|_{M(r_1,\ell_1)}=0,
 \quad
 \lim_{t\to +0}\|v(t)-S(D_2t)\nu\|_{M(r_2,\ell_2)}=0,
\end{align}
where $r_1\in[\alpha_A^{-1}r_1^*,r_1^*)$, $r_2\in[\beta_A^{-1}r_2^*,r_2^*)$, $\ell_1\in[1,\alpha_Ar_1/r_1^*]$, and $\ell_2\in[1,\beta_Ar_2/r_2^*]$.
\end{theorem}
Notice that, in case (A), we have $r_1^*\ge r_2^*>1$ by $p\le q$ and $\alpha_A\le r_1^*$ by $\beta_A\le r_2^*$. 
\begin{theorem}
\label{Theorem:1.2} 
Consider case~{\rm (F)}. 
Assume $\mu$, $\nu\in{\mathcal M}_{{\rm ul}}$. 
Then there exists a solution $(u,v)$ to problem~\eqref{eq:P} in ${\mathbb R}^N\times(0,T)$ 
for some $T\in(0,\infty)$ 
such that 
\begin{equation}
\label{eq:1.10}
\sup_{t\in(0,T)}\left\{\|u(t)\|_{L^1_{{\rm ul}}}+t^{\frac{N}{2}}\|u(t)\|_{L^\infty}\right\}<\infty,
\quad
\sup_{t\in(0,T)} \left\{\|v(t)\|_{L^1_{{\rm ul}}}+t^{\frac{N}{2}}\|v(t)\|_{L^\infty}\right\}<\infty.
\end{equation}
Furthermore, 
\begin{equation}
\label{eq:1.11}
\lim_{t\to+0}\left(\|u(t)-S(D_1t)\mu\|_{L^1_{{\rm ul}}}+\|v(t)-S(D_2t)\nu\|_{L^1_{{\rm ul}}}\right)=0.
\end{equation}
\end{theorem}
We discuss the optimality of Theorems~\ref{Theorem:1.1} and \ref{Theorem:1.2} in Section~7.

Next, we introduce weak Zygmund type spaces to obtain our sufficient conditions for the existence of solutions 
to problem~\eqref{eq:P} in cases~(B)--(E).
Throughout this paper, let $\Phi$ be a non-decreasing function in $[0,\infty)$ with the following properties:
\begin{itemize}
 \item[($\Phi$1)] $\Phi(0)=1$;
  \item[($\Phi$2)] there exists $C>0$ such that
  $\Phi(a^2)\le C\Phi(a)$ for $a\ge 0$; 
  \item[($\Phi$3)]  
  for any $\delta>0$, there exist $C_\delta>0$ and $\tau_\delta>0$ such that 
  $$
  \tau_2^{-\delta}\Phi(\tau_2)\le C_\delta \tau_1^{-\delta}\Phi(\tau_1)\quad\mbox{if}\quad \tau_\delta\le\tau_1\le\tau_2.
  $$
\end{itemize}
For any $r\in[1,\infty]$ and $\alpha\in[0,\infty)$, 
we define weak Zygmund type spaces $L^{r,\infty}\Phi(L)^\alpha$ and ${\mathfrak L}^{r,\infty}\Phi({\mathfrak L})^\alpha$ by
$$
L^{r,\infty}\Phi(L)^\alpha\coloneqq\{f\in {\mathcal L}\,:\,\|f\|_{L^{r,\infty}\Phi(L)^\alpha}<\infty\},
\,\,\,\,
{\mathfrak L}^{r,\infty}\Phi({\mathfrak L})^\alpha\coloneqq\{f\in {\mathcal L}\,:\,\|f\|_{{\mathfrak L}^{r,\infty}\Phi({\mathfrak L})^\alpha}<\infty\},
$$
respectively, where 
\begin{equation*}
\begin{array}{ll}
\|f\|_{L^{r,\infty}\Phi(L)^\alpha}\coloneqq\displaystyle{\sup_{s>0}}\,\left\{s\Phi(s^{-1})^\alpha f^{*}(s)^{r}\right\}^{\frac{1}{r}}\quad & \mbox{if}\quad r<\infty,\vspace{3pt}\\
\|f\|_{{\mathfrak L}^{r,\infty}\Phi({\mathfrak L})^\alpha}\coloneqq\displaystyle{\sup_{s>0}}\,\left\{s\Phi(s^{-1})^\alpha (|f|^r)^{**}(s)\right\}^{\frac{1}{r}}\quad & \mbox{if}\quad r<\infty,\vspace{3pt}\\
\|f\|_{L^{r,\infty}\Phi(L)^\alpha}\coloneqq L^\infty,\quad \|f\|_{{\mathfrak L}^{r,\infty}\Phi({\mathfrak L})^\alpha}\coloneqq L^\infty \quad & \mbox{if}\quad r=\infty. 
\end{array}
\end{equation*}
Here
$$
f^{**}(s)\coloneqq\frac{1}{s}\int_0^s f^*(\tau)\,d\tau,\quad s\in(0,\infty).
$$
Similarly to \cite{IIK}*{Lemma~2.1}, 
we see that ${\mathfrak L}^{r,\infty}\Phi({\mathfrak L})^\alpha$ is a Banach space equipped with the norm $\|\cdot\|_{{\mathfrak L}^{r,\infty}\Phi({\mathfrak L})^\alpha}$.
$L^{r,\infty}\Phi(L)^\alpha$ is also a Banach space if $r>1$ (see Lemma~\ref{Lemma:3.9}). 
Furthermore, 
$$
L^{r,\infty}\Phi(L)^0=L^{r,\infty},\qquad 
{\mathfrak L}^{r,\infty}\Phi({\mathfrak L})^0=L^r,\qquad
{\mathfrak L}^{r,\infty}\Phi({\mathfrak L})^\alpha\subset L^{r,\infty}\Phi(L)^\alpha.
$$
In the case of $\Phi(\tau)=\log(e+\tau)$, we write 
$$
L^{r,\infty}(\log L)^\alpha:=L^{r,\infty}\Phi(L)^\alpha,\quad
{\mathfrak L}^{r,\infty}(\log {\mathfrak L})^\alpha:={\mathfrak L}^{r,\infty}\Phi({\mathfrak L})^\alpha,
$$
for simplicity.

Next, we define uniformly local weak Zygmund type spaces $L^{r,\infty}_{{\rm ul}}\Phi(L)^\alpha$ and ${\mathfrak L}^{r,\infty}_{{\rm ul}}\Phi({\mathfrak L})^\alpha$. 
For any $R\in(0,\infty]$, set 
$$
\|f\|_{\Phi,r,\alpha;R}\coloneqq\sup_{x\in{\mathbb R}^N}\|f\chi_{B(x,R)}\|_{L^{r,\infty}\Phi(L)^\alpha},
\quad
|||f|||_{\Phi,r,\alpha;R}\coloneqq\sup_{x\in{\mathbb R}^N}\|f\chi_{B(x,R)}\|_{{\mathfrak L}^{r,\infty}\Phi({\mathfrak L})^\alpha}.
$$
Then 
$\|f\|_{L^{r,\infty}\Phi(L)^\alpha}=\|f\|_{\Phi,r,\alpha;\infty}$ and $\|f\|_{{\mathfrak L}^{r,\infty}\Phi({\mathfrak L})^\alpha}=|||f|||_{\Phi,r,\alpha;\infty}$. 
We write 
$$
\|f\|_{\Phi,r,\alpha}:=\|f\|_{\Phi,r,\alpha;1},
\quad 
|||f|||_{\Phi,r,\alpha}:=|||f|||_{\Phi,r,\alpha;1},
$$
for simplicity. 
Then we define 
$$
L_{{\rm ul}}^{r,\infty}\Phi(L)^\alpha\coloneqq\{f\in {\mathcal L}\,:\,\|f\|_{\Phi,r,\alpha}<\infty\},
\,\,\,\,
{\mathfrak L}_{{\rm ul}}^{r,\infty}\Phi({\mathfrak L})^\alpha\coloneqq\{f\in {\mathcal L}\,:\,|||f|||_{\Phi,r,\alpha}<\infty\}.
$$
We remark that 
\begin{equation}
\label{eq:1.12}
{\mathfrak L}_{{\rm ul}}^{r,\infty}\Phi({\mathfrak L})^0=L_{{\rm ul}}^r. 
\end{equation}
In the case of $\Phi(\tau)=\log(e+\tau)$, 
we write
\begin{equation*}
\begin{array}{ll}
L_{{\rm ul}}^{r,\infty}(\log L)^\alpha:=L^{r,\infty}_{{\rm ul}}\Phi(L)^\alpha,\quad
 & {\mathfrak L}_{{\rm ul}}^{r,\infty}(\log {\mathfrak L})^\alpha:={\mathfrak L}_{{\rm ul}}^{r,\infty}\Phi({\mathfrak L})^\alpha,\vspace{3pt}\\
\|\cdot\|_{r,\alpha;R}:=\|\cdot\|_{\Phi,r,\alpha;R},\quad 
 & |||\cdot|||_{r,\alpha;R}:=|||\cdot|||_{\Phi,r,\alpha;R},\\
 \|\cdot\|_{r,\alpha}:=\|\cdot\|_{\Phi,r,\alpha},\quad 
 & |||\cdot|||_{r,\alpha}:=|||\cdot|||_{\Phi,r,\alpha},
\end{array}
\end{equation*}
for simplicity. 
\vspace{3pt}%

Now we are ready to state our main results in cases~(B)--(E).
\begin{theorem}
\label{Theorem:1.3} 
Consider case~{\rm (B)}. 
Let
\begin{equation}
\label{eq:1.13}
\alpha_B\coloneqq\frac{q+1}{p+1}\frac{p}{pq-1},\quad
\beta_B\coloneqq\frac{1}{pq-1}.
\end{equation}
For any $T_*\in(0,\infty)$, 
there exists $\delta_B>0$ such that if $(\mu,\nu)\in{\mathcal L}\times{\mathcal L}$ satisfies
\begin{equation}
\label{eq:1.14}
\|\mu\|_{\frac{q+1}{p+1},\alpha_B; T^{\frac{1}{2}}}+|||\nu|||_{1,\beta_B; T^{\frac{1}{2}}}\le\delta_B
\end{equation}
for some $T\in(0,T_*]$, 
then there exists a solution $(u,v)$ to problem~\eqref{eq:P} in ${\mathbb R}^N\times(0,T)$ 
such that 
\begin{equation}
\label{eq:1.15}
\begin{split}
 & \sup_{t\in(0,T)}\|u(t)\|_{\frac{q+1}{p+1},\alpha_B;T^{\frac{1}{2}}}
 +\sup_{t\in(0,T)} \left\{t^{\frac{N}{2}\frac{p+1}{q+1}}\left[\log\left(e+\frac{1}{t}\right)\right]^{\frac{p}{pq-1}}\|u(t)\|_{L^\infty}\right\}<\infty,\\
 & \sup_{t\in(0,T)}|||v(t)|||_{1,\beta_B; T^{\frac{1}{2}}}
 +\sup_{t\in(0,T)}\left\{t^{\frac{N}{2}}\left[\log\left(e+\frac{1}{t}\right)\right]^{\frac{1}{pq-1}}\|v(t)\|_{L^\infty}\right\}<\infty.
\end{split}
\end{equation}
Furthermore, 
\begin{equation}
\label{eq:1.16}
\lim_{t\to +0} \|u(t)-S(D_1t)\mu\|_{\frac{q+1}{p+1},\alpha;T^{\frac{1}{2}}}=0,
\quad
\lim_{t\to +0}|||v(t)-S(D_2t)\nu|||_{1,\beta;T^{\frac{1}{2}}}=0,
\end{equation}
for $\alpha\in[0,\alpha_B)$ and $\beta\in[0,\beta_B)$. 
\end{theorem}
\begin{theorem}
\label{Theorem:1.4} 
Consider case~{\rm (C)}. 
For any $T_*\in(0,\infty)$, 
there exists $\delta_C>0$ such that if $(\mu,\nu)\in{\mathcal L}\times{\mathcal L}$ satisfies
\begin{equation}
\label{eq:1.17}
|||\mu|||_{1,\frac{N}{2};T^{\frac{1}{2}}}+|||\nu|||_{1,\frac{N}{2};T^{\frac{1}{2}}}\le\delta_C
\end{equation}
for some $T\in(0,T_*]$, then there exists a solution $(u,v)$ to problem~\eqref{eq:P} in ${\mathbb R}^N\times(0,T)$ 
such that 
\begin{align*}
 & \sup_{0<t<T}\left\{|||u(t)|||_{1,\frac{N}{2};T^{\frac{1}{2}}}+|||v(t)|||_{1,\frac{N}{2};T^{\frac{1}{2}}}\right\}<\infty,\\
 & \sup_{0<t<T}\left\{t^{\frac{N}{2}}\left[\log\left(e+\frac{1}{t}\right)\right]^{\frac{N}{2}}\left(\|u(t)\|_{L^\infty}+\|v(t)\|_{L^\infty}\right)\right\}<\infty.
\end{align*}
Furthermore, the solution~$(u,v)$ satisfies 
$$
\lim_{t\to +0} |||u(t)-S(D_1t)\mu|||_{1,\gamma; T^{\frac{1}{2}}}=0,
\quad
\lim_{t\to +0}|||v(t)-S(D_2t)\nu|||_{1,\gamma; T^{\frac{1}{2}}}=0,
$$
for $\gamma\in[0,N/2)$.
\end{theorem}
\begin{theorem}
\label{Theorem:1.5} 
Consider case~{\rm (D)}. 
Let $\Phi$ be a non-decreasing function in $[0,\infty)$ with properties~{\rm ($\Phi$1)}--{\rm ($\Phi$3)} 
such that 
\begin{equation}
\label{eq:1.18}
\int_0^1 s^{-1}\Phi(s^{-1})^{-q}\,ds<\infty.
\end{equation}
Let 
\begin{equation}
\label{eq:1.19}
\mu\in L_{{\rm ul}}^{\frac{Nq}{N+2},\infty}\Phi(L)^{\frac{Nq}{N+2}},\qquad \nu\in {\mathcal M}_{{\rm ul}}.
\end{equation}
Then there exists a solution $(u,v)$ to problem~\eqref{eq:P} in ${\mathbb R}^N\times(0,T)$ 
for some $T>0$ such that 
\begin{equation*}
\begin{split}
 & \sup_{t\in(0,T)}\left\{\|u(t)\|_{\Phi,\frac{Nq}{N+2},\frac{Nq}{N+2}}+t^{\frac{N+2}{2q}}\Phi(t^{-1})\|u(t)\|_{L^\infty}\right\}<\infty,\\
 & \sup_{t\in(0,T)}\left\{\|v(t)\|_{L^1_{\rm ul}}+t^{\frac{N}{2}}\|v(t)\|_{L^\infty}\right\}<\infty.
\end{split}
\end{equation*}
Furthermore, 
$$
\lim_{t\to +0}|||u(t)-S(D_1t)\mu|||_{\Phi,\frac{Nq}{N+2},\frac{Nq}{N+2}}=0,
\quad
\lim_{t\to +0}\|v(t)-S(D_2t)\nu\|_{L^1_{{\rm ul}}}=0.
$$
\end{theorem}
\begin{theorem}
\label{Theorem:1.6} 
Consider case~{\rm (E)}. 
Let $\Phi$ be a non-decreasing function in $[0,\infty)$ with properties~{\rm ($\Phi$1)}--{\rm ($\Phi$3)} and satisfy \eqref{eq:1.18}. 
Let 
\begin{equation}
\label{eq:1.20}
\mu\in {\mathfrak L}_{{\rm ul}}^{1,\infty}\Phi({\mathfrak L}),\qquad \nu\in {\mathcal M}_{{\rm ul}}.
\end{equation}
Then there exists a solution $(u,v)$ to problem~\eqref{eq:P} in ${\mathbb R}^N\times(0,T)$ 
for some $T>0$ such that 
\begin{equation*}
\begin{split}
 & \sup_{t\in(0,T)}\left\{|||u(t)|||_{\Phi,1,1}+t^{\frac{N+2}{2q}}\Phi(t^{-1})\|u(t)\|_{L^\infty}\right\}<\infty,\\
 & \sup_{t\in(0,T)}\left\{\|v(t)\|_{L^1_{\rm ul}}+t^{\frac{N}{2}}\|v(t)\|_{L^\infty}\right\}<\infty.
\end{split}
\end{equation*}
Furthermore, 
$$
\lim_{t\to +0}|||u(t)-S(D_1t)\mu|||_{\Phi,1,1}=0,
\quad
\lim_{t\to +0}\|v(t)-S(D_2t)\nu\|_{L^1_{{\rm ul}}}=0.
$$
\end{theorem}
Similarly to Theorems~\ref{Theorem:1.1}--\ref{Theorem:1.2},  
in Section~7, 
we discuss the optimality of Theorems~\ref{Theorem:1.3}--\ref{Theorem:1.6}.
\medskip

The rest of this paper is organized as follows. 
In Section~2 we treat cases~(A) and (F), and prove Theorems~\ref{Theorem:1.1} and \ref{Theorem:1.2}.
In Section~3 we establish decay estimates of $S(t)\varphi$ in uniformly local weak type Zygmund spaces 
$L_{{\rm ul}}^{r,\infty}\Phi(L)^\alpha$ and ${\mathfrak L}_{{\rm ul}}^{r,\infty}\Phi({\mathfrak L})^\alpha$. 
In Section~4 we treat case~(B) and prove Theorem~\ref{Theorem:1.3} 
using $L_{{\rm ul}}^{r,\infty}(\log L)^\alpha$ and ${\mathfrak L}_{{\rm ul}}^{r,\infty}(\log {\mathfrak L})^\alpha$. 
In Section~5 we treat case~(C) and prove Theorem~\ref{Theorem:1.4} using ${\mathfrak L}_{{\rm ul}}^{r,\infty}(\log {\mathfrak L})^\alpha$. 
In Section~6 we treat cases~(D) and (E) and prove Theorems~\ref{Theorem:1.5} and \ref{Theorem:1.6} 
using $L_{{\rm ul}}^{r,\infty}\Phi(L)^\alpha$ and ${\mathfrak L}_{{\rm ul}}^{r,\infty}\Phi({\mathfrak L})^\alpha$. 
In Section~7, taking into the account of Proposition~\ref{Proposition:1.1}, we discuss the optimality of Theorems~\ref{Theorem:1.1}--\ref{Theorem:1.6}.  
\section{Proofs of Theorems~\ref{Theorem:1.1} and \ref{Theorem:1.2}}
This section is divided into three subsections. 
In Section~2.1 we construct approximate solutions to problem~\eqref{eq:P}. 
In Section~2.2 we introduce similar transformation of solutions to problem~\eqref{eq:P}. 
In Section~2.3 we prove Theorems~\ref{Theorem:1.1} and \ref{Theorem:1.2}. 
In all that follows we will use $C$ to denote generic positive constants and point out that $C$  
may take different values  within a calculation. 
For any positive functions $f_1$ and $f_2$ in $(0,\infty)$, we write 
$$
\mbox{$f_1\asymp f_2$ for $s>0$}\quad\mbox{if}\quad
\mbox{$Cf_2(s)\le f_1(s)\le Cf_2(s)$ for $s>0$}.
$$
\subsection{Approximate solutions}
Let $\mu$, $\nu\in{\mathcal M}$. 
Set 
$$
u_0(x,t)\coloneqq[S(D_1t)\mu](x),
\quad 
v_0(x,t)\coloneqq[S(D_2t)\nu](x),
\quad (x,t)\in{\mathbb R}^N\times(0,\infty).
$$ 
For $n=1,2,\dots$, we define the functions $u_n$ and $v_n$ in ${\mathbb R}^N\times(0,\infty)$ inductively by 
\begin{equation}
\label{eq:2.1}
\begin{split}
u_n(x,t) & \coloneqq[S(D_1t)\mu](x)+\int_0^t [S(D_1(t-s))v_{n-1}(s)^p](x)\,ds,\\
v_n(x,t) & \coloneqq[S(D_2t)\nu](x)+\int_0^t [S(D_2(t-s))u_{n-1}(s)^q](x)\,ds,
\end{split}
\end{equation}
for almost all $(x,t)\in{\mathbb R}^N\times(0,\infty)$.
By induction we see that
\begin{equation}
\label{eq:2.2}
\begin{split}
 & 0\le u_0(x,t)\le u_1(x,t)\le\cdots\le u_n(x,t)\le\cdots,\\
 & 0\le v_0(x,t)\le v_1(x,t)\le\cdots\le v_n(x,t)\le\cdots,
\end{split}
\end{equation}
for almost all $(x,t)\in{\mathbb R}^N\times(0,\infty)$.  
Then we can define the limits
\begin{equation}
\label{eq:2.3}
u(x,t)\coloneqq\lim_{n\to\infty}u_n(x,t),\quad v(x,t)\coloneqq\lim_{n\to\infty}v_n(x,t), 
\end{equation}
for almost all $(x,t)\in{\mathbb R}^N\times(0,\infty)$, and see that 
$(u,v)$ satisfies integral system~\eqref{eq:1.1} in ${\mathbb R}^N\times(0,\infty)$. 
If $u$ and $v$ are finite almost everywhere in ${\mathbb R}^N\times(0,T)$ for some $T\in(0,\infty]$, 
then $(u,v)$ is a solution to problem~\eqref{eq:P} in ${\mathbb R}^N\times(0,T)$. 

Assume that there exists a supersolution~$(\overline{u},\overline{v})$ to problem~\eqref{eq:P} in ${\mathbb R}^N\times(0,T)$ 
for some $T\in(0,\infty]$. Similarly to \eqref{eq:2.2}, by induction we see that 
\begin{align*}
 & 0\le u_0(x,t)\le u_1(x,t)\le\cdots\le u_n(x,t)\le\cdots\le \overline{u}(x,t)<\infty,\\
 & 0\le v_0(x,t)\le v_1(x,t)\le\cdots\le v_n(x,t)\le\cdots\le\overline{v}(x,t)<\infty,
\end{align*}
for almost all $(x,t)\in{\mathbb R}^N\times(0,T)$. 
Then $(u,v)$ defined by \eqref{eq:2.3} is a solution to problem~\eqref{eq:P} in ${\mathbb R}^N\times(0,T)$ such that
$$
0\le u(x,t)\le\overline{u}(x,t)<\infty,\quad 0\le v(x,t)\le\overline{v}(x,t)<\infty,
$$
for almost all $(x,t)\in{\mathbb R}^N\times(0,T)$. 
\subsection{Transformations of solutions}
Let $(u,v)$ be a solution to problem~{\rm (P)} in ${\mathbb R}^N\times(0,T)$ for some $T\in(0,\infty)$.
Let $k>0$. Set
$$
\hat{u}(x,t)\coloneqq T^{\frac{p+1}{pq-1}}u(k T^{1/2}x,Tt),
\qquad
\hat{v}(x,t)\coloneqq T^{\frac{q+1}{pq-1}}v(k T^{1/2}x,Tt),
$$
for $x\in {\mathbb R}^N$ and $t\in(0,1)$.
Then $(\hat{u},\hat{v})$ satisfies
\begin{equation*}
\left\{
\begin{array}{ll}
\partial_t \hat{u}=D_1k^{-2}\Delta\hat{u}+\hat{v}^p & \quad\mbox{in}\quad{\mathbb R}^N\times(0,1),\vspace{3pt}\\
\partial_t \hat{v}=D_2k^{-2}\Delta\hat{v}+\hat{u}^q & \quad\mbox{in}\quad{\mathbb R}^N\times(0,1),\vspace{3pt}\\
(\hat{u}(\cdot,0),\hat{v}(\cdot,0))=(\hat{\mu},\hat{\nu}) & \quad\mbox{in}\quad{\mathbb R}^N.
\end{array}
\right.
\end{equation*}
Here $\hat{\mu}$ and $\hat{\nu}$ are Radon measure in ${\mathbb R}^N$ such that
$$
\hat{\mu}(K)=k^{-N}T^{\frac{p+1}{pq-1}-\frac{N}{2}}\mu(k T^{\frac{1}{2}}K),
\quad
\hat{\nu}(K)=k^{-N}T^{\frac{q+1}{pq-1}-\frac{N}{2}}\nu(k T^{\frac{1}{2}}K),
$$
for Borel sets $K$ in ${\mathbb R}^N$.
In particular, setting 
$$
k=\max\{D_1,D_2\}^{\frac{1}{2}},
$$
we see that problem~\eqref{eq:P} is transformed to problem~\eqref{eq:P} with $\max\{D_1,D_2\}=1$. 
\subsection{Proofs of Theorems~\ref{Theorem:1.1} and \ref{Theorem:1.2}}
We recall some properties in uniformly local Morrey spaces. 
It follows from \eqref{eq:1.4} that 
\begin{equation}
\label{eq:2.4}
\begin{array}{ll}
\|f\|_{M(r,\alpha;R)}\le \|f\|_{M(r,\beta;R)}\quad & \mbox{if}\quad \alpha\le\beta,\vspace{3pt}\\
\|f^k\|_{M(r,\alpha;R)}=\|f\|^k_{M(kr,k\alpha;R)}\quad & \mbox{if}\quad k>0.
\end{array}
\end{equation}
For any $k\ge 1$, there exists $C>0$ such that 
\begin{equation}
\label{eq:2.5}
\|f\|_{M(r,\alpha;kR)}\le C\|f\|_{M(r,\alpha;R)},\quad R\in(0,\infty]
\end{equation}
(see e.g., \cite{IS}*{Lemma~2.1}).
Furthermore, we have:
\begin{Lemma}
\label{Lemma:2.1}
{\rm (1)}
	Let $1\le r_1\le r_2\le\infty$ and $\alpha\in[1,r_2/r_1]$. 
  	Then there exists $C_1>0$ such that 
  	\begin{equation}
  	\label{eq:2.6}
  		\sup_{t\in(0,R^2)}\left\{t^{\frac{N}{2}\left(\frac{1}{r_1}-\frac{1}{r_2}\right)}
  	\|S(t)\varphi\|_{M(r_2,\alpha;R)}\right\}\le C_2\|\varphi\|_{M(r_1,1;R)},\quad \varphi\in M(r_1,1;R),
  	\end{equation}
  	for $R\in(0,\infty]$.
\vspace{3pt}
\newline
{\rm (2)}
Let $1\le r\le\infty$ and $\alpha\in[1,r]$. 
Then there exists $C_2>0$ such that 
	\begin{equation}
  	\label{eq:2.7}
  		\sup_{t\in(0,1)}\left\{t^{\frac{N}{2}\left(1-\frac{1}{r}\right)}
  		\|S(t)\mu\|_{M(r,\alpha)}\right\}\le C\|\mu\|_{{\mathcal M}_{{\rm ul}}},\quad \mu\in {\mathcal M}_{{\rm ul}}.
  	\end{equation}
\end{Lemma}
{\bf Proof.}
We prove Lemma~\ref{Lemma:2.1}~(1). 
The proof is divided into two steps.
\newline
{\bf Step.1}
We prove inequality~\eqref{eq:2.6} with $R=\infty$ using the following decay estimate. 
\begin{itemize}
  \item 
  For any $1\le r\le q\le \infty$, there exists $C>0$ such that 
  \begin{equation}
  \label{eq:2.8}
  \sup_{x\in{\mathbb R}^N}\|S(t)\varphi\|_{L^q(B(x,R))}\le Ct^{-\frac{N}{2}\left(\frac{1}{r}-\frac{1}{q}\right)}\sup_{x\in{\mathbb R}^N}\|\varphi\|_{L^{r}(B(x,R))}
  \end{equation}
  for $t\in(0,R^2)$ and $R>0$. (See \cite{MT}*{Corollary 3.1}.)
\end{itemize}
Let $1\le r_1\le r_2\le\infty$ and $1\le\alpha\le r_2/r_1$. 
By \eqref{eq:2.8} with $r=1$, $q=\infty$, and $R=t^{1/2}$ 
we have
\begin{equation}
\label{eq:2.9}
\begin{split}
 & |B(z,\sigma)|^{\frac{1}{r_2}-\frac{1}{\alpha}}
\|S(t)\varphi\|_{L^\alpha(B(z,\sigma))}\\
 & \le |B(z,\sigma)|^{\frac{1}{r_2}}\|S(t)\varphi\|_\infty
\le Ct^{\frac{N}{2r_2}}\|S(t)\varphi\|_\infty\\
 & \le Ct^{\frac{N}{2r_2}}\cdot Ct^{-\frac{N}{2}}
\sup_{x\in{\mathbb R}^N}\|\varphi\|_{L^1(B(x,t^{1/2}))}\\
 & \le Ct^{\frac{N}{2r_2}}\cdot Ct^{-\frac{N}{2r_1}}
\sup_{x\in{\mathbb R}^N}\left\{|B(x,t^{1/2})|^{\frac{1}{r_1}-1}\|\varphi\|_{L^1(B(x,t^{1/2}))}\right\}\\
 & \le Ct^{-\frac{N}{2}\left(\frac{1}{r_1}-\frac{1}{r_2}\right)}\|\varphi\|_{M(r_1,1;\infty)}
\end{split}
\end{equation}
for $z\in{\mathbb R}^N$, $t>0$, and $\sigma\in(0,t^{1/2})$. 
Furthermore, 
by \eqref{eq:2.8} with $r=q=1$ and $R=\sigma$ and with $r=1$, $q=\infty$, and $R=t^{1/2}$
we have
\begin{equation}
\label{eq:2.10}
\begin{split}
 & |B(z,\sigma)|^{\frac{1}{r_2}-\frac{1}{\alpha}}\|S(t)\varphi\|_{L^\alpha(B(z,\sigma))}\\
 & \le \left(|B(z,\sigma)|^{\frac{\alpha}{r_2}-1}\|S(t)\varphi\|_{L^1(B(z,\sigma))}\right)^{\frac{1}{\alpha}}\|S(t)\varphi\|_\infty^{1-\frac{1}{\alpha}}\\
 & \le \left(C|B(z,\sigma)|^{\frac{\alpha}{r_2}-1}\sup_{x\in{\mathbb R}^N}\|\varphi\|_{L^1(B(x,\sigma))}\right)^{\frac{1}{\alpha}}
 \left(Ct^{-\frac{N}{2}}\sup_{x\in{\mathbb R}^N}\|\varphi\|_{L^1(B(x,t^{1/2}))}\right)^{1-\frac{1}{\alpha}}\\
 & \le \left(C|B(z,\sigma)|^{\frac{\alpha}{r_2}-\frac{1}{r_1}}\sup_{x\in{\mathbb R}^N}\left\{|B(x,\sigma)|^{\frac{1}{r_1}-1}\|\varphi\|_{L^1(B(x,\sigma))}\right\}\right)^{\frac{1}{\alpha}}\\
 & \qquad\qquad\times\left(Ct^{-\frac{N}{2r_1}}\sup_{x\in{\mathbb R}^N}\left\{|B(x,t^{1/2})|^{\frac{1}{r_1}-1}\|\varphi\|_{L^1(B(x,t^{1/2}))}\right\}\right)^{1-\frac{1}{\alpha}}\\
  & \le C\left(|B(z,t^{1/2})|^{\frac{\alpha}{r_2}-\frac{1}{r_1}}\|\varphi\|_{M(r_1,1;\infty)}\right)^{\frac{1}{\alpha}}
 \left(Ct^{-\frac{N}{2r_1}}\|\varphi\|_{M(r_1,1;\infty)}\right)^{1-\frac{1}{\alpha}}\\
  & \le Ct^{-\frac{N}{2}\left(\frac{1}{r_1}-\frac{1}{r_2}\right)}\|\varphi\|_{M(r_1,1;\infty)}
\end{split}
\end{equation}
for $z\in{\mathbb R}^N$, $t>0$, and $\sigma\in(t^{1/2},\infty)$. Here we used the relation $\alpha/r_2\le 1/r_1$. 
Combining \eqref{eq:2.9} and \eqref{eq:2.10}, we obtain 
$$
|B(z,\sigma)|^{\frac{1}{r_2}-\frac{1}{\alpha}}\|S(t)\varphi\|_{L^\alpha(B(z,\sigma))}
\le Ct^{-\frac{N}{2}\left(\frac{1}{r_1}-\frac{1}{r_2}\right)}\|\varphi\|_{M(r_1,1;\infty)}
$$
for $z\in{\mathbb R}^N$ and $\sigma\in(0,\infty)$. 
This implies \eqref{eq:2.6} with $R=\infty$.
(See also \cite{S2}*{Proposition~4.1} for another proof of \eqref{eq:2.6} with $R=\infty$.) 
\vspace{3pt}
\newline
{\bf Step 2.}
Let $1\le r_1\le r_2\le\infty$ and $\alpha\in[1,r_2/r_1]$.
Let $R\in(0,\infty]$. 
For the proof of \eqref{eq:2.6} with $R<\infty$,  
it suffices to find $C>0$ such that 
\begin{equation}
\label{eq:2.11}
t^{\frac{N}{2}\left(\frac{1}{r_1}-\frac{1}{r_2}\right)}
\left\|\chi_{B(z,R)}S(t)\varphi\right\|_{M(r_2,\alpha;\infty)}\le C|||\varphi|||_{M(r_1,1;R)}
\end{equation}
for $z\in{\mathbb R}^n$ and $0<t\le R^2$. 
Then, by translating if necessary, we have only to consider the case of $z=0$. 

The proof is a modification of the proofs of \cite{HI01}*{Theorem~1.2} and \cite{IIK}*{Proposition~3.2}. 
By Besicovitch's covering lemma 
we can find an integer $m$ depending only on $n$ and 
a set $\{x_{k,i}\}_{k=1,\dots,m,\,i\in{\mathbb N}}\subset{\mathbb R}^n\setminus B(0,10R)$ such that 
\begin{equation}
\label{eq:2.12}
B_{k,i}\cap B_{k,j}=\emptyset\quad\mbox{if $i\not=j$}
\qquad\mbox{and}\qquad
{\mathbb R}^n\setminus B(0,10R)\subset\bigcup_{k=1}^m\bigcup_{i=1}^\infty B_{k,i},
\end{equation}
where $B_{k,i}:=\overline{B(x_{k,i},R)}$. 
Then 
\begin{equation}
\label{eq:2.13}
\left|\left[S(t)\varphi\right](x)\right|\le |u_0(x,t)|
+\sum_{k=1}^m\sum_{i=1}^\infty |u_{k,i}(x,t)|,
\quad(x,t)\in{\mathbb R}^n\times(0,R^2),
\end{equation}
where 
$$
u_0(x,t):=[S(t)(\varphi\chi_{B(0,10R)})](x),
\quad
u_{k,i}(x,t):=[S(t)(\varphi\chi_{B_{k,i}})](x).
$$
By \eqref{eq:2.5} and \eqref{eq:2.11} with $R=\infty$ we have 
\begin{equation}
\label{eq:2.14}
\begin{split}
 & \|\chi_{B(0,R)}u_0(t)\|_{M(r_2,\alpha,\infty)}\le\|u_0(t)\|_{M(r_2,\alpha;\infty)}\\
 & \qquad\quad
 \le Ct^{-\frac{N}{2}\left(\frac{1}{r_1}-\frac{1}{r_2}\right)}\left\|\varphi\chi_{B(0,10R)}\right\|_{M(r_1,1;\infty)}\\
 & \qquad\quad
 \le Ct^{-\frac{N}{2}\left(\frac{1}{r_1}-\frac{1}{r_2}\right)}\|\varphi\|_{M(r_1,1;10R)}\\
 & \qquad\quad
 \le Ct^{-\frac{N}{2}\left(\frac{1}{r_1}-\frac{1}{r_2}\right)}\|\varphi\|_{M(r_1,1;R)},\quad t\in(0,R^2].
\end{split}
\end{equation}
Let $k=1,\dots,m$ and $i{\in \mathbb{N}}$. 
Then we see that
\begin{equation}
\label{eq:2.15}
\begin{split}
|u_{k,i}(x,t)| & \le C\int_{B(x_{k,i},R)}G(x-y,t)|\varphi(y)|\,dy\\
 & =C\int_{\mathbb R^n}G(x-z-x_{k,i},t)\varphi_{k,i}(z)\,dz
\end{split}
\end{equation}
for $(x,t)\in{\mathbb R}^n\times(0,\infty)$, where 
$\varphi_{k,i}(x)=|\varphi(x+x_{k,i})|\chi_{B(0,R)}$. 
It follows from $|x_{k,i}|\ge 10R$  that 
$$
\frac{|x-z-x_{k,i}|}{t^{1/2}}\ge
\frac{|x_{k,i}|-|x-z|}{t^{1/2}}\ge\frac{|x_{k,i}|}{2t^{1/2}}+\frac{5R-2|x-z|}{t^{1/2}}+\frac{|x-z|}{t^{1/2}}\ge\frac{|x_{k,i}|}{2R}+\frac{|x-z|}{t^{1/2}}
$$
for $x$, $z\in B(0,R)$ and $t\in(0,R^2)$. 
This implies that 
\begin{equation}
\label{eq:2.16}
G(x-z-x_{k,i},t) 
\le (4\pi t)^{-\frac{N}{2}}\exp\left(-\frac{|x_{k,i}|^2}{16R^2}-\frac{|x-z|^2}{4t}\right) \le \exp\left(-\frac{|x_{k,i}|^2}{16R^2}\right) G(x-z,t)
\end{equation}
for $x$, $z\in B(0,R)$ and $t\in(0,R^2)$. 
We observe from \eqref{eq:2.15} and \eqref{eq:2.16} that 
$$
|u_{k,i}(x,t)|\le C\exp\left(-\frac{|x_{k,i}|^2}{16R^2}\right)[S(t)\varphi_{k,i}](x)
$$
for $x\in B(0,R)$ and $t\in(0,R^2)$. 
Then, by \eqref{eq:2.11} with $R=\infty$  we obtain 
\begin{equation}
\label{eq:2.17}
\begin{split}
 & \|u_{k,i}(t)\chi_{B(0,R)}\|_{M(r_2,\alpha;\infty)}\\
 & \le C\exp\left(-\frac{|x_{k,i}|^2}{16R^2}\right)\|S(t)\varphi_{k,i}\|_{M(r_2,\alpha;\infty)}\\
 & \le C\exp\left(-\frac{|x_{k,i}|^2}{16R^2}\right) t^{-\frac{N}{2}\left(\frac{1}{r_1}-\frac{1}{r_2}\right)}\|\varphi_{k,i}\|_{M(r_1,1;\infty)}\\
 & \le C\exp\left(-\frac{|x_{k,i}|^2}{16R^2}\right) t^{-\frac{N}{2}\left(\frac{1}{r_1}-\frac{1}{r_2}\right)}\|\varphi\chi_{B(x_{k,i},R)}\|_{M(r_1,1;\infty)}\\
 & \le C\exp\left(-\frac{|x_{k,i}|^2}{16R^2}\right) t^{-\frac{N}{2}\left(\frac{1}{r_1}-\frac{1}{r_2}\right)}\|\varphi\|_{M(r_1,1;R)}
\end{split}
\end{equation}
for $t\in(0,R^2)$. 

On the other hand, since
$$
\frac{|y|}{2}\le\frac{1}{2}\left(|x_{k,i}|+R\right)\le|x_{k,i}|\quad\mbox{for}\quad y\in B_{k,i},
$$
we have 
$$
\frac{1}{|B_{k,i}|}\int_{B_{k,i}} \exp\left(-\frac{|y|^2}{64R^2}\right) \,dy\ge \exp\left(-\frac{|x_{k,i}|^2}{16R^2}\right).
$$
Then, by \eqref{eq:2.12} we see that 
\begin{equation}
\label{eq:2.18}
\begin{split}
\sum_{i=1}^\infty \exp\left(-\frac{|x_{k,i}|^2}{16R^2}\right) 
 & \le CR^{-N}\sum_{i=1}^\infty\int_{B_{k,i}} \exp\left(-\frac{|y|^2}{64R^2}\right) \,dy\\
 & \le CR^{-N}\int_{{\mathbb R}^n}  \exp\left(-\frac{|y|^2}{64R^2}\right) \,dy\,dy
\le C
\end{split}
\end{equation}
for $R>0$. 
Combining \eqref{eq:2.13}, \eqref{eq:2.14}, \eqref{eq:2.17}, and \eqref{eq:2.18}
we obtain 
\begin{align*}
 & t^{\frac{N}{2}\left(\frac{1}{r_1}-\frac{1}{r_2}\right)}
 \left\|\chi_{B(0,R)}S(t)\varphi\right\|_{M(r_2,\alpha;\infty)}\\
 & \le C\|\varphi\|_{M(r_1,1;R)}+C\|\varphi\|_{M(r_1,1;R)}
 \sum_{k=1}^m\sum_{i=1}^\infty \exp\left(-\frac{|x_{k,i}|^2}{16R^2}\right) 
 \le C\|\varphi\|_{M(r_1,1;R)}
\end{align*}
for $t\in(0,R^2)$. 
This implies \eqref{eq:2.11} with $z=0$. 
Thus \eqref{eq:2.6} with $R<\infty$ holds, and 
Lemma~\ref{Lemma:2.1}~(1) follows. 
Similarly, we obtain Lemma~\ref{Lemma:2.1}~(2), and the proof of Lemma~\ref{Lemma:2.1} is complete. 
$\Box$\vspace{3pt}

We prove Theorems~\ref{Theorem:1.1} and \ref{Theorem:1.2}. 
In cases~(A) and (F), 
following the arguments in \cite{FI02}*{Section~3} and \cite{IKS}, 
we construct a supersolution to problem~\eqref{eq:P} to find a solution to problem~\eqref{eq:P}. 
\vspace{3pt}
\newline
{\bf Proof of Theorem~\ref{Theorem:1.1}.}
Consider case~(A), that is, 
$$
\frac{q+1}{pq-1}<\frac{N}{2}.
$$
Let $D=\min\{D_1,D_2\}$ and $D'\coloneqq\max\{D_1,D_2\}$. 
By Section~2.2 
it suffices to consider the case of $D'=1$. 
Then 
\begin{equation}
\label{eq:2.19}
G(x,D_it)=(4\pi D_it)^{-\frac{N}{2}}\exp\left(-\frac{|x|^2}{4D_it}\right)
\le D^{-\frac{N}{2}}G(x,t)
\end{equation}
in ${\mathbb R}^N\times(0,\infty)$, where $i=1,2$. 

Let $\delta_A>0$ be small enough and assume \eqref{eq:1.6}.
Set
\begin{equation}
\label{eq:2.20}
\begin{aligned}
&
w(x,t)\coloneqq\left[S(t)\left(\mu^{\alpha_A}+\nu^{\beta_A}\right)\right](x),
\\
&
\overline{u}(x,t)\coloneqq2D^{-\frac{N}{2}}w(x,t)^{\frac{1}{{\alpha_A}}},
\,\,\,
\overline{v}(x,t)\coloneqq2D^{-\frac{N}{2}}w(x,t)^{\frac{1}{{\beta_A}}},
\end{aligned}
\end{equation}
for $(x,t)\in{\mathbb R}^N\times(0,\infty)$. 
It follows from the semigroup property of $S(t)$ that 
\begin{equation}
\label{eq:2.21}
	w(x,t)=[S(t-s)w(s)](x),\quad x\in{\mathbb R}^N,\,0\le s\le t.
\end{equation}
Since 
\begin{equation}
\label{eq:2.22}
\alpha_A\beta_A^{-1}r_2^*=\frac{q+1}{p+1}\frac{N}{2}\frac{pq-1}{q+1}=\frac{N}{2}\frac{pq-1}{p+1}=r_1^*,
\end{equation}
it follows from \eqref{eq:1.6} and \eqref{eq:2.4} that 
$$
\left\|\mu^{\alpha_A}+\nu^{\beta_A}\right\|_{M(\beta_A^{-1}r_2^*,1;T^{\frac{1}{2}})}
\le\|\mu\|_{M(r_1^*,\alpha_A;T^{\frac{1}{2}})}^{\alpha_A}+\|\nu\|_{M(r_2^*,{\beta_A};T^{\frac{1}{2}})}^{\beta_A}\le\delta_A. 
$$
This together with \eqref{eq:2.6} implies that 
\begin{equation}
\label{eq:2.23}
 \|w(t)\|_{M(r,\eta;T^{\frac{1}{2}})}\le C\delta_A t^{-\frac{N}{2}\left(\frac{{\beta_A}}{r_2^*}-\frac{1}{r}\right)},
 \quad t\in(0,T), 
\end{equation}
for $r\in[\beta_A^{-1}r_2^*,\infty]$ and $\eta\in[1,\beta_A r/r_2^*]$.

We prove that $(\overline{u},\overline{v})$ is a supersolution to problem~\eqref{eq:P} in ${\mathbb R}^N\times(0,T)$. 
It follows from~\eqref{eq:1.5} that
\begin{align*}
\frac{q}{\alpha_A}=\frac{q(p+1)}{\beta_A(q+1)}>1,\quad
 & -\frac{N\beta_A}{2r_2^*}\left(\frac{q}{\alpha_A}-1\right)+1=\frac{N\beta_A}{2r_2^*}\left(-\frac{q}{\alpha_A}+1+\frac{2r_2^*}{N\beta_A}\right)\\
 & =\frac{N\beta_A}{2r_2^*}\left(-\frac{1}{\beta_A}\frac{pq+q}{q+1}+1+\frac{1}{\beta_A}\frac{pq-1}{q+1}\right)
 =\frac{N\beta_A}{2r_2^*}\left(1-\frac{1}{\beta_A}\right)>0.
\end{align*}
These together with \eqref{eq:2.19}, \eqref{eq:2.21}, and \eqref{eq:2.23} with $r=\infty$ imply that 
\begin{equation}
\label{eq:2.24}
\begin{split}
 & \int_0^t [S(D_2(t-s))\overline{u}(s)^q](x)\,ds\\
 & \le D^{-\frac{N}{2}}\int_0^t [S(t-s)\overline{u}(s)^q](x)\,ds\le C\int_0^t \left[S(t-s)w(s)^{\frac{q}{\alpha_A}}\right](x)\,ds\\
 & \le C\int_0^t \left[S(t-s)\|w(s)\|_{L^\infty}^{\frac{q}{\alpha_A}-1}w(s)\right](x)\,ds 
 =Cw(x,t)\int_0^t \|w(s)\|_{L^\infty}^{\frac{q}{\alpha_A}-1}\,ds\\
 & \le C\delta_A^{\frac{q}{\alpha_A}-1} w(x,t)\int_0^t s^{-\frac{N\beta_A}{2r_2^*}\left(\frac{q}{\alpha_A}-1\right)}\,ds
 =C\delta_A^{\frac{q}{\alpha_A}-1} w(x,t)\int_0^t s^{-1+\frac{N\beta_A}{2r_2^*}\left(1-\frac{1}{\beta_A}\right)}\,ds\\
 & \le C\delta_A^{\frac{q}{\alpha_A}-1} t^{\frac{N\beta_A}{2r_2^*}\left(1-\frac{1}{\beta_A}\right)}w(x,t)
 \quad\mbox{in ${\mathbb R}^N\times(0,T)$}.
\end{split}
\end{equation}
Taking small enough $\delta_A>0$ if necessary, 
by Jensen's inequality, \eqref{eq:2.19}, \eqref{eq:2.23}, and \eqref{eq:2.24} we obtain 
\begin{equation}
\label{eq:2.25}
\begin{split}
 & [S(D_2t)\nu](x)+\int_0^t [S(D_2(t-s))\overline{u}(s)^q](x)\,ds\\
 & \le D^{-\frac{N}{2}}[S(t)\nu](x)+C\delta_A^{\frac{q}{\alpha_A}-1} t^{\frac{N\beta_A}{2r_2^*}\left(1-\frac{1}{\beta_A}\right)}w(x,t)\\
 & \le D^{-\frac{N}{2}}\left[S(t)\nu^{\beta_A}\right](x)^{\frac{1}{{\beta_A}}}+C\delta_A^{\frac{q}{\alpha_A}-1} t^{\frac{N\beta_A}{2r_2^*}\left(1-\frac{1}{\beta_A}\right)}\|w(t)\|_{L^\infty}^{1-\frac{1}{{\beta_A}}}w(x,t)^{\frac{1}{{\beta_A}}}\\
 & \le \frac{1}{2}\overline{v}(x,t)+C\delta_A^{\frac{q}{\alpha_A}-\frac{1}{\beta_A}}w(t)^{\frac{1}{{\beta_A}}}
 \le \overline{v}(x,t)\quad\mbox{in}\quad {\mathbb R}^N\times(0,T).
\end{split}
\end{equation}
Here we used the relation 
\begin{equation}
  \notag 
  \frac{q}{\alpha_A}-\frac{1}{\beta_A}
  =
  \frac{1}{\alpha_A}\left(
    q-\frac{q+1}{p+1}
  \right)
  =
  \frac{1}{\alpha_A}\frac{pq-1}{p+1}
  > 0. 
\end{equation}

On the other hand, 
it follows from \eqref{eq:1.5} that 
\begin{equation}
\begin{split}
\label{eq:2.26}
 & -\frac{N\beta_A}{2r_2^*}\left(\frac{p}{\beta_A}-1\right)+1
=\frac{N\beta_A}{2r_2^*}\left(-\frac{p}{\beta_A}+1+\frac{2r_2^*}{N\beta_A}\right)\\
 & =\frac{N\beta_A}{2r_2^*}\left(-\frac{1}{\beta_A}\frac{pq+p}{q+1}+1+\frac{1}{\beta_A}\frac{pq-1}{q+1}\right)
=\frac{N\beta_A}{2r_2^*}\left(1-\frac{1}{\beta_A}\frac{p+1}{q+1}\right)\\
 & =\frac{N\beta_A}{2r_2^*}\left(1-\frac{1}{\alpha_A}\right)>0.
\end{split}
\end{equation}
Then, similarly to \eqref{eq:2.24}, in the case of $p>{\beta_A}$, 
we have 
\begin{equation}
\label{eq:2.27}
\begin{split}
 & \int_0^t [S(D_1(t-s))\overline{v}(s)^p](x)\,ds\\
 & \le D^{-\frac{N}{2}}\int_0^t [S(t-s)\overline{v}(s)^p](x)\,ds\le C\int_0^t [S(t-s)w(s)^{\frac{p}{{\beta_A}}}](x)\,ds\\
 & \le Cw(x,t)\int_0^t \|w(s)\|_{L^\infty}^{\frac{p}{{\beta_A}}-1}\,ds
 \le C\delta_A^{\frac{p}{{\beta_A}}-1} w(x,t)\int_0^t s^{-1+\frac{N\beta_A}{2r_2^*}\left(1-\frac{1}{\alpha_A}\right)}\,ds\\
 & \le C\delta_A^{\frac{p}{{\beta_A}}-1} t^{\frac{N\beta_A}{2r_2^*}\left(1-\frac{1}{\alpha_A}\right)}w(x,t)
 \quad\mbox{in}\quad {\mathbb R}^N\times(0,T).
\end{split}
\end{equation}
Taking small enough $\delta_A>0$ if necessary, 
by Jensen's inequality, \eqref{eq:2.19}, \eqref{eq:2.23}, and \eqref{eq:2.27} 
we obtain 
\begin{equation}
\label{eq:2.28}
\begin{split}
 & [S(D_1t)\mu](x)+\int_0^t [S(D_1(t-s))\overline{v}(s)^p](x)\,ds\\
 & \le D^{-\frac{N}{2}}[S(t)\mu](x)+C\delta_A^{\frac{p}{{\beta_A}}-1} t^{\frac{N\beta_A}{2r_2^*}\left(1-\frac{1}{\alpha_A}\right)}w(x,t)\\
 & \le D^{-\frac{N}{2}}\left[S(t)\mu^{\alpha_A}\right](x)^{\frac{1}{\alpha_A}}
 +C\delta_A^{\frac{p}{{\beta_A}}-1} t^{\frac{N\beta_A}{2r_2^*}\left(1-\frac{1}{\alpha_A}\right)}\|w(t)\|_{L^\infty}^{1-\frac{1}{\alpha_A}}w(x,t)^{\frac{1}{{\alpha_A}}}\\
 & \le \frac{1}{2}\overline{u}(x,t)+C\delta_A^{\frac{p}{{\beta_A}}-\frac{1}{\alpha_A}}w(x,t)^{\frac{1}{\alpha_A}}
 \le \overline{u}(x,t)
\end{split}
\end{equation}
in ${\mathbb R}^N\times(0,T)$ in the case of $p>{\beta_A}$. 
Here we used the relation
\begin{equation}
\label{eq:2.29}
\frac{p}{{\beta_A}}-\frac{1}{\alpha_A}=\frac{1}{\beta_A}\left(p-\frac{p+1}{q+1}\right)=\frac{1}{\beta_A}\frac{pq-1}{q+1}>0.
\end{equation}
In the case of $p\le{\beta_A}$, 
it follows from Jensen's inequality, \eqref{eq:2.19}, and \eqref{eq:2.21} implies that 
\begin{equation}
\label{eq:2.30}
\begin{split}
 & \int_0^t [S(D_1(t-s))\overline{v}(s)^p](x)\,ds\\
 & \le D^{-\frac{N}{2}}\int_0^t [S(t-s)\overline{v}(s)^p](x)\,ds\le C\int_0^t [S(t-s)w(s)^{\frac{p}{{\beta_A}}}](x)\,ds\\
 & \le C\int_0^t [S(t-s)w(s)](x)^{\frac{p}{{\beta_A}}}\,ds=Ctw(x,t)^{\frac{p}{{\beta_A}}}
 \quad\mbox{in}\quad {\mathbb R}^N\times(0,T).
\end{split}
\end{equation}
Then, similarly to \eqref{eq:2.28}, 
by \eqref{eq:2.26} and \eqref{eq:2.29}, 
taking small enough $\delta_A>0$ if necessary, 
we obtain  
\begin{equation}
\label{eq:2.31}
\begin{split}
 & [S(D_1t)\mu](x)+\int_0^t [S(D_1(t-s))\overline{v}(s)^p](x)\,ds\\
 & \le D^{-\frac{N}{2}}[S(t)\mu](x)+Ctw(x,t)^{\frac{p}{{\beta_A}}-\frac{1}{\alpha_A}}w(x,t)^{\frac{1}{\alpha_A}}\\
 & \le D^{-\frac{N}{2}}[S(t)\mu^{\alpha_A}](x)^{\frac{1}{\alpha_A}}
 +C\delta_A^{\frac{p}{{\beta_A}}-\frac{1}{\alpha_A}} t^{1-\frac{N\beta_A}{2r_2^*}\left(\frac{p}{\beta_A}-\frac{1}{\alpha_A}\right)}w(x,t)^{\frac{1}{\alpha_A}}\\
 & \le \frac{1}{2}\overline{u}(x,t)+C\delta_A^{\frac{p}{{\beta_A}}-\frac{1}{\alpha_A}}w(x,t)^{\frac{1}{\alpha_A}}
 \le \overline{u}(x,t)
\end{split}
\end{equation}
in ${\mathbb R}^N\times(0,T)$ in the case of $p\le{\beta_A}$. 
Therefore we deduce from \eqref{eq:2.25}, \eqref{eq:2.28}, and \eqref{eq:2.31} that 
$(\overline{u},\overline{v})$ is a supersolution to problem~\eqref{eq:P} in ${\mathbb R}^N\times(0,T)$. 
Then, by the arguments in Section~2.1 
we find a solution $(u,v)$ to problem~\eqref{eq:P} in ${\mathbb R}^N\times(0,T)$ such that 
\begin{equation}
\label{eq:2.32}
0\le u(x,t)\le\overline{u}(x,t),\quad 0\le v(x,t)\le\overline{v}(x,t),\quad
(x,t)\in{\mathbb R}^N\times(0,T). 
\end{equation}
These together with \eqref{eq:2.4}, \eqref{eq:2.20}, \eqref{eq:2.22}, and \eqref{eq:2.23} imply that 
\begin{align*}
& \|u(t)\|_{M(r_1^*,\alpha_A; T^{\frac{1}{2}})}^{\alpha_A}+\|v(t)\|_{M(r_2^*, \beta_A; T^{\frac{1}{2}})}^{\beta_A}\\
& \le C\|w(t)^{\frac{1}{\alpha_A}}\|_{M(r_1^*,\alpha_A;T^{\frac{1}{2}})}^{\alpha_A}
+C\|w(t)^{\frac{1}{{\beta_A}}}\|_{M(r_2^*,{\beta_A};T^{\frac{1}{2}})}^{\beta_A}
\le C\|w(t)\|_{M(\beta_A^{-1}r_2^*,1;T^{\frac{1}{2}})}\le C,
\\
& t^{\frac{N}{2r_1^*}}\|u(t)\|_{L^\infty}+t^{\frac{N}{2r_2^*}}\|v(t)\|_{L^\infty}\le 
t^{\frac{N}{2r_1^*}}\|w(t)\|_{L^\infty}^{\frac{1}{\alpha_A}}+t^{\frac{N}{2r_2^*}}\|w(t)\|_{L^\infty}^{\frac{1}{{\beta_A}}}\le C,
\end{align*}
for $t\in(0,T)$. Thus \eqref{eq:1.7} and \eqref{eq:1.8} hold. 

It remains to prove \eqref{eq:1.9} for $r_1\in[1,r_1^*)$, $r_2\in[1,r_2^*)$, 
$\ell_1\in[1,\alpha_Ar_1/r_1^*]$, and $\ell_2\in[1,\beta_Ar_2/r_2^*]$. 
Since 
$$
r_1^*\ge r_2^*,\quad
pr_1^*=\frac{N}{2}\frac{p(pq-1)}{p+1}>\frac{N}{2}\frac{p(pq-1)}{p+pq}=r_2^*,
$$
and $\|f\|_{M(m_1,\ell)}\le C\|f\|_{M(m_2,\ell)}$ for $f\in M(m_2,\ell)$ if $1\le m_1\le m_2<\infty$,
it suffices to consider the case of 
\begin{equation}
\label{eq:2.33}
1\le r_2<r_1<r_1^*,\quad
r_2^*<pr_1,
\quad
\beta_A^{-1}r_2^*<r_2<r_2^*.
\end{equation} 
By \eqref{eq:2.5}, \eqref{eq:2.23}, \eqref{eq:2.24}, and \eqref{eq:2.32} we have
\begin{equation}
\label{eq:2.34}
\begin{split}
 & \|v(t)-S(D_2t)\nu\|_{M(r_2,\ell_2)}
 \le C\left\|\,\int_0^t [S(D_2(t-s))\overline{u}(s)^q](x)\,ds\,\right\|_{M(r_2,\ell_2;T^{\frac{1}{2}})} \\
 & \le Ct^{\frac{N\beta_A}{2r_2^*}\left(1-\frac{1}{\beta_A}\right)}\|w(t)\|_{M(r_2,\ell_2;T^{\frac{1}{2}})}\\
 & \le Ct^{\frac{N\beta_A}{2r_2^*}\left(1-\frac{1}{\beta_A}\right)}\cdot Ct^{-\frac{N}{2}\left(\frac{\beta_A}{r_2^*}-\frac{1}{r_2}\right)}
\le Ct^{\frac{N}{2}\left(\frac{1}{r_2}-\frac{1}{r_2^*}\right)}\to 0
\end{split}
\end{equation}
as $t\to +0$. 
Furthermore, 
since $r_1>r_2>\beta_A^{-1}r_2^*$ (see \eqref{eq:2.33}), 
if $p>{\beta_A}$, 
then, by \eqref{eq:2.5}, \eqref{eq:2.22}, \eqref{eq:2.23}, \eqref{eq:2.27}, and \eqref{eq:2.32} we obtain 
\begin{equation}
\label{eq:2.35}
\begin{split}
 & \|u(t)-S(D_1t)\mu\|_{M(r_1,\ell_1)}
\le C\left\|\,\int_0^t [S(D_1(t-s))\overline{v}(s)^p](x)\,ds\,\right\|_{M(r_1,\ell_1;T^{\frac{1}{2}})}\\
 & \le Ct^{\frac{N\beta_A}{2r_2^*}\left(1-\frac{1}{\alpha_A}\right)}\|w(t)\|_{M(r_1,\ell_1; T^{\frac{1}{2}})}\\
 & \le Ct^{\frac{N\beta_A}{2r_2^*}\left(1-\frac{1}{\alpha_A}\right)}\cdot Ct^{-\frac{N}{2}\left(\frac{{\beta_A}}{r_2^*}-\frac{1}{r_1}\right)}
\le Ct^{\frac{N}{2}\left(\frac{1}{r_1}-\frac{1}{r_1^*}\right)}\to 0
\end{split}
\end{equation}
as $t\to +0$. 
If $p\le \beta_A$, 
by \eqref{eq:2.4}, \eqref{eq:2.5}, \eqref{eq:2.23}, \eqref{eq:2.30}, \eqref{eq:2.32}, and \eqref{eq:2.33} 
we have 
\begin{equation}
\label{eq:2.36}
\begin{split}
 & \|u(t)-S(D_1t)\mu\|_{M(r_1,\ell_1)}
\le C\left\|\,\int_0^t [S(D_1(t-s))\overline{v}(s)^p](x)\,ds\,\right\|_{M(r_1,\ell_1;T^{\frac{1}{2}})}\\
 & \le Ct\|w(t)^{\frac{p}{{\beta_A}}}\|_{M(r_1,\ell_1;T^{\frac{1}{2}})}
 \le Ct\|w(t)\|_{M(\beta_A^{-1}pr_1,\beta_A^{-1}p\ell_1 T^{\frac{1}{2}})}^{\beta_A^{-1}p}\\
 & \le Ct^{1-\frac{N}{2}\left(\frac{\beta_A}{r_2^*}-\frac{\beta_A}{pr_1}\right)\frac{p}{\beta_A}}
 =Ct^{\frac{N}{2}\left(\frac{1}{r_1}+\frac{2}{N}-\frac{2}{N}\frac{p(q+1)}{pq-1}\right)}
 =Ct^{\frac{N}{2}\left(\frac{1}{r_1}-\frac{1}{r_1^*}\right)}\to 0
\end{split}
\end{equation}
as $t\to +0$. 
By \eqref{eq:2.34}, \eqref{eq:2.35}, and \eqref{eq:2.36} 
we obtain \eqref{eq:1.9}. 
Thus Theorem~\ref{Theorem:1.1} follows.
$\Box$ \vspace{5pt}

\noindent
{\bf Proof of Theorem~\ref{Theorem:1.2}.}
Consider case~(F), that is, 
$$
\frac{q+1}{pq-1}>\frac{N}{2}\quad\mbox{and}\quad q<1+\frac{2}{N}.
$$
 Assume $\mu$, $\nu\in{\mathcal M}_{{\rm ul}}$. 
Let $D\coloneqq\min\{D_1, D_2\}$ and $D'\coloneqq\max\{D_1, D_2\}$. 
Similarly to the proof of Theorem~\ref{Theorem:1.1}, 
we can assume, without loss of generality, that $D'=1$. 

Set 
$$
w(x,t)\coloneqq2D^{-\frac{N}{2}}[S(t)(\mu+\nu)](x)+2t,
\quad
(x,t)\in{\mathbb R}^N\times(0,\infty).
$$
It follows that 
\begin{equation}
\label{eq:2.37}
[S(t-s)w(s)](x)=2D^{-\frac{N}{2}}[S(t)(\mu+\nu)](x)+2s\le w(x,t)
\end{equation}
for $x\in{\mathbb R}^N$ and $0<s<t$. 
By \eqref{eq:2.7} with $\alpha=r$, 
for any $r\in[1,\infty]$, we have 
\begin{equation}
\label{eq:2.38}
\|w(t)\|_{L^r_{{\rm ul}}}\le C(t^{-\frac{N}{2}\left(1-\frac{1}{r}\right)}+ t)\le Ct^{-\frac{N}{2}\left(1-\frac{1}{r}\right)},
\quad t\in(0,1).
\end{equation}

We prove that $(w,w)$ is a supersolution to problem~\eqref{eq:P} in ${\mathbb R}^N\times(0,T)$ 
for some $T\in(0,1)$. 
Since $1<q<1+2/N$, it follows from \eqref{eq:2.19}, \eqref{eq:2.37}, and \eqref{eq:2.38} that 
\begin{equation}
\label{eq:2.39}
\begin{split}
 & \int_0^t [S(D_2(t-s))w(s)^q](x)\,ds
 \le D^{-\frac{N}{2}}\int_0^t [S(t-s)w(s)^q](x)\,ds\\
 & \le D^{-\frac{N}{2}}w(x,t)\int_0^t \|w(s)\|_{L^\infty}^{q-1}\,ds
 \le CD^{-\frac{N}{2}}t^{1-\frac{N}{2}(q-1)}w(x,t)
\end{split}
\end{equation}
for $(x,t)\in{\mathbb R}^N\times(0,1)$. 
Taking small enough $T\in(0,1)$,  
by \eqref{eq:2.39} we have 
\begin{equation}
\label{eq:2.40}
\begin{split}
 & S(D_2t)\nu+\int_0^t [S(D_2(t-s))w(s)^q](x)\,ds\\
 & \le D^{-\frac{N}{2}}S(t)\nu+CD^{-\frac{N}{2}}t^{1-\frac{N}{2}(q-1)}w(x,t)\\
 & \le \left(\frac{1}{2}+CD^{-\frac{N}{2}}T^{1-\frac{N}{2}(q-1)}\right)w(x,t)
\le w(x,t),
\quad (x,t)\in{\mathbb R}^N\times(0,T).
\end{split}
\end{equation}
On the other hand, 
it follows from $0<p\le q$ that $a^p\le(a+1)^p\le Ca^q+1$ for $a\ge 0$. 
Then, similarly to \eqref{eq:2.39}, 
we have 
\begin{equation}
\label{eq:2.41}
\begin{split}
\int_0^t [S(D_1(t-s))w(s)^p](x)\,ds
 & \le t+C\int_0^t [S(D_1(t-s))w(s)^q](x)\,ds\\
 & \le t+CD^{-\frac{N}{2}}t^{1-\frac{N}{2}(q-1)}w(x,t)
\end{split}
\end{equation}
for $(x,t)\in{\mathbb R}^N\times(0,1)$. 
Taking small enough $T\in(0,1)$ if necessary,  
by \eqref{eq:2.41} we see that 
\begin{align*}
 & S(D_1t)\mu+\int_0^t [S(D_1(t-s))w(s)^p](x)\,ds\\
 & \le D^{-\frac{N}{2}}S(t)\mu+t+CD^{-\frac{N}{2}}t^{1-\frac{N}{2}(q-1)}w(x,t)\\
 & \le
 \left(\frac{1}{2}+CD^{-\frac{N}{2}}T^{1-\frac{N}{2}(q-1)}\right)w(x,t)\le w(x,t), 
 \quad (x,t)\in{\mathbb R}^N\times(0,T).
\end{align*}
This together with \eqref{eq:2.40} implies that $(w,w)$ is a supersolution to problem~\eqref{eq:P} in ${\mathbb R}^N\times(0,T)$. 
By the arguments in Section~2.1 we find a solution to problem~\eqref{eq:P} in ${\mathbb R}^N\times(0,T)$ such that 
\begin{equation}
\label{eq:2.42}
0\le u(x,t)\le w(x,t),\quad 0\le v(x,t)\le w(x,t),\quad
(x,t)\in{\mathbb R}^N\times(0,T). 
\end{equation}
Then \eqref{eq:1.10} follows from \eqref{eq:2.38}.
Furthermore, we deduce from \eqref{eq:2.38}, \eqref{eq:2.39}, \eqref{eq:2.41}, and \eqref{eq:2.42} that 
$$
\|u(t)-S(D_1t)\mu\|_{L^1_{{\rm ul}}}+\|v(t)-S(D_2t)\mu\|_{L^1_{{\rm ul}}}\le Ct^{1-\frac{N}{2}(q-1)}\|w(t)\|_{L^1_{{\rm ul}}}+Ct\to 0
$$
as $t\to +0$. Thus \eqref{eq:1.11} holds, and the proof of Theorem~\ref{Theorem:1.2} is complete.
$\Box$
\section{Decay estimates in weak Zygmund type spaces}
In this section we obtain some properties of our weak Zygmund type spaces 
$L^{r,\infty}\Phi(L)^\alpha$, ${\mathfrak L}^{r,\infty}\Phi({\mathfrak L})^\alpha$, 
$L_{{\rm ul}}^{r,\infty}\Phi(L)^\alpha$, and ${\mathfrak L}_{{\rm ul}}^{r,\infty}\Phi({\mathfrak L})^\alpha$. 
Furthermore, 
we develop the arguments in \cite{IIK}*{Section~3} to establish decay estimates of $S(t)\varphi$ in our weak Zygmund type spaces. 
Throughout this paper, for any $r\in[1,\infty]$, we denote by $r'$ the H\"older conjugate of $r$, that is, 
$r'=r/(r-1)$ if $r\in(1,\infty)$, $r'=\infty$ if $r=1$, and $r'=1$ if $r=\infty$.
\subsection{Preliminary lemmas}
We recall some properties of the non-increasing rearrangement $f^*$ 
and its averaging $f^{**}$ for $f\in{\mathcal L}$.
\begin{itemize}
  \item[(a)] 
  Since $f^*$ is non-increasing in $(0,\infty)$, it follows that 
  \begin{equation}
  \label{eq:3.1}
  f^{**}(s)\ge f^*(s),\quad s\in(0,\infty).
  \end{equation}
  \item[(b)] 
  For any $r\in[1,\infty)$, 
  Jensen's inequality together with \eqref{eq:1.3} and \eqref{eq:3.1} yields
  $$
  (f^{**}(s))^r\le \frac{1}{s}\int_0^s f^*(s)^r\,ds=\frac{1}{s}\int_0^s (|f|^r)^*(s)=(|f|^r)^{**}(s),\quad s\in(0,\infty).
  $$
  \item[(c)] 
  It follows from \cite{BS}*{Chapter 2, Proposition 3.3} that 
  \begin{equation}
  \label{eq:3.2}
  f^{**}(s)=\frac{1}{s}\int_0^s f^*(\tau)\,d\tau=\frac{1}{s}\sup_{|E|=s}\int_E |f(x)|\,dx,\quad s\in(0,\infty).
 \end{equation}
  \item[(d)]
  {(O'Neil's inequality)} For any $f_1$, $f_2\in{\mathcal L}$, it follows from \cite{ONeil}*{Lemma~1.6} that 
  \begin{equation}
  \label{eq:3.3}
  (f_1*f_2)^{**}(s)\le\int_s^\infty  f_1^{**}(\tau)f_2^{**}(\tau)\,d\tau,\quad s\in(0,\infty),
  \end{equation}
  where
  $(f_1*f_2)(x)=\int_{{\mathbb R}^N}f_1(x-y)f_2(y)\,dy$ for almost all $x\in{\mathbb R}^N$.
  \item[(e)]
  For any $f_1$, $f_2\in {\mathcal L}$,
  it follows from \cite{ONeil}*{Theorem 3.3} that
  \begin{equation}
  \label{eq:3.4}
  (f_1f_2)^{**}(s)\le\frac{1}{s}\int_0^s f_1^*(\tau)f_2^*(\tau)\,d\tau,\quad s\in(0,\infty).
  \end{equation}
\end{itemize}
Then, for any $r\in[1,\infty)$ and $\alpha\in[0,\infty)$, we have 
\begin{equation}
\label{eq:3.5}
\begin{split}
 & \|f\|_{{\mathfrak L}^{r,\infty}\Phi({\mathfrak L})^\alpha}=\sup_{s>0}\,\left\{s\Phi(s^{-1})^\alpha (|f|^r)^{**}(s)\right\}^{\frac{1}{r}}\\
 & =\sup_{s>0}\,\left\{\Phi(s^{-1})^\alpha \sup_{|E|=s}\int_E |f(x)|^r\,dx\right\}^{\frac{1}{r}}\\
 & =\sup_{s>0}\,\left\{\Phi(s^{-1})^\alpha \int_0^s (|f|^r)^*(\tau)\,d\tau\right\}^{\frac{1}{r}}
 =\sup_{s>0}\,\left\{\Phi(s^{-1})^\alpha\int_0^s f^*(\tau)^r\,d\tau\right\}^{\frac{1}{r}}\\
 & \ge\sup_{s>0}\,\left\{s\Phi(s^{-1})^\alpha f^*(s)^r\right\}^{\frac{1}{r}}
 =\|f\|_{L^{r,\infty}\Phi(L)^\alpha}.
\end{split}
\end{equation}
Furthermore, we have:
\begin{lemma}
\label{Lemma:3.1}
Let $\Phi$ be a non-decreasing function in $[0,\infty)$ with properties \rm{($\Phi$1)}--\rm{($\Phi$3)}. 
Let $r\in[1,\infty)$ and $\alpha\ge 0$. 
Then
$$
\||f|^k\|_{\Phi, r,\alpha;R}=\|f\|_{\Phi, kr,\alpha;R}^k,
\qquad
||||f|^k|||_{r,\alpha;R}=|||f|||_{kr,\alpha;R}^k,
$$
for $f\in{\mathcal L}$, $k>0$ with $kr\ge 1$, and $R\in(0,\infty]$. 
\end{lemma}
{\bf Proof.}
Let $f\in{\mathcal L}$ and $k>0$ with $kr\ge 1$. 
It follows from \eqref{eq:1.3} and \eqref{eq:3.5} that 
\begin{align*}
\||f|^k\|_{L^{r,\infty}\Phi(L)^\alpha} & =\sup_{s>0}\,\left\{s\Phi(s^{-1})^\alpha (|f|^k)^*(s)^r\right\}^{\frac{1}{r}}\\
 & =\sup_{s>0}\,\left\{s\Phi(s^{-1})^\alpha f^*(s)^{kr}\right\}^{\frac{1}{r}}
 =\|f\|^k_{L^{kr,\infty}\Phi(L)^{\alpha}},\\
\||f|^k\|_{{\mathfrak L}^{r,\infty}\Phi({\mathfrak L})^\alpha} & 
=\sup_{s>0}\,\left\{\Phi(s^{-1})^\alpha\int_0^s ((|f|^k)^r)^*(\tau)\,d\tau\right\}^{\frac{1}{r}}\\
 & =\sup_{s>0}\,\left\{\Phi(s^{-1})^\alpha\int_0^s (|f|^{kr})^*(\tau)\,d\tau\right\}^{\frac{1}{r}}
 =\|f\|_{{\mathfrak L}^{kr,\infty}\Phi({\mathfrak L})^{\alpha}}^k.
\end{align*}
These imply the desired relations with $R=\infty$. 
Furthermore, for any $R\in(0,\infty)$, 
\begin{align*}
& \||f|^k\|_{\Phi, r,\alpha;R}=\sup_{x\in{\mathbb R}^N}\| |f|^k\chi_{B(x,R)}\|_{L^{r,\infty}\Phi(L)^\alpha}
=\sup_{x\in{\mathbb R}^N}\| |f|\chi_{B(x,R)}\|_{L^{kr,\infty}\Phi(L)^{\alpha}}^k
=|||f|||_{\Phi, kr,\alpha;R}^r,\\
& ||||f|^k|||_{\Phi, r,\alpha;R}=\sup_{x\in{\mathbb R}^N}\| |f|^k\chi_{B(x,R)}\|_{{\mathfrak L}^{r,\infty}\Phi({\mathfrak L})^\alpha}
=\sup_{x\in{\mathbb R}^N}\| |f|\chi_{B(x,R)}\|_{{\mathfrak L}^{kr,\infty}\Phi({\mathfrak L})^{\alpha}}^k
=|||f|||_{\Phi, kr,\alpha;R}^k. 
\end{align*}
Thus Lemma~\ref{Lemma:3.1} follows. 
$\Box$
\begin{lemma}
\label{Lemma:3.2}
Let $\Phi$ be a non-decreasing function in $[0,\infty)$ with properties \rm{($\Phi$1)}--\rm{($\Phi$3)}. 
Let $r\in[1,\infty]$ and $\alpha_1$, $\alpha_2\ge 0$ be such that 
\begin{equation}
\label{eq:3.6}
\alpha=\frac{\alpha_1}{r}+\frac{\alpha_2}{r'}.
\end{equation}
Then 
\begin{equation}
\label{eq:3.7}
|||f_1f_2|||_{\Phi,1,\alpha;R}\le |||f_1|||_{\Phi,r,\alpha_1;R} |||f_2|||_{\Phi,r',\alpha_2;R}
\end{equation}
for $f_1$, $f_2\in{\mathcal L}$ and $R\in(0,\infty]$. 
Furthermore, for any $R\in(0,\infty)$, 
there exists $C>0$ such that 
\begin{equation}
\label{eq:3.8}
|||f|||_{\Phi,r_1,\alpha;R} \le C|||f|||_{\Phi,r_2,\beta;R}
\end{equation}
for $f\in{\mathcal L}$, $1\le r_1\le r_2\le\infty$, and $0\le\alpha\le\beta<\infty$. 
\end{lemma}
{\bf Proof.}
It suffices to consider $r\in(1,\infty)$. 
Let $\alpha_1$, $\alpha_2\ge 0$ satisfy \eqref{eq:3.6}. 
Let $f_1$, $f_2\in {\mathcal L}$. 
It follows from H\"older's inequality, \eqref{eq:3.4}, and \eqref{eq:3.5} that
\begin{align*}
 & \|f_1f_2\|_{{\mathfrak L}^{1,\infty}\Phi({\mathfrak L})^\alpha}
 =\sup_{s>0}\,\left\{s\Phi(s^{-1})^\alpha (f_1f_2)^{**}(s)\right\}\\
 & \le\sup_{s>0}\,\left\{\Phi(s^{-1})^\alpha\int_0^s f_1^*(\tau) f_2^*(\tau)\,d\tau\right\}\\
 & \le\sup_{s>0}\,\left\{\Phi(s^{-1})^\alpha\left(\int_0^s f_1^*(\tau)^{r}\,d\tau\right)^{\frac{1}{r}}\left(\int_0^s f_2^*(\tau)^{r'}\,d\tau\right)^{\frac{1}{r'}}\right\}\\
 & \le\sup_{s>0}\,\left\{\Phi(s^{-1})^{\alpha_1}\int_0^s f_1^*(\tau)^{r}\,d\tau\right\}^{\frac{1}{r}}\,\,
 \sup_{s>0}\,\left\{\Phi(s^{-1})^{\alpha_2}\int_0^s f_2^*(\tau)^{r'}\,d\tau\right\}^{\frac{1}{r'}}\\
 & =\|f_1\|_{{\mathfrak L}^{r,\infty}\Phi({\mathfrak L})^{\alpha_1}}\|f_2\|_{{\mathfrak L}^{r',\infty}\Phi({\mathfrak L})^{\alpha_2}}.
\end{align*}
Then
\begin{align*}
|||f_1f_2|||_{\Phi,1,\alpha;R}
 & =\sup_{x\in{\mathbb R}^n}\|f_1f_2\chi_{B(x,R)}\|_{{\mathfrak L}^{1,\infty}\Phi({\mathfrak L})^\alpha}\\
 & \le \sup_{x\in{\mathbb R}^n}\|f_1\chi_{B(x,R)}\|_{{\mathfrak L}^{r,\infty}\Phi({\mathfrak L})^{\alpha_1}}\cdot 
  \sup_{x\in{\mathbb R}^n}\|f_2\chi_{B(x,R)}\|_{{\mathfrak L}^{r',\infty}\Phi({\mathfrak L})^{\alpha_2}}\\
 & =|||f_1|||_{\Phi,r,\alpha_1;R}|||f_2|||_{\Phi,r',\alpha_2;R}
\end{align*}
for $R\in(0,\infty]$. 
This implies \eqref{eq:3.7}. 

Let $R\in(0,\infty)$.
It follows from the monotonicity and ($\Phi$1) that $\Phi(\tau)\ge 1$ for $\tau\in[0,\infty)$. 
Then, by Lemma~\ref{Lemma:3.1} and \eqref{eq:3.7} we have
\begin{align*}
|||f|||_{\Phi,r_1,\alpha;R} & =||||f|^{r_1}|||_{\Phi,1,\alpha;R}^{\frac{1}{r_1}}
\le ||||f|^{r_1}|||_{\Phi,\frac{r_2}{r_1},\alpha;R}^{\frac{1}{r_1}}|||1|||_{\Phi,\left(\frac{r_2}{r_1}\right)',\alpha;R}^{\frac{1}{r_1}}
\le C|||f|||_{\Phi,r_2,\alpha;R}\\
 & =C\sup_{x\in{\mathbb R}^N}\sup_{s>0}\,\left\{s\Phi(s^{-1})^\alpha (|f\chi_{B(x,R)}|^{r_2})^{**}\right\}^{\frac{1}{r_2}}\\
 & \le C\sup_{x\in{\mathbb R}^N}\sup_{s>0}\,\left\{s\Phi(s^{-1})^\beta (|f\chi_{B(x,R)}|^{r_2})^{**}\right\}^{\frac{1}{r_2}}
 = C|||f|||_{\Phi,r_2,\beta;R} 
\end{align*}
for $f\in{\mathcal L}$, $1\le r_1\le r_2\le\infty$, and $0\le\alpha\le\beta<\infty$.
Thus \eqref{eq:3.8} holds, and the proof of Lemma~\ref{Lemma:3.2} is complete. 
$\Box$
\vspace{5pt}

Next, we recall the following two lemmas on Hardy's inequality.
(See \cite{Muckenhoup}*{Theorems~1 and 2}.)  
\begin{lemma}
\label{Lemma:3.3}
Let $r\in[1,\infty]$. 
Let $U$ and $V$ be locally integrable functions in $[0,\infty)$. 
Then there exists $C>0$ such that 
$$
\|U\tilde{f}\|_{L^r((0,\infty))}
\le C\|Vf\|_{L^r((0,\infty))}
\quad\mbox{with}\quad
\tilde{f}(s)\coloneqq\int_0^s f(\tau)\,d\tau
$$
holds for locally integrable functions $f$ in $[0,\infty)$ 
if and only if 
$$
\sup_{s>0}\,
\left\{\|U\|_{L^r((s,\infty))}\|V^{-1}\|_{L^{r'}((0,s))}\right\}<\infty.
$$
\end{lemma}
\begin{lemma}
\label{Lemma:3.4}
Let $r\in[1,\infty]$. 
Let $U$ and $V$ be locally integrable functions in $[0,\infty)$. 
Then there exists $C>0$ such that 
$$
\|U\hat{f}\|_{L^r((0,\infty))}
\le C\|Vf\|_{L^r((0,\infty))}
\quad\mbox{with}\quad
\hat{f}(s)\coloneqq\int_s^\infty f(\tau)\,d\tau
$$
holds for locally integrable functions $f$ in $(0,\infty)$ with $f\in L^1((1,\infty))$ if and only if 
$$
\sup_{s>0}\left\{\|U\|_{L^r((0,s))}\|V^{-1}\|_{L^{r'}((s,\infty))}\right\}<\infty.
$$
\end{lemma}
\subsection{Decay estimates}
In this subsection we prove the following proposition on 
decay estimates of $S(t)\varphi$ in weak Zygmund type spaces $L^{r,\infty}\Phi(L)^\alpha$ and ${\mathfrak L}^{r,\infty}\Phi({\mathfrak L})^\alpha$. 
\begin{proposition}
\label{Proposition:3.1}
Let $\Phi$ be a non-decreasing function in $[0,\infty)$ with properties \rm{($\Phi$1)}--\rm{($\Phi$3)}. 
Let $1\le r_1\le r_2\le \infty$ and $\alpha$, $\beta\ge 0$. 
Assume that $\alpha\le\beta$ if $r_1=r_2$. 
\begin{itemize}
  \item[{\rm (1)}] 
  There exists $C_1>0$ such that
  $$
  \left\|S(t)\varphi\right\|_{{\mathfrak L}^{r_2,\infty}\Phi({\mathfrak L})^\beta}
  \le C_1t^{-\frac{N}{2}\left(\frac{1}{r_1}-\frac{1}{r_2}\right)}\Phi(t^{-1})^{-\frac{\alpha}{r_1}+\frac{\beta}{r_2}}\|\varphi\|_{{\mathfrak L}^{r_1,\infty}\Phi({\mathfrak L})^\alpha},
  \quad t>0,
  $$
  for $\varphi\in {\mathfrak L}^{r_1,\infty}\Phi({\mathfrak L})^\alpha$.
  \item[{\rm (2)}] 
  Let $r_1>1$. There exists $C_2>0$ such that
  \begin{equation}
  \label{eq:3.9}
  \left\|S(t)\varphi\right\|_{L^{r_2,\infty}\Phi(L)^\beta}
  \le C_2t^{-\frac{N}{2}\left(\frac{1}{r_1}-\frac{1}{r_2}\right)}\Phi(t^{-1})^{-\frac{\alpha}{r_1}+\frac{\beta}{r_2}}\|\varphi\|_{L^{r_1,\infty}\Phi(L)^\alpha},
  \quad t>0,
  \end{equation}
  for $\varphi\in L^{r_1,\infty}\Phi(L)^\alpha$. 
  \item[{\rm (3)}] 
  Assume that $1<r_1<r_2$. Then there exists $C_3>0$ such that
   $$
   \left\|S(t)\varphi\right\|_{{\mathfrak L}^{r_2,\infty}\Phi({\mathfrak L})^\beta}
   \le C_3t^{-\frac{N}{2}\left(\frac{1}{r_1}-\frac{1}{r_2}\right)}\Phi(t^{-1})^{-\frac{\alpha}{r_1}+\frac{\beta}{r_2}}\|\varphi\|_{L^{r_1,\infty}\Phi(L)^\alpha},\quad t>0,
   $$
   for $\varphi\in L^{r_1,\infty}\Phi(L)^\alpha$.
\end{itemize}
\end{proposition}
At the end of this subsection, 
as an application of Proposition~\ref{Proposition:3.1}, 
we establish decay estimates of $S(t)\varphi$ in uniformly local weak Zygmund type spaces 
$L_{{\rm ul}}^{r,\infty}\Phi(L)^\alpha$ and ${\mathfrak L}_{{\rm ul}}^{r,\infty}\Phi({\mathfrak L})^\alpha$. 

For the proof of Proposition~\ref{Proposition:3.1}, 
we prepare the following four lemmas on $\Phi$.
\begin{lemma}
\label{Lemma:3.5}
Assume the same conditions as in Proposition~{\rm \ref{Proposition:3.1}}. 
\begin{itemize}
  \item[{\rm (1)}] 
  For any fixed $k>0$, 
  $$
  \Phi(a+k)\asymp \Phi(ka)\asymp \Phi(a^k)\asymp \Phi(a)
  $$
  for $a\in(0,\infty)$. 
  \item[{\rm (2)}] 
  Let $\alpha\in{\mathbb R}$ and $\delta>0$. 
  Then there exists $C>0$ such that 
  $$
  \tau_1^\delta \Phi(\tau_1^{-1})^\alpha\le C\tau_2^\delta \Phi(\tau_2^{-1})^\alpha,
  \qquad
  \tau_1^{-\delta} \Phi(\tau_1^{-1})^\alpha\ge C^{-1}\tau_2^{-\delta} \Phi(\tau_2^{-1})^\alpha,
  $$
  for $\tau_1$, $\tau_2\in(0,\infty)$ with $\tau_1\le\tau_2$.
\end{itemize}
\end{lemma}
{\bf Proof.}
We prove assertion~(1). 
It suffices to consider the case where $k>1$ and $a$ is large enough. 
Let $\ell$ be a natural number such that $k\le 2^\ell$. 
Since $\Phi$ is non-decreasing in $[0,\infty)$, by~$(\Phi 2)$ we see that 
$$
\Phi(a)\le\Phi(a+k)\le \Phi(ka)\le \Phi(a^k)
\le \Phi(a^{2^\ell})\le C\Phi(a^{2^{\ell-1}})\le\cdots\le C\Phi(a)
$$
for large enough $a$. Thus assertion~(1) follows. 

We prove assertion~(2). 
Since $\Phi$ is non-decreasing in $[0,\infty)$, 
for any $\alpha\in{\mathbb R}$ and $\delta>0$, 
by~$(\Phi3)$ we find $\tau_*>0$ such that 
the desired inequalities hold for 
$0<\tau_1\le \tau_2\le \tau_*$.
In particular, we have 
\begin{equation}
  \label{eq:3.10}
  \tau^\delta \Phi(\tau^{-1})^\alpha\le C\tau_*^\delta\Phi(\tau_*^{-1})^\alpha, 
  \qquad 
  \tau^{-\delta} \Phi(\tau^{-1})^\alpha\ge C^{-1} \tau_*^{-\delta}\Phi(\tau_*^{-1})^\alpha,
\end{equation}
for $0<\tau\le \tau_*$. 
On the other hand, it follows from the monotonicity of $\Phi$ and ($\Phi$1) that
$$
C^{-1}\le \Phi(\tau^{-1})\le C,\quad \tau\in[\tau_*,\infty). 
$$ 
Then we observe from \eqref{eq:3.10} that 
\begin{align*}
 & \tau_1^\delta \Phi(\tau_1^{-1})^\alpha\le C\tau_*^\delta\Phi(\tau_*^{-1})^\alpha\le C\tau_2^\delta\Phi(\tau_2^{-1})^\alpha
 \quad\mbox{if}\quad \tau_1\le\tau_*\le\tau_2,\\
 & \tau_1^\delta \Phi(\tau_1^{-1})^\alpha\le C\tau_2^\delta\Phi(\tau_2^{-1})^\alpha\quad\mbox{if}\quad \tau_*\le\tau_1\le\tau_2.
\end{align*}
Similarly, we have
\begin{align*}
 & \tau_1^{-\delta} \Phi(\tau_1^{-1})^\alpha\ge C\tau_*^{-\delta}\Phi(\tau_*^{-1})^\alpha\ge C\tau_2^{-\delta}\Phi(\tau_2^{-1})^\alpha
\quad\mbox{if}\quad \tau_1\le\tau_*\le\tau_2,\\
 & \tau_1^{-\delta} \Phi(\tau_1^{-1})^\alpha\ge C\tau_2^{-\delta}\Phi(\tau_2^{-1})^\alpha\quad\mbox{if}\quad \tau_*\le\tau_1\le\tau_2.
\end{align*}
Thus assertion~(2) follows. 
The proof is complete.
$\Box$
\begin{lemma}
\label{Lemma:3.6}
Assume the same conditions as in Proposition~{\rm \ref{Proposition:3.1}}. 
\begin{itemize}
  \item[{\rm (1)}] 
  Let $q>-1$ and $\alpha\in{\mathbb R}$. Then there exists $C_1>0$ such that 
  $$
  \int_0^s \tau^q\Phi(\tau^{-1})^\alpha\,d\tau\le C_1s^{q+1}\Phi(s^{-1})^\alpha,
  \quad s>0.
  $$
  \item[{\rm (2)}]  
  Let $q<-1$ and $\alpha\in{\mathbb R}$. Then there exists $C_2>0$ such that 
  $$
  \int_s^\infty \tau^q\Phi(\tau^{-1})^\alpha\,d\tau\le C_2s^{q+1}\Phi(s^{-1})^\alpha,
  \quad s>0.
  $$
\end{itemize}
\end{lemma}
{\bf Proof.} 
We prove assertion~(1). 
Let $\delta>0$ be such that $q-\delta>-1$. 
By Lemma~\ref{Lemma:3.5}~(2) we have
\begin{align*}
\int_0^s \tau^q\Phi(\tau^{-1})^\alpha\,d\tau
 & =\int_0^s \tau^{q-\delta}\cdot\tau^\delta\Phi(\tau^{-1})^\alpha\,d\tau\\
 & \le C s^\delta\Phi(s^{-1})^\alpha\int_0^s \tau^{q-\delta}\,d\tau
 \le Cs^{q+1}\Phi(s^{-1})^\alpha,\quad s>0.
\end{align*}
Thus assertion~(1) follows. 

We prove assertion~(2). 
Let $\epsilon>0$ be such that $q+\epsilon<-1$. 
Similarly to the proof of assertion~(1), 
by Lemma~\ref{Lemma:3.5}~(2) we see that 
\begin{align*}
\int_s^\infty \tau^q\Phi(\tau^{-1})^\alpha\,d\tau
 & =\int_s^\infty \tau^{q+\epsilon}\cdot\tau^{-\epsilon}\Phi(\tau^{-1})^\alpha\,d\tau\\
 & \le Cs^{-\epsilon}\Phi(s^{-1})^\alpha\int_s^\infty \tau^{q+\epsilon}\,d\tau
 \le Cs^{q+1}\Phi(s^{-1})^\alpha,\quad s>0.
\end{align*}
Thus assertion~(2) follows. The proof is complete.
$\Box$
\begin{lemma}
\label{Lemma:3.7}
Assume the same conditions as in Proposition~{\rm \ref{Proposition:3.1}}. 
\begin{itemize}
  \item[{\rm (1)}] 
  Let $1\le r<\infty$ and $\alpha\ge 0$. Then  
  $$
  \sup_{s>0}\left\{s^{\frac{1}{r}}\Phi(s^{-1})^{\frac{\alpha}{r}} f^{**}(s)\right\}\le
  \|f\|_{{\mathfrak L}^{r,\infty}\Phi({\mathfrak L})^\alpha},
  \quad
  f\in {\mathcal L}.
  $$
  \item[{\rm (2)}]  
  Let $1<r<\infty$ and $\alpha\ge 0$. Then there exists $C>0$ such that 
  $$
  \sup_{s>0}\left\{s^{\frac{1}{r}}\Phi(s^{-1})^{\frac{\alpha}{r}}f^{**}(s)\right\}\le
  C\|f\|_{L^{r,\infty}\Phi(L)^\alpha},
  \quad
  f\in {\mathcal L}.
  $$
\end{itemize}
 \end{lemma}
{\bf Proof.}
Let $f\in{\mathcal L}$. 
For any $r\in[1,\infty)$, it follows from Jensen's inequality and \eqref{eq:1.3} that 
\begin{align*}
 & \sup_{s>0}\left\{s^{\frac{1}{r}}\Phi(s^{-1})^{\frac{\alpha}{r}}  f^{**}(s)\right\}\\
 & \le \sup_{s>0}\left\{s^{\frac{1}{r}}\Phi(s^{-1})^{\frac{\alpha}{r}} \left(s^{-1}\int_0^s f^*(\tau)^r\,d\tau\right)^{\frac{1}{r}}\right\}
=\sup_{s>0}\left\{\Phi(s^{-1})^\alpha \int_0^s (|f|^r)^*(\tau)\,d\tau\right\}^{\frac{1}{r}}\\
 & =\sup_{s>0}\left\{s\Phi(s^{-1})^\alpha (|f|^r)^{**}(s)\right\}^{\frac{1}{r}}
 =\|f\|_{{\mathfrak L}^{r,\infty}\Phi({\mathfrak L})^\alpha},
\end{align*}
which implies assertion~(1). 

Let $r\in(1,\infty)$, and set 
$U(\tau)\coloneqq\tau^{\frac{1}{r}-1}\Phi(\tau^{-1})^{\frac{\alpha}{r}}$ and 
$V(\tau)\coloneqq\tau^{\frac{1}{r}}\Phi(\tau^{-1})^{\frac{\alpha}{r}}$ for $\tau>0$. 
It follows from Lemma~\ref{Lemma:3.5}~(2) and Lemma~\ref{Lemma:3.6}~(1) that
$$
\sup_{s>0}\left\{\|U\|_{L^\infty((s,\infty))}\int_0^s |V(\tau)|^{-1}\,d\tau\right\}
\le \sup_{s>0}\left\{Cs^{\frac{1}{r}-1}\Phi(s^{-1})^{\frac{\alpha}{r}}\cdot Cs^{1-\frac{1}{r}}\Phi(s^{-1})^{-\frac{\alpha}{r}}\right\}<\infty.
$$
This together with Lemma~\ref{Lemma:3.3} with $r=\infty$ implies that 
\begin{align*}
 & \sup_{s>0}\left\{s^{\frac{1}{r}}\Phi(s^{-1})^{\frac{\alpha}{r}} f^{**}(s)\right\}
=\sup_{s>0}\left\{U(s)\int_0^s f^*(s)\,ds\right\}\\
 & \le C\sup_{s>0}\,\{V(s)f^*(s)\}
=C\sup_{s>0}\left\{s^{\frac{1}{r}}\Phi(s^{-1})^{\frac{\alpha}{r}} f^*(s)\right\}\\
 & =C\sup_{s>0}\,\left\{s\Phi(s^{-1})^\alpha f^*(s)^r\right\}^{\frac{1}{r}}
=C\|f\|_{L^{r,\infty}\Phi(L)^\alpha},
\end{align*}
which implies assertion~(2). 
Thus Lemma~\ref{Lemma:3.7} follows.
$\Box$

\begin{lemma}
\label{Lemma:3.8}
Assume the same conditions as in Proposition~{\rm \ref{Proposition:3.1}}. 
Let $1\le r\le q<\infty$ and $\gamma\in{\mathbb R}$. 
Then there exists $C>0$ such that
\begin{equation}
\label{eq:3.11}
\int_0^\infty\tau^{q\left(1-\frac{1}{r}\right)}\Phi(\tau^{-1})^\gamma g_t^*(\tau)^q\,d\tau
\le Ct^{-\frac{Nq}{2}\left(\frac{1}{r}-\frac{1}{q}\right)}\Phi(t^{-1})^\gamma,\quad t>0,
\end{equation}
where $g_t(x)\coloneqq G(x,t)$. 
\end{lemma}
{\bf Proof.}
For any $t>0$, it follows that
\begin{equation}
\label{eq:3.12}
g_t^*(s)=(4\pi t)^{-\frac{N}{2}}\exp\left(-\frac{(\omega_N^{-1}s)^{\frac{2}{N}}}{4t}\right),
\quad s>0.
\end{equation}
Then 
\begin{equation}
\label{eq:3.13}
\begin{split}
I & \coloneqq\int_0^\infty \tau^{q\left(1-\frac{1}{r}\right)}\Phi(\tau^{-1})^\gamma g_t^*(\tau)^q\,d\tau\\
 & \,\le Ct^{-\frac{Nq}{2}}\int_0^\infty \tau^{q\left(1-\frac{1}{r}\right)}\Phi(\tau^{-1})^\gamma
 \exp\left(-\frac{\tau^{\frac{2}{N}}}{Ct}\right)\,d\tau\\
 & \,\le Ct^{-\frac{Nq}{2}\left(\frac{1}{r}-\frac{1}{q}\right)}
 \int_0^\infty \xi^{Nq\left(1-\frac{1}{r}\right)+N-1}e^{-C^{-1}\xi^2}
 \Phi\left((t^{1/2}\xi)^{-N}\right)^\gamma\,d\xi,
 \quad t>0.
\end{split}
\end{equation}
Let $\epsilon>0$ be small enough. 
Then, by Lemma~\ref{Lemma:3.5} we have 
\begin{equation}
\label{eq:3.14}
\begin{split}
\Phi\left((t^{1/2}\xi)^{-N}\right)^\gamma
 & \le C(t^{1/2}\xi)^{-\epsilon}(t^{1/2}\xi)^\epsilon \Phi\left((t^{1/2}\xi)^{-1}\right)^\gamma\\
 & \le C(t^{1/2}\xi)^{-\epsilon}(t^{1/2})^\epsilon \Phi\left((t^{1/2})^{-1}\right)^\gamma
\le C\xi^{-\epsilon}\Phi(t^{-1})^\gamma,
\quad \xi\in(0,1].
\end{split}
\end{equation}
Similarly, we see that 
\begin{equation}
\label{eq:3.15}
\begin{split}
\Phi\left((t^{1/2}\xi)^{-N}\right)^\gamma
 & \le C(t^{1/2}\xi)^\epsilon(t^{1/2}\xi)^{-\epsilon}\Phi\left((t^{1/2}\xi)^{-1}\right)^\gamma\\
 & \le C(t^{1/2}\xi)^\epsilon (t^{1/2})^{-\epsilon} \Phi\left((t^{1/2})^{-1}\right)^\gamma
\le C\xi^\epsilon \Phi(t^{-1})^\gamma,
\quad \xi\in(1,\infty).
\end{split}
\end{equation}
Combining \eqref{eq:3.13}, \eqref{eq:3.14}, and \eqref{eq:3.15}, 
we obtain 
\begin{align*}
I & \le Ct^{-\frac{Nq}{2}\left(\frac{1}{r}-\frac{1}{q}\right)}\Phi(t^{-1})^\gamma
\int_0^\infty \xi^{Nq\left(1-\frac{1}{r}\right)+N-1}(\xi^{-\epsilon}+\xi^\epsilon)e^{-C^{-1}\xi^2}\,d\xi\\
 & \le Ct^{-\frac{Nq}{2}\left(\frac{1}{r}-\frac{1}{q}\right)}\Phi(t^{-1})^\gamma,
 \quad t>0. 
\end{align*}
Thus \eqref{eq:3.11} holds, and the proof is complete.
$\Box$
\vspace{5pt}

Now we are ready to prove Proposition~\ref{Proposition:3.1}. 
We first prove Proposition~\ref{Proposition:3.1}~(1) and (3). 
\vspace{3pt}
\newline
{\bf Proof of Proposition~\ref{Proposition:3.1}~(1) and (3).}
The proof is divided into the following three cases: 
$$
1\le r_1<r_2<\infty;\qquad 1\le r_1=r_2<\infty;\qquad 1\le r_1\le r_2=\infty.
$$ 
\underline{Step 1.}
Consider the case of $1\le r_1<r_2<\infty$. 
By \eqref{eq:3.5} it suffices to prove 
\begin{equation}
\label{eq:3.16}
\left\|S(t)\varphi\right\|_{{\mathfrak L}^{r_2,\infty}\Phi({\mathfrak L})^\beta}
\le C_1t^{-\frac{N}{2}\left(\frac{1}{r_1}-\frac{1}{r_2}\right)}\Phi(t^{-1})^{-\frac{\alpha}{r_1}+\frac{\beta}{r_2}}\|\varphi\|_{X^{r_1,\alpha}},
\quad t>0,
\end{equation}
where 
$X^{r,\alpha}\coloneqq{\mathfrak L}^{1,\infty}\Phi({\mathfrak L})^\alpha$ if $r=1$ and 
$X^{r,\alpha}\coloneqq L^{r,\infty}\Phi(L)^\alpha$ if $r>1$. 
It follows from \eqref{eq:3.1}, \eqref{eq:3.3}, and \eqref{eq:3.5} that 
\begin{align*}
\left\|S(t)\varphi\right\|_{{\mathfrak L}^{r_2,\infty}\Phi({\mathfrak L})^\beta}^{r_2}
 & =\sup_{s>0}\,\left\{\Phi(s^{-1})^\beta\int_0^s \left(\left(S(t)\varphi\right)^*(\tau)\right)^{r_2}\,d\tau\right\}\\
 & \le \sup_{s>0}\,\left\{\Phi(s^{-1})^\beta\int_0^s \left(\left(S(t)\varphi\right)^{**}(\tau)\right)^{r_2}\,d\tau\right\}\\
 & \le \sup_{s>0}\,\left\{\Phi(s^{-1})^\beta\int_0^s \left(\int_\tau^\infty g_t^{**}(\eta)\varphi^{**}(\eta)\,d\eta\right)^{r_2}\,d\tau\right\},
 \quad t>0.
\end{align*}
Since $\Phi(s^{-1})^\beta$ is non-increasing for $s\in(0,\infty)$, we have
\begin{equation}
\label{eq:3.17}
\left\|S(t)\varphi\right\|_{{\mathfrak L}^{r_2,\infty}\Phi({\mathfrak L})^\beta}^{r_2}
\le\int_0^\infty \left(\Phi(\tau^{-1})^{\frac{\beta}{r_2}}\int_\tau^\infty g_t^{**}(\eta)\varphi^{**}(\eta)\,d\eta\right)^{r_2}\,d\tau,
\quad t>0.
\end{equation}
Set
$U(\tau)\coloneqq\Phi(\tau^{-1})^{\frac{\beta}{r_2}}$ and 
$V(\tau)\coloneqq\tau\Phi(\tau^{-1})^{\frac{\beta}{r_2}}$ for $\tau>0$. 
It follows from Lemma~\ref{Lemma:3.6} that
\begin{align*}
 & \sup_{s>0}\,\left(\int_0^s |U(\tau)|^{r_2}\,d\tau\right)^{\frac{1}{r_2}}\left(\int_s^\infty |V(\tau)|^{-r_2'}\,d\tau\right)^{\frac{1}{r_2'}}\\
 & \le \sup_{s>0}\,\left\{Cs^{\frac{1}{r_2}}\Phi(\tau^{-1})^{\frac{\beta}{r_2}}\cdot Cs^{-1+\frac{1}{r_2'}}\Phi(\tau^{-1})^{-\frac{\beta}{r_2}}\right\}<\infty.
\end{align*}
Then, by Lemma~\ref{Lemma:3.4}, Lemma~\ref{Lemma:3.7}, and \eqref{eq:3.17} we have 
\begin{equation}
\label{eq:3.18}
\begin{split}
 & \left\|S(t)\varphi\right\|_{{\mathfrak L}^{r_2,\infty}\Phi({\mathfrak L})^\beta}^{r_2}
\le C\int_0^\infty \left(\tau\Phi(\tau^{-1})^{\frac{\beta}{r_2}} g_t^{**}(\tau)\varphi^{**}(\tau)\right)^{r_2}\,d\tau\\
 & \le C\sup_{s>0}\,\left\{s^{\frac{1}{r_1}}\Phi(s^{-1})^{\frac{\alpha}{r_1}} \varphi^{**}(s)\right\}^{r_2}
 \int_0^\infty\bigg( \tau^{1-\frac{1}{r_1}}\Phi(\tau^{-1})^{-\frac{\alpha}{r_1}+{\frac{\beta}{r_2}}}g_t^{**}(\tau)\bigg)^{r_2}\,d\tau\\
 & \le C\|\varphi\|_{X^{r_1,\alpha}}^{r_2}\int_0^\infty\bigg( \tau^{1-\frac{1}{r_1}}\Phi(\tau^{-1})^\gamma g_t^{**}(\tau)\bigg)^{r_2}\,d\tau\\
 & =C\|\varphi\|_{X^{r_1,\alpha}}^{r_2}\int_0^\infty\bigg( \tau^{-\frac{1}{r_1}}\Phi(\tau^{-1})^\gamma \int_0^\tau g_t^*(s)\,ds\bigg)^{r_2}\,d\tau, 
 \quad t>0,
\end{split}
\end{equation}
where ${\gamma}=-\frac{\alpha}{r_1}+\frac{\beta}{r_2}$.  

Set
$\tilde{U}(\tau)\coloneqq\tau^{-\frac{1}{r_1}}\Phi(\tau^{-1})^\gamma$ 
and $\tilde{V}(\tau)\coloneqq\tau^{1-\frac{1}{r_1}}\Phi(\tau^{-1})^\gamma$ for $\tau>0$. 
Since $r_2>r_1$ and $r_2'<r_1'$, 
by Lemma~\ref{Lemma:3.6} we have 
\begin{equation}
\label{eq:3.19}
\begin{split}
  & \sup_{s>0}
  \,
  \left\{
  \left(\int_s^\infty |\tilde{U}(\tau)|^{r_2}\,d\tau\right)^{\frac{1}{r_2}}
  \left(\int_0^s |\tilde{V}(\tau)|^{-r_2'}\,{d\tau}\right)^{\frac{1}{r_2'}}
  \right\}
  \\
  & =\sup_{s>0}\,
  \left\{
  \left(\int_s^\infty \tau^{-\frac{r_2}{r_1}}\Phi(\tau^{-1})^{r_2\gamma}\,d\tau\right)^{\frac{1}{r_2}}
  \left(\int_0^s \tau^{-\frac{r_2'}{r_1'}}\Phi(\tau^{-1})^{-r_2'{\gamma}}\,{d\tau}\right)^{\frac{1}{r_2'}}
  \right\}
  \\
  & \le \sup_{s>0}\,\left\{C s^{\frac{1}{r_2}-\frac{1}{r_1}}\Phi(s^{-1})^\gamma \cdot Cs^{\frac{1}{r_2'}-\frac{1}{r_1'}}\Phi(s^{-1})^{-{\gamma}}\right\}<\infty.
\end{split}
\end{equation}
Applying Lemma~\ref{Lemma:3.3} to \eqref{eq:3.18}, by \eqref{eq:3.19}
we obtain 
$$
\left\|S(t)\varphi\right\|_{{\mathfrak L}^{q,\infty}\Phi({\mathfrak L})^\beta}^{r_2}
\le C\|\varphi\|_{X^{r_1,\alpha}}^{r_2}
\int_0^\infty\left(\tau^{1-\frac{1}{r_1}}\Phi(\tau^{-1})^\gamma g_t^*(\tau)\right)^{r_2}\,d\tau, 
\quad t>0. 
$$
This together with Lemma~\ref{Lemma:3.8} implies that
$$
\left\|S(t)\varphi\right\|_{{\mathfrak L}^{r_2,\infty}\Phi({\mathfrak L})^\beta}^{r_2}
\le Ct^{-\frac{Nr_2}{2}\left(\frac{1}{r_1}-\frac{1}{r_2}\right)}\Phi(t^{-1})^{r_2\gamma}\|\varphi\|_{X^{r_1,\alpha}}^{r_2},
\quad t>0. 
$$
Thus inequality~\eqref{eq:3.16} holds, and Proposition~\ref{Proposition:3.1}~(1) and (3) hold in the case of $1\le r_1<r_2<\infty$. 
\newline
\underline{Step 2.}
Consider the case of $1\le r_1=r_2<\infty$. Set $r\coloneqq r_1=r_2$. 
It follows from Jensen's inequality that 
\begin{align*}
\left|\left[S(t)\varphi\right](x)\right|^r
\le\int_{{\mathbb R}^N} g_t(x-y)|\varphi(y)|^r\,dy,
\quad (x,t)\in{\mathbb R}^N\times(0,\infty).
\end{align*}
This together with \eqref{eq:3.3} implies that
\begin{equation}
\label{eq:3.20}
\begin{split}
\|S(t)\varphi\|_{{\mathfrak L}^{r,\infty}\Phi({\mathfrak L})^\beta}^r
 & =\sup_{s>0}\, \left\{s\Phi(s^{-1})^\beta  (|S(t)\varphi|^r)^{**}(s)\right\}\\
 & \le \sup_{s>0}\, \left\{s\Phi(s^{-1})^\beta \int_s^\infty g_t^{**}(\tau)(|\varphi|^r)^{**}(s)\,d\tau\right\},
 \quad t>0. 
\end{split}
\end{equation}
Set 
$\hat{U}(\tau)\coloneqq\tau\Phi(\tau^{-1})^{\beta}$ and $\hat{V}(\tau)\coloneqq\tau^2\Phi(\tau^{-1})^{\beta}$ for $\tau>0$. 
By Lemma~\ref{Lemma:3.5}~(2) and Lemma~\ref{Lemma:3.6}~(2) 
we have 
\begin{equation}
\label{eq:3.21}
\sup_{s>0}\,\left\{\|\hat{U}\|_{L^\infty((0,s))}\int_s^\infty |\hat{V}(\tau)|^{-1}\,d\tau\right\}\le \sup_{s>0}\left\{Cs\Phi(s^{-1})^{\beta}\cdot Cs^{-1}\Phi(s^{-1})^{-\beta}\right\}<\infty.
\end{equation}
Applying Lemma~\ref{Lemma:3.4} with $r=\infty$, by \eqref{eq:3.20} and \eqref{eq:3.21} we obtain 
\begin{equation}
\label{eq:3.22}
\begin{split}
 & \|S(t)\varphi\|_{{\mathfrak L}^{r,\infty}\Phi({\mathfrak L})^\beta}^r
 \le C\sup_{s>0}\, \left\{s^2\Phi(s^{-1})^{\beta} g_t^{**}(s)(|\varphi|^r)^{**}(s)\right\}\\
 & \qquad\quad
 \le C\sup_{s>0}\, \left\{s\Phi(s^{-1})^{\beta-\alpha}g_t^{**}(s)\right\}
\cdot
\sup_{s>0}\left\{s\Phi(s^{-1})^\alpha (|\varphi|^r)^{**}(s)\right\}\\
& \qquad\quad
=C\sup_{s>0}\, \left\{\Phi(s^{-1})^{\beta-\alpha}\int_0^s g_t^*(\tau)\,d\tau\right\}\|\varphi\|_{{\mathfrak L}^{r,\infty}\Phi({\mathfrak L})^\alpha}^r,
\quad t>0.
\end{split}
\end{equation}
Furthermore, 
since $\alpha\le\beta$, 
$\Phi(t^{-1})^{\beta-\alpha}$ is non-increasing in $(0,\infty)$, 
by Lemma~\ref{Lemma:3.8} we have
$$
\sup_{s>0}\, \left\{\Phi(s^{-1})^{\beta-\alpha}\int_0^s g_t^*(\tau)\,d\tau\right\}\\
\le\int_0^\infty \Phi(\tau^{-1})^{\beta-\alpha}g_t^*(\tau)\,d\tau
\le C\Phi(t^{-1})^{\beta-\alpha},
\quad t>0.
$$
This together with \eqref{eq:3.22} implies that 
$$
\|S(t)\varphi\|_{{\mathfrak L}^{r,\infty}\Phi({\mathfrak L})^\beta}^r
\le C\Phi(t^{-1})^{\beta-\alpha}\|\varphi\|^r_{{\mathfrak L}^{r,\infty}\Phi({\mathfrak L})^\alpha},
\quad t>0. 
$$
Thus Proposition~\ref{Proposition:3.1}~(1) holds in the case of $1\le r_1=r_2<\infty$. 
\newline
\underline{Step 3.}
It remains to consider the case of $1\le r_1\le r_2=\infty$. 
Let $X^{r,\alpha}$ be as in Step 1.
If $r_1=r_2=\infty$, then 
\begin{equation}
\label{eq:3.23}
\|S(t)\varphi\|_{L^\infty}\le \|\varphi\|_{L^\infty}\int_{{\mathbb R}^N}g_t(y)\,dy
\le \|\varphi\|_{L^\infty}, \quad t>0,
\end{equation}
and Proposition~\ref{Proposition:3.1}~(1) follows. 
If $1\le r_1<r_2=\infty$, 
by \eqref{eq:3.16} with $r_2=2r_1$ we have 
\begin{align*}
\|S(t)\varphi\|_{L^\infty} & =\left\|S\left(\frac{t}{2}\right)S\left(\frac{t}{2}\right)\varphi\right\|_{L^\infty}
\le Ct^{-\frac{N}{4r_1}}\left\|S\left(\frac{t}{2}\right)\varphi\right\|_{L^{2r_1}}\\
 & =Ct^{-\frac{N}{4r_1}}\left\|S\left(\frac{t}{2}\right)\varphi\right\|_{{\mathfrak L}^{2r_1,\infty}\Phi({\mathfrak L})^0}
 \le Ct^{-\frac{N}{4r_1}}\cdot Ct^{-\frac{N}{2}\left(\frac{1}{r_1}-\frac{1}{2r_1}\right)}
\Phi(t^{-1})^{-\frac{\alpha}{r_1}}\|\varphi\|_{X^{r_1,\alpha}}\\
 & =Ct^{-\frac{N}{2r_1}}\Phi(t^{-1})^{-\frac{\alpha}{r_1}}\|\varphi\|_{X^{r_1,\alpha}},
 \quad t>0.
\end{align*}
Thus Proposition~\ref{Proposition:3.1}~(1) and (3) hold in the case of $1\le r_1<r_2=\infty$. 
Therefore the proof of Proposition~\ref{Proposition:3.1}~(1) and (3) is complete. 
$\Box$\vspace{5pt}

\noindent
{\bf Proof of Proposition~\ref{Proposition:3.1}.}
It remains to prove Proposition~\ref{Proposition:3.1}~(2).
It follows from \eqref{eq:3.1} and \eqref{eq:3.3} that 
\begin{equation*}
	\begin{aligned}
	\left\|S(t)\varphi\right\|_{L^{r_2,\infty}\Phi(L)^\beta}
 	& 
	=\sup_{s>0}\,\left\{s\Phi(s^{-1})^\beta \left(S(t)\varphi\right)^*(s)^{r_2}\right\}^{\frac{1}{r_2}}
	\\
 	& 
	\le \sup_{s>0}\,\left\{s\Phi(s^{-1})^\beta \left(S(t)\varphi\right)^{**}(s)^{r_2} \right\}^{\frac{1}{r_2}}
	\\
	 & 
	 \le \sup_{s>0}\,\left\{s^{\frac{1}{r_2}}\Phi(s^{-1})^{\frac{\beta}{r_2}}
	 \int_s^\infty g_t^{**}(\eta)\varphi^{**}(\eta)\,d\eta\right\},
	 \quad t>0.
	\end{aligned}
\end{equation*}
Set $U(\tau)\coloneqq \tau^{\frac{1}{r_2}}\Phi(\tau^{-1})^{\frac{\beta}{r_2}}$
and $V(\tau)\coloneqq \tau^{1+\frac{1}{r_2}}\Phi(\tau^{-1})^{\frac{\beta}{r_2}}$ for $\tau>0$. 
By Lemma~\ref{Lemma:3.5}~(2) and Lemma~\ref{Lemma:3.6}~(2) we have 
$$
\sup_{s>0}\,\bigg\{\|U\|_{L^\infty((0,s))}\int_s^\infty|V(\tau)|^{-1}\,d\tau\bigg\}
\le \sup_{s>0}\,\bigg\{Cs^\frac{1}{r_2}\Phi(s^{-1})^{\frac{\beta}{r_2}}
 \cdot C s^{-\frac{1}{r_2}}\Phi(s^{-1})^{-\frac{\beta}{r_2}}\bigg\}<\infty.
$$
Then, by Lemma~\ref{Lemma:3.4} with $r=\infty$ and Lemma~\ref{Lemma:3.7}~(2) we obtain 
\begin{equation}
\label{eq:3.24}
\begin{split}
	 & \left\|S(t)\varphi\right\|_{L^{r_2,\infty}\Phi(L)^\beta}
	\le C\sup_{s>0}\,\left\{s^{1+\frac{1}{r_2}}\Phi(s^{-1})^{\frac{\beta}{r_2}}
	g_t^{**}(s)\varphi^{**}(s)\right\}
	\\
 	& \qquad\quad
	\le C\sup_{s>0}\,\left\{s^{\frac{1}{r_1}}\Phi(s^{-1})^{\frac{\alpha}{r_1}} \varphi^{**}(s)\right\}
	\cdot\sup_{s>0}\left\{s^{1+\frac{1}{r_2}-\frac{1}{r_1}}\Phi(s^{-1})^{\frac{\beta}{r_2}-\frac{\alpha}{r_1}}g_t^{**}(s)\right\}\\
	& \qquad\quad
	\le C\|\varphi\|_{L^{r_1,\infty}\Phi(L)^\alpha}\sup_{s>0}\left\{s^{\frac{1}{r_2}-\frac{1}{r_1}}\Phi(s^{-1})^{\frac{\beta}{r_2}-\frac{\alpha}{r_1}}\int_0^s g_t^{*}(\tau)\,d\tau\right\},
	\quad t>0. 
\end{split}
\end{equation}	

Consider the case of $1<r_1<r_2<\infty$. 
Set 
$$
\hat{U}(\tau)\coloneqq\tau^{\frac{1}{r_2}-\frac{1}{r_1}}\Phi(\tau^{-1})^{\frac{\beta}{r_2}-\frac{\alpha}{r_1}},
\quad
\hat{V}(\tau)\coloneqq\tau^{1+\frac{1}{r_2}-\frac{1}{r_1}}\Phi(\tau^{-1})^{\frac{\beta}{r_2}-\frac{\alpha}{r_1}},
$$
for $\tau>0$. 
By Lemma~\ref{Lemma:3.5}~(2) and Lemma~\ref{Lemma:3.6}~(1) we have
\begin{align*}
& \sup_{s>0}\left\{\|\hat{U} \|_{L^\infty((s,\infty))}
\int_0^s |\hat{V}(\tau)|^{-1}\, d\tau
\right\}\\
& \le \sup_{s>0}\left\{Cs^{\frac{1}{r_2}-\frac{1}{r_1}}\Phi(s^{-1})^{\frac{\beta}{r_2}-\frac{\alpha}{r_1}}\cdot 
Cs^{-\frac{1}{r_2}+\frac{1}{r_1}}\Phi(s^{-1})^{-\frac{\beta}{r_2}+\frac{\alpha}{r_1}}\right\}<\infty.
\end{align*}
This together with Lemma~\ref{Lemma:3.3} with $r=\infty$ implies that 
\begin{equation}
\label{eq:3.25}
\sup_{s>0}\left\{s^{\frac{1}{r_2}-\frac{1}{r_1}}\Phi(s^{-1})^{\frac{\beta}{r_2}-\frac{\alpha}{r_1}}\int_0^sg_t^{*}(\tau)\,d\tau\right\}
\le C\sup_{s>0}\left\{s^{1+\frac{1}{r_2}-\frac{1}{r_1}}\Phi(s^{-1})^{\frac{\beta}{r_2}-\frac{\alpha}{r_1}}g_t^*(s)\right\}
\end{equation}	
for $t>0$. 
On the other hand, since $\Phi$ is non-decreasing in $[0,\infty)$, 
it follows from $(\Phi 1)$ and $(\Phi 2)$ that 
\begin{equation}
  \notag 
  \Phi(ab) \le \Phi((\max\{a,b\})^2) 
  \le 
  C \Phi(\max\{a,b\}) 
  \le 
  C \Phi(\max\{a,b\})\Phi(\min\{a,b\}) 
  = 
  C \Phi(a)\Phi(b) 
\end{equation}
for $a$, $b\ge 0$. 
Then, by Lemma~\ref{Lemma:3.5}~(1) and \eqref{eq:3.12} we have
\begin{align*}
 & \sup_{s>0}\left\{s^{1+\frac{1}{r_2}-\frac{1}{r_1}}\Phi(s^{-1})^{\frac{\beta}{r_2}-\frac{\alpha}{r_1}}g_t^*(s)\right\}\\
 & =\sup_{\eta>0}\left\{\Phi\left(\omega_N^{-1}(4t\eta)^{-\frac{N}{2}}\right)^{\frac{\beta}{r_2}-\frac{\alpha}{r_1}}
\left(\omega_N(4t\eta)^{\frac{N}{2}}\right)^{1+\frac{1}{r_2}-\frac{1}{r_1}}
(4\pi t)^{-\frac{N}{2}}e^{-\eta}\right\}\\
 & \le Ct^{-\frac{N}{2}\left(\frac{1}{r_1}-\frac{1}{r_2}\right)}\Phi(t^{-1})^{\frac{\beta}{r_2}-\frac{\alpha}{r_1}}
\sup_{\eta>0}\left\{\eta^{\frac{N}{2}\left(1+\frac{1}{r_2}-\frac{1}{r_1}\right)}\Phi(\eta^{-1})^{\frac{\beta}{r_2}-\frac{\alpha}{r_1}} e^{-\eta}\right\}\\
 & \le  Ct^{-\frac{N}{2}\left(\frac{1}{r_1}-\frac{1}{r_2}\right)}\Phi(t^{-1})^{\frac{\beta}{r_2}-\frac{\alpha}{r_1}},
 \quad t>0.
\end{align*}
This together with \eqref{eq:3.24} and \eqref{eq:3.25} implies \eqref{eq:3.9} in the case of $1<r_1<r_2<\infty$. 

Consider the case of $1<r_1=r_2<\infty$. Set $r\coloneqq r_1=r_2$. Let $\alpha\le\beta$. 
Since $\Phi(t^{-1})^{\beta-\alpha}$ is non-increasing in $(0,\infty)$, 
it follows from Lemma~\ref{Lemma:3.8} that 
$$
\sup_{s>0}\left\{\Phi(s^{-1})^{\frac{\beta}{r}-\frac{\alpha}{r}}\int_0^s g_t^{*}(\tau)\,d\tau\right\}
\le\int_0^\infty \Phi(\tau^{-1})^{\frac{\beta}{r}-\frac{\alpha}{r}}g_t^{*}(\tau)\,d\tau
\le C\Phi(t^{-1})^{\frac{\beta}{r}-\frac{\alpha}{r}},
\quad t>0.
$$
This together with \eqref{eq:3.24} implies \eqref{eq:3.9} in the case of $1<r_1=r_2<\infty$. 
Furthermore, in the case of $1<r_1<r_2=\infty$, 
similarly to Step~3 in the proof of Proposition~\ref{Proposition:3.1}~(1) and (3), 
we have 
\begin{align*}
\|S(t)\varphi\|_{L^\infty}
 & \le Ct^{-\frac{N}{4r_1}}\left\|S\left(\frac{t}{2}\right)\varphi\right\|_{L^{2r_1,\infty}}
 =Ct^{-\frac{N}{4r_1}}\left\|S\left(\frac{t}{2}\right)\varphi\right\|_{L^{2r_1,\infty}\Phi(L)^0}\\
 & \le Ct^{-\frac{N}{4r_1}}\cdot Ct^{-\frac{N}{2}\left(\frac{1}{r_1}-\frac{1}{2r_1}\right)}\Phi(t^{-1})^{-\frac{\alpha}{r_1}}\|\varphi\|_{L^{r_1,\infty}\Phi(L)^\alpha}\\
  & \le Ct^{-\frac{N}{2r_1}}\Phi(t^{-1})^{-\frac{\alpha}{r_1}}\|\varphi\|_{L^{r_1,\infty}\Phi(L)^\alpha},
  \quad t>0.
 \end{align*}
Thus \eqref{eq:3.9} holds in the case of $1<r_1<r_2=\infty$. 
In addition, if $1<r_1=r_2=\infty$, 
\eqref{eq:3.9} follows from \eqref{eq:3.23}. 
Thus \eqref{eq:3.9} holds, and the proof of Proposition~\ref{Proposition:3.1} is complete.
$\Box$
\vspace{5pt}

Furthermore, we apply the same arguments as in the proof of \cite{IIK}*{Proposition~3.2} 
together with Proposition~\ref{Proposition:3.1}, and obtain the following proposition. 
\begin{proposition}
\label{Proposition:3.2}
Let $\Phi$ be a non-decreasing function in $[0,\infty)$ with properties \rm{($\Phi$1)}--\rm{($\Phi$3)}. 
Let $1\le r_1\le r_2\le \infty$, $\alpha$, $\beta\ge 0$, and $R_*\in(0,\infty)$. 
Assume that $\alpha\le\beta$ if $r_1=r_2$. 
\begin{itemize}
  \item[{\rm (1)}] 
  There exists $C_1>0$ such that 
  $$
  |||S(t)\varphi|||_{\Phi,r_2,\beta; R}
  \le C_1t^{-\frac{N}{2}\left(\frac{1}{r_1}-\frac{1}{r_2}\right)}\Phi(t^{-1})^{-\frac{\alpha}{r_1}+\frac{\beta}{r_2}}|||\varphi|||_{\Phi,r_1,\alpha; R}
  $$
  for $\varphi\in {\mathfrak L}^{r_1,\infty}_{{\rm ul}}\Phi({\mathfrak L})^\alpha$, $R\in(0,R_*]$, and $t\in(0,R^2)$.
  \item[{\rm (2)}]
  Let $r_1>1$.  There exists $C_2>0$ such that 
  $$
  \|S(t)\varphi\|_{\Phi,r_2,\beta; R}
  \le C_2t^{-\frac{N}{2}\left(\frac{1}{r_1}-\frac{1}{r_2}\right)}\Phi(t^{-1})^{-\frac{\alpha}{r_1}+\frac{\beta}{r_2}}\|\varphi\|_{\Phi,r_1,\alpha; R}
  $$
  for $\varphi\in L^{r_1,\infty}_{{\rm ul}}\Phi(L)^\alpha$, $R\in(0,R_*]$, and $t\in(0,R^2)$.  
  \item[{\rm (3)}] 
  Let $1<r_1<r_2$. There exists $C_3>0$ such that 
  $$
  |||S(t)\varphi|||_{\Phi,r_2,\beta; R}
  \le C_3t^{-\frac{N}{2}\left(\frac{1}{r_1}-\frac{1}{r_2}\right)}\Phi(t^{-1})^{-\frac{\alpha}{r_1}+\frac{\beta}{r_2}}\|\varphi\|_{\Phi,r_1,\alpha; R}
  $$
  for $\varphi\in L^{r_1,\infty}_{{\rm ul}}\Phi(L)^\alpha$, $R\in(0,R_*]$, and $t\in(0,R^2)$.
\end{itemize}
\end{proposition}

At the end of this section, 
we apply Hardy's inequality again to show that $L^{r,\infty}\Phi(L)^\alpha$ are Banach spaces if $r>1$. 
\begin{lemma}
\label{Lemma:3.9}
Let $\Phi$ be a non-decreasing function in $[0,\infty)$ with properties \rm{($\Phi$1)}--\rm{($\Phi$3)}. 
Let $r\in(1,\infty)$ and $\alpha\in[0,\infty)$. 
For any $f\in{\mathcal L}$, 
set 
\begin{equation}
\label{eq:3.26}
\|f\|_{L^{r,\infty}\Phi(L)^\alpha}':=
\sup_{s>0}\,\left\{s^{\frac{1}{r}}\Phi(s^{-1})^{\frac{\alpha}{r}}f^{**}(s)\right\}.		
\end{equation}
Then there exists $C>0$ such that 
$$
\|f\|_{L^{r,\infty}\Phi(L)^\alpha}\le \|f\|_{L^{r,\infty}\Phi(L)^\alpha}'\le C\|f\|_{L^{r,\infty}\Phi(L)^\alpha},
\quad f\in{\mathcal L}.
$$
Furthermore, 
$L^{r,\infty}\Phi(L)^\alpha$ is a Banach space equipped with the norm $\|\cdot\|'_{L^{r,\infty}\Phi(L)^\alpha}$.
\end{lemma}
{\bf Proof.}
Let $r\in(1,\infty)$ and $\alpha\in[0,\infty)$. 
It follows from \eqref{eq:3.1} and \eqref{eq:3.26} that 
$$
\|f\|_{L^{r,\infty}\Phi(L)^\alpha}'
\ge\sup_{s>0}\left\{s^{\frac{1}{r}}\Phi(s^{-1})^{\frac{\alpha}{r}}f^*(s)\right\}=\|f\|_{L^{r,\infty}\Phi(L)^\alpha}
$$
for $f\in{\mathcal L}$. 
Furthermore, it follows from Lemma~\ref{Lemma:3.7}~(2) that
$$
\|f\|_{L^{r,\infty}\Phi(L)^\alpha}'\le C\|f\|_{L^{r,\infty}\Phi(L)^\alpha}
$$
for $f\in{\mathcal L}$. 
On the other hand,
we observe from \eqref{eq:3.2} that 
\begin{align*}
  & \|f\|_{L^{r,\infty}\Phi(L)^\alpha}'=
\sup_{s>0}\sup_{|E|=s}\left\{s^{-1+\frac{1}{r}}\Phi(s^{-1})^{\frac{\alpha}{r}}\int_E |f(x)|\,dx\right\},
\quad f\in{\mathcal L}.
\end{align*}
Then we easily see that $\|\cdot\|_{L^{r,\infty}\Phi(L)^\alpha}'$ is a norm of $L^{r,\infty}\Phi(L)^\alpha$. 
In addition, we see 
that $L^{r,\infty}\Phi(L)^\alpha$ is a Banach spaces equipped with the norm $\|\cdot\|_{L^{q,\infty}\Phi(L)^\alpha}'$.
Thus Lemma~\ref{Lemma:3.9} follows. 
$\Box$
\vspace{5pt}
\newline
Then we have:
\begin{lemma}
\label{Lemma:3.10}
Let $\Phi$ be a non-decreasing function in $[0,\infty)$ with properties \rm{($\Phi$1)}--\rm{($\Phi$3)}. 
Let $r\in(1,\infty)$ and $\alpha\in[0,\infty)$. 
For any $f\in{\mathcal L}$ and $R\in(0,\infty]$, set 
$$
 \|f\|'_{\Phi,r,\alpha; R}:=\sup_{z\in{\mathbb R}^N}\|f\chi_{B(z,R)}\|'_{L^{r,\infty}\Phi(L)^\alpha}.
$$
Then there exists $C>0$ such that 
$$
\|f\|_{\Phi, r,\alpha; R}\le \|f\|'_{\Phi, r,\alpha; R}\le C\|f\|_{\Phi, r,\alpha; R},
\quad f\in{\mathcal L},\,\, R\in(0,\infty].
$$
Furthermore,
$$
\|f+g\|'_{\Phi, r,\alpha; R}\le \|f\|'_{\Phi, r,\alpha; R}+\|g\|'_{\Phi, r,\alpha; R},
\quad f, g\in{\mathcal L},\,\, R\in(0,\infty].
$$
\end{lemma}
\section{Proof of Theorem~\ref{Theorem:1.3}}
We consider case~(B), that is, 
\begin{equation}
\label{eq:4.1}
\frac{q+1}{pq-1}=\frac{N}{2}\quad\mbox{and}\quad p<q,
\end{equation}
and prove Theorem~\ref{Theorem:1.3} using uniformly local weak Zygmund type spaces 
$L_{{\rm ul}}^{r,\infty}(\log L)^\alpha$ and ${\mathfrak L}_{{\rm ul}}^{r,\infty}(\log {\mathfrak L})^\alpha$. 
Throughout this section, we set $\Phi(\tau)\coloneqq\log(e+\tau)$ for $\tau\ge 0$. 
Then $(\Phi1)$--$(\Phi3)$ hold. 

Recalling $pq>1$, we set 
\begin{equation}
\label{eq:4.2}
r\in\left(\frac{q+1}{p+1},q\right).
\end{equation}
Let $\alpha_*\in(0,\beta_B)$. 
It follows from \eqref{eq:4.1} that
\begin{gather}
\label{eq:4.3}
 pr>\frac{pq+p}{p+1}>1,\\
\label{eq:4.4}
 -\frac{N}{2}p+1=-\frac{N}{2}p+\frac{N}{2}\frac{pq-1}{q+1}=-\frac{N}{2}\frac{p+1}{q+1},\\
\label{eq:4.5}
p=\frac{p+1}{q+1}+\frac{2}{N}<1+\frac{2}{N}.
\end{gather}
Let $T_*\in(0,\infty)$. 
For any $T\in(0,T_*]$,
by Proposition~\ref{Proposition:3.2} and Lemma~\ref{Lemma:3.10} we find $C_*>0$ such that 
\begin{equation}
\label{eq:4.6}
\begin{split}
 & \|S(D_1t)\mu\|'_{\frac{q+1}{p+1},\alpha_B; T^{\frac{1}{2}}}\le C_*\|\mu\|_{\frac{q+1}{p+1},\alpha_B; T^{\frac{1}{2}}},\\
 & |||S(D_1t)\mu|||_{r,\alpha_*; T^{\frac{1}{2}}}
 \le C_*t^{-\frac{N}{2}\left(\frac{p+1}{q+1}-\frac{1}{r}\right)}\Phi(t^{-1})^{-p\beta_B+\frac{\alpha_*}{r}}\|\mu\|_{\frac{q+1}{p+1},\alpha_B; T^{\frac{1}{2}}},\\
 & \|S(D_1t)\mu\|_{L^\infty}
 \le C_*t^{-\frac{N}{2}{\frac{p+1}{q+1}}}\Phi(t^{-1})^{-p\beta_B}\|\mu\|_{\frac{q+1}{p+1},\alpha_B; T^{\frac{1}{2}}},\\
 & |||S(D_2t)\nu)|||_{1,\beta_B; T^{\frac{1}{2}}}\le C_*|||\nu|||_{1,\beta_B; T^{\frac{1}{2}}},\\
 & \|S(D_2t)\nu)\|_{L^\infty}\le C_*t^{-\frac{N}{2}}\Phi(t^{-1})^{-\beta_B}|||\nu|||_{1,\beta_B; T^{\frac{1}{2}}},\quad t\in(0,T),
\end{split}
\end{equation}
where $\alpha_B$ and $\beta_B$ are as in Theorem~\ref{Theorem:1.3}, that is, 
$\alpha_B=\frac{q+1}{p+1}\frac{p}{pq-1}$ and $\beta_B=\frac{1}{pq-1}$. 
Then we have:
\begin{lemma}
\label{Lemma:4.1}
Consider case~{\rm (B)}. 
Let $\{(u_n,v_n)\}$ be as in Section~{\rm 2.1}. 
Let $r$ and $\alpha_*$ be as in the above. 
Let  
\begin{equation}
\label{eq:4.7}
0<\epsilon<\frac{pq-1}{p}.
\end{equation}
For any $T_*\in(0,\infty)$, there exists $\delta>0$ with the following property: 
if 
\begin{equation}
\label{eq:4.8}
\|\mu\|_{\frac{q+1}{p+1},\alpha_B; T^{\frac{1}{2}}}\le\delta,\qquad
|||\nu|||_{1,\beta_B; T^{\frac{1}{2}}}\le\delta^{q-\epsilon},
\end{equation}
for some $T\in(0,T_*]$, 
then 
\begin{align}
\label{eq:4.9}
 & \sup_{t\in(0,T)}\|u_n(t)\|'_{\frac{q+1}{p+1},\alpha_B; T^{\frac{1}{2}}}\le 2C_*\delta,\\
\label{eq:4.10}
 & \sup_{t\in(0,T)}\left\{t^{\frac{N}{2}\left(\frac{p+1}{q+1}-\frac{1}{r}\right)}\Phi(t^{-1})^{p\beta_B-\frac{\alpha_*}{r}}
 |||u_n(t)|||_{r,\alpha_*; T^{\frac{1}{2}}}\right\}\le 2C_*\delta,\\
\label{eq:4.11}
 & \sup_{t\in(0,T)}\left\{t^{\frac{N}{2}\frac{p+1}{q+1}}\Phi(t^{-1})^{p\beta_B}\|u_n(t)\|_{L^\infty}\right\}\le 2C_*\delta,\\
\label{eq:4.12}
 & \sup_{t\in(0,T)}|||v_n(t)|||_{1,\beta_B; T^{\frac{1}{2}}}\le 2C_*\delta^{q-\epsilon},\\
 \label{eq:4.13} 
 & \sup_{t\in(0,T)} \left\{t^{\frac{N}{2}}\Phi(t^{-1})^{\beta_B}\|v_n(t)\|_{L^\infty}\right\}\le 2C_*\delta^{q-\epsilon},
\end{align}
for $n=0,1,2,\dots$, where $C_*$ is as in \eqref{eq:4.6}. 
Furthermore, there exists $C>0$ such that
\begin{align}
\label{eq:4.14}
 & \sup_{t\in(0,T)}\left\{\Phi(t^{-1})^{p\beta_B-\frac{p+1}{q+1}\alpha}
\left\|\,\int_0^t S(D_1(t-s))v_n(s)^p\,ds\,\right\|'_{\frac{q+1}{p+1},\alpha;T^{\frac{1}{2}}}\right\}\le C\delta^{(q-\epsilon)p},\\
\label{eq:4.15}
 & \sup_{t\in(0,T)}\left\{\Phi(t^{-1})^{\beta_B-\beta}\,\biggr|\biggr|\biggr|\int_0^t S(D_2(t-s))u_n(s)^q\,ds\,\biggr|\biggr|\biggr|_{1,\beta;T^{\frac{1}{2}}}\right\}\le C\delta^q,
 \end{align}
 for $\alpha\in[0,\alpha_B]$, $\beta\in[\alpha_*,\beta_B]$, and $n=0,1,2,\dots$. 
\end{lemma}
{\bf Proof.}
Let $T_*\in(0,\infty)$, and assume \eqref{eq:4.8} for some $T\in(0,T_*]$. 
By induction we prove \eqref{eq:4.9}--\eqref{eq:4.15} for $n=0,1,2,\dots$. 
It follows from \eqref{eq:4.6} that \eqref{eq:4.9}--\eqref{eq:4.13} hold for $n=0$. 
We assume that \eqref{eq:4.9}--\eqref{eq:4.13} hold for some $n=n_*\in\{0,1,2,\dots\}$. 
\vspace{3pt}
\newline
\underline{Step 1.} 
Let 
$$
\ell\in\left\{\frac{q+1}{p+1},r,\infty\right\},\quad \gamma\in[0,\infty). 
$$
Set 
\begin{equation*}
\|\,\cdot\,\|_{X_{\ell,\gamma}}\coloneqq
\left\{
\begin{array}{ll} 
\|\,\cdot\,\|'_{\frac{q+1}{p+1},\gamma;T^{\frac{1}{2}}}\quad & \mbox{if}\quad \ell=\frac{q+1}{p+1},\vspace{3pt}\\
|||\,\cdot\,|||_{r,\gamma;T^{\frac{1}{2}}} & \mbox{if}\quad \ell=r,\vspace{3pt}\\
\|\,\cdot\,\|_{L^\infty} & \mbox{if}\quad \ell=\infty,
\end{array}
\right.
\end{equation*}
for simplicity. We find $C_1=C_1(N,p,q,r)>0$ such that 
\begin{equation}
\label{eq:4.16}
\left\|\,\int_0^t S(D_1(t-s))v_{n_*}(s)^p\,ds\,\right\|_{X_{\ell,\gamma}}
\le C_1C_*^p\delta^{(q-\epsilon)p}t^{-\frac{N}{2}\left(\frac{p+1}{q+1}-\frac{1}{\ell}\right)}\Phi(t^{-1})^{\frac{\gamma}{\ell}-p\beta_B}
\end{equation}
for $t\in(0,T)$. 

If $p\ge 1$, 
thanks to \eqref{eq:4.4} and \eqref{eq:4.5}, 
by \eqref{eq:4.12} and \eqref{eq:4.13} with $n=n_*$ we apply Proposition~\ref{Proposition:3.2}
and Lemmas~\ref{Lemma:3.1}, \ref{Lemma:3.6}, and~\ref{Lemma:3.10} to obtain
\begin{equation}
\label{eq:4.17}
\begin{split}
 & \left\|\,\int_0^{t/2} S(D_1(t-s))v_{n_*}(s)^p\,ds\,\right\|_{X_{\ell,\gamma}}
 \le \int_0^{t/2} \|S(D_1(t-s))v_{n_*}(s)^p\|_{X_{\ell,\gamma}}\,ds\\
 & \le C\int_0^{t/2}(t-s)^{-\frac{N}{2}\left(1-\frac{1}{\ell}\right)}
 \Phi((t-s)^{-1})^{\frac{\gamma}{\ell}-\beta_B}
 |||v_{n_*}(s)^p|||_{1,\beta_B;T^{\frac{1}{2}}}\,ds\\
 & \le Ct^{-\frac{N}{2}\left(1-\frac{1}{\ell}\right)}\Phi(t^{-1})^{\frac{\gamma}{\ell}-\beta_B}
\int_0^{t/2}\|v_{n_*}(s)\|_{L^\infty}^{p-1} |||v_{n_*}(s)|||_{1,\beta_B;T^{\frac{1}{2}}}\,ds\\
 & \le CC_*^p\delta^{(q-\epsilon)p}t^{-\frac{N}{2}\left(1-\frac{1}{\ell}\right)}\Phi(t^{-1})^{\frac{\gamma}{\ell}-\beta_B}
 \int_0^{t/2} s^{-\frac{N(p-1)}{2}}\Phi(s^{-1})^{-(p-1)\beta_B}\,ds\\
 & \le CC_*^p\delta^{(q-\epsilon)p}t^{-\frac{N}{2}\left(1-\frac{1}{\ell}\right)}\Phi(t^{-1})^{\frac{\gamma}{\ell}-\beta_B}
 \cdot t^{-\frac{N}{2}(p-1)+1}\Phi(t^{-1})^{-(p-1)\beta_B}\\
 & =CC_*^p\delta^{(q-\epsilon)p}t^{-\frac{N}{2}\left(\frac{p+1}{q+1}-\frac{1}{\ell}\right)}\Phi(t^{-1})^{\frac{\gamma}{\ell}-p\beta_B},
 \quad t\in(0,T).
\end{split}
\end{equation}
Similarly, if $0< p<1$, thanks to \eqref{eq:4.4}, 
by \eqref{eq:4.12} with $n=n_*$
we apply Proposition~\ref{Proposition:3.2}, Lemma~\ref{Lemma:3.1}, and Lemma~\ref{Lemma:3.10} to obtain
\begin{equation}
\label{eq:4.18}
\begin{split}
 & \left\|\,\int_0^{t/2} S(D_1(t-s))v_{n_*}(s)^p\,ds\,\right\|_{X_{\ell,\gamma}}\le \int_0^{t/2} \|S(D_1(t-s))v_{n_*}(s)^p\|_{X_{\ell,\gamma}}\,ds\\
 & \le C\int_0^{t/2}(t-s)^{-\frac{N}{2}\left(p-\frac{1}{\ell}\right)}\Phi((t-s)^{-1})^{\frac{\gamma}{\ell}-p\beta_B}
 |||v_{n_*}(s)^p|||_{\frac{1}{p},\beta_B;T^{\frac{1}{2}}}\,ds\\
 & \le Ct^{-\frac{N}{2}\left(p-\frac{1}{\ell}\right)}\Phi(t^{-1})^{\frac{\gamma}{\ell}-p\beta_B}
\int_0^{t/2}|||v_{n_*}(s)|||_{1,\beta_B;T^{\frac{1}{2}}}^p\,ds\\
 & \le CC_*^p\delta^{(q-\epsilon)p}t^{-\frac{N}{2}\left(p-\frac{1}{\ell}\right)+1}\Phi(t^{-1})^{\frac{\gamma}{\ell}-p\beta_B}\\
 & =CC_*^p\delta^{(q-\epsilon)p}t^{-\frac{N}{2}\left(\frac{p+1}{q+1}-\frac{1}{\ell}\right)}\Phi(t^{-1})^{\frac{\gamma}{\ell}-p\beta_B},
 \quad t\in(0,T).
\end{split}
\end{equation}
On the other hand, 
by \eqref{eq:4.3} we find $\ell_*\in(1,\ell)$ such that 
$$
\frac{N}{2}\left(\frac{1}{\ell_*}-\frac{1}{\ell}\right)<1,\qquad p\ell_*>1.
$$
Then, thanks to \eqref{eq:4.4}, 
by \eqref{eq:4.12} and \eqref{eq:4.13} with $n=n_*$  we apply Proposition~\ref{Proposition:3.2}, Lemma~\ref{Lemma:3.6}, and Lemma~\ref{Lemma:3.10} to obtain
\begin{equation*}
\begin{split}
 & \left\|\,\int_{t/2}^t S(D_1(t-s))v_{n_*}(s)^p\,ds\,\right\|_{X_{\ell,\gamma}}\le \int_{t/2}^t \|S(D_1(t-s))v_{n_*}(s)^p\|_{X_{\ell,\gamma}}\,ds\\
 & \le C\int_{t/2}^t (t-s)^{-\frac{N}{2}\left(\frac{1}{\ell_*}-\frac{1}{\ell}\right)}
 \Phi((t-s)^{-1})^{\frac{\gamma}{\ell}-\frac{\beta_B}{\ell_*}}
 |||v_{n_*}(s)^p|||_{\ell_*,\beta_B;T^{\frac{1}{2}}}\,ds\\
 & \le C\int_{t/2}^t (t-s)^{-\frac{N}{2}\left(\frac{1}{\ell_*}-\frac{1}{\ell}\right)}
 \Phi((t-s)^{-1})^{\frac{\gamma}{\ell}-\frac{\beta_B}{\ell_*}}
 \|v_{n_*}(s)\|_{L^\infty}^{p-\frac{1}{\ell_*}}|||v_{n_*}(s)|||_{1,\beta_B;T^{\frac{1}{2}}}^{\frac{1}{\ell_*}}\,ds\\
 & \le CC_*^p\delta^{(q-\epsilon)p}\left(t^{-\frac{N}{2}}\Phi(t^{-1})^{-\beta_B}\right)^{p-\frac{1}{\ell_*}}\\
 & \qquad\quad
 \times\int_{t/2}^t (t-s)^{-\frac{N}{2}\left(\frac{1}{\ell_*}-\frac{1}{\ell}\right)}
 \Phi((t-s)^{-1})^{\frac{\gamma}{\ell}-\frac{\beta_B}{\ell_*}}\,ds\\
 & \le CC_*^p\delta^{(q-\epsilon)p}t^{-\frac{N}{2}\left(p-\frac{1}{\ell_*}\right)}\Phi(t^{-1})^{-\beta_B\left(p-\frac{1}{\ell_*}\right)}
 t^{-\frac{N}{2}\left(\frac{1}{\ell_*}-\frac{1}{\ell}\right)+1}\Phi(t^{-1})^{\frac{\gamma}{\ell}-\frac{\beta_B}{\ell_*}}\\
 & =CC_*^p\delta^{(q-\epsilon)p}t^{-\frac{N}{2}\left(\frac{p+1}{q+1}-\frac{1}{\ell}\right)}\Phi(t^{-1})^{\frac{\gamma}{\ell}-p\beta_B}, 
 \quad t\in(0,T). 
 \end{split}
\end{equation*}
This together with \eqref{eq:4.17} and \eqref{eq:4.18} implies \eqref{eq:4.16}. 
Furthermore, applying \eqref{eq:4.16} with $\ell=(q+1)/(p+1)$ 
and $\gamma=\alpha\in [0,\alpha_B]$, 
we obtain \eqref{eq:4.14} with $n=n_*$.
\vspace{5pt}
\newline
\underline{Step 2.} 
We prove that \eqref{eq:4.9}--\eqref{eq:4.11} hold with $n=n_*+1$. 
Let $\delta>0$ be small enough. 
By Lemma~\ref{Lemma:3.10}, \eqref{eq:4.6}, \eqref{eq:4.8}, and \eqref{eq:4.16} with $\ell=(q+1)/(p+1)$ and $\gamma=\alpha_B$ 
we have 
\begin{align*}
\|u_{n_*+1}(t)\|'_{\frac{q+1}{p+1},\alpha_B;T^{\frac{1}{2}}}
 & \le \|S(D_1t)\mu\|'_{\frac{q+1}{p+1},\alpha_B;T^{\frac{1}{2}}}+\left\|\,\int_0^t S(D_1(t-s))v_{n_*}(s)^p\,ds\,\right\|'_{\frac{q+1}{p+1},\alpha_B;T^{\frac{1}{2}}}\\
 & \le C_*\delta+CC_*^p\delta^{(q-\epsilon)p}\le 2C_*\delta,\quad t\in(0,T).
\end{align*}
Here we used the relations $\alpha_B(p+1)/(q+1)-p\beta_B=0$ (see \eqref{eq:1.13})
and  $(q-\epsilon)p>1$ (see~\eqref{eq:4.7}).
Similarly, by \eqref{eq:4.6}, \eqref{eq:4.8}, and \eqref{eq:4.16} with $\ell=r$ and $\gamma=\alpha_*$ 
we have 
\begin{align*}
|||u_{n_*+1}(t)|||_{r,\alpha_*;T^{\frac{1}{2}}}
 & \le |||S(D_1t)\mu|||_{r,\alpha_*;T^{\frac{1}{2}}}+\biggr|\biggr|\biggr|\,\int_0^t S(D_1(t-s))v_{n_*}(s)^p\,ds\,\biggr|\biggr|\biggr|_{r,\alpha_*;T^{\frac{1}{2}}}\\
 & \le \left(C_*\delta+CC_*^p\delta^{(q-\epsilon)p}\right)t^{-\frac{N}{2}\left(\frac{p+1}{q+1}-\frac{1}{r}\right)}\Phi(t^{-1})^{\frac{\alpha_*}{r}-p\beta_B}\\
 & \le 2C_*\delta t^{-\frac{N}{2}\left(\frac{p+1}{q+1}-\frac{1}{r}\right)}\Phi(t^{-1})^{\frac{\alpha_*}{r}-p\beta_B},\quad t\in(0,T).
\end{align*}
Furthermore, by \eqref{eq:4.6}, \eqref{eq:4.8}, and \eqref{eq:4.16} with $\ell=\infty$ 
we have
\begin{align*}
\|u_{n_*+1}(t)\|_{L^\infty}
 & \le\|S(D_1t)\mu\|_{L^\infty}+\left\|\,\int_0^t S(D_1(t-s))v_{n_*}(s)^p\,ds\,\right\|_{L^\infty}\\
 & \le \left(C_*\delta+CC_*^p\delta^{(q-\epsilon)p}\right)t^{-\frac{N}{2}\frac{p+1}{q+1}}\Phi(t^{-1})^{-p\beta_B}\\
 & \le 2C_*\delta t^{-\frac{N}{2}\frac{p+1}{q+1}}\Phi(t^{-1})^{-p\beta_B},\quad t\in(0,T).
\end{align*}
These imply that \eqref{eq:4.9}--\eqref{eq:4.11} hold with $n=n_*+1$. 
\vspace{5pt}
\newline
\underline{Step 3.} 
Let $m\in[1,\infty]$ and $\eta\in[0,\infty)$ be such that $\eta\ge\alpha_*$ if $m=1$. 
We find $C_2=C_2(N,p,q,r,\alpha_*)>0$ such that 
\begin{equation}
\label{eq:4.19}
\biggr|\biggr|\biggr|\,\int_0^t S(D_2(t-s))u_{n_*}(s)^q\,ds\,\biggr|\biggr|\biggr|_{m,\eta;T^{\frac{1}{2}}}
\le C_2C_*^q\delta^qt^{-\frac{N}{2}\left(1-\frac{1}{m}\right)}\Phi(t^{-1})^{\frac{\eta}{m}-\beta_B}
\end{equation}
for $t\in(0,T)$. 
Set $m_*\coloneqq 1$ if $m=1$. 
If $m>1$, let $m_*\in[1,m)$ be such that 
$$
\frac{N}{2}\left(\frac{1}{m_*}-\frac{1}{m}\right)<1. 
$$
By Proposition~\ref{Proposition:3.2}, Lemma~\ref{Lemma:3.1}, and \eqref{eq:4.2} we have
\begin{align*}
 & \biggr|\biggr|\biggr|\,\int_0^t S(D_2(t-s))u_{n_*}(s)^q\,ds\,\biggr|\biggr|\biggr|_{m,\eta; T^{\frac{1}{2}}}
 \le\int_0^t |||S(D_2(t-s))u_{n_*}(s)^q|||_{m,\eta; T^{\frac{1}{2}}}\,ds\\
 & \le C\int_0^{t/2}(t-s)^{-\frac{N}{2}\left(1-\frac{1}{m}\right)}\Phi((t-s)^{-1})^{\frac{\eta}{m}-\alpha_*}|||u_{n_*}(s)^q|||_{1,\alpha_*; T^{\frac{1}{2}}}\,ds\\
 & \qquad\quad
 +C\int_{t/2}^t(t-s)^{-\frac{N}{2}\left(\frac{1}{m_*}-\frac{1}{m}\right)}\Phi((t-s)^{-1})^{\frac{\eta}{m}-\frac{\alpha_*}{m_*}}|||u_{n_*}(s)^q|||_{m_*,\alpha_*;T^{\frac{1}{2}}}\,ds\\
 & \le Ct^{-\frac{N}{2}\left(1-\frac{1}{m}\right)}\Phi(t^{-1})^{\frac{\eta}{m}-\alpha_*}
 \int_0^{t/2}\|u_{n_*}(s)\|_{L^\infty}^{q-r}|||u_{n_*}(s)|||^r_{r,\alpha_*; T^{\frac{1}{2}}}\,ds\\
 & \qquad
 +C\sup_{s\in(t/2,t)}\left\{\|u_{n_*}(s)\|_{L^\infty}^{q-\frac{r}{m_*}}|||u_{n_*}(s)|||_{r,\alpha_*; T^{\frac{1}{2}}}^{\frac{r}{m_*}}\right\}\\
 & \qquad\qquad
 \times\int_{t/2}^t(t-s)^{-\frac{N}{2}\left(\frac{1}{m_*}-\frac{1}{m}\right)}\Phi((t-s)^{-1})^{\frac{\eta}{m}-\frac{\alpha_*}{m_*}}\,ds,\quad t\in(0,T).
 \end{align*}
Furthermore, 
since $\Phi(\tau)=\log(e+\tau)$, $0<T\le T_*<\infty$, and $\alpha_*<\beta_B$, 
by \eqref{eq:4.10} and \eqref{eq:4.11} with $n=n_*$ we obtain
\begin{align*}
 & \int_0^{t/2}\|u_{n_*}(s)\|_{L^\infty}^{q-r}|||u_{n_*}(s)|||^r_{r,\alpha_*;T^{\frac{1}{2}}}\,ds\\
 & \le CC_*^q\delta^q\int_0^{t/2} \left(s^{-\frac{N}{2}\frac{p+1}{q+1}}\Phi(s^{-1})^{-p\beta_B}\right)^{q-r}
 \left(s^{-\frac{N}{2}\left(\frac{p+1}{q+1}-\frac{1}{r}\right)}\Phi(s^{-1})^{-p\beta_B+\frac{\alpha_*}{r}}\right)^r\,ds\\
 & \le CC_*^q\delta^q\int_0^{t/2} s^{-1}\Phi(s^{-1})^{-1-\beta_B+\alpha_*}\,ds
 \le CC_*^q\delta^q\Phi(t^{-1})^{-\beta_B+\alpha_*},\quad t\in(0,T).
\end{align*}
Here we used relations
\begin{equation}
\label{eq:4.20}
\begin{split}
 & -\frac{N}{2}\frac{p+1}{q+1}(q-r)-\frac{N}{2}\left(\frac{p+1}{q+1}-\frac{1}{r}\right)r
=-\frac{N}{2}\frac{pq+q}{q+1}+\frac{N}{2}=-\frac{N}{2}\frac{pq-1}{q+1}=-1,\\
 & -p\beta_B(q-r)-p\beta_B r+\alpha_*=-pq\beta_B+\alpha_*=-\frac{pq}{pq-1}+\alpha_*=-1-\beta_B+\alpha_*.
\end{split}
\end{equation}
The first relation (resp.~the second relation) follows from \eqref{eq:4.1} (resp.~\eqref{eq:1.13}).
Similarly, we see that
\begin{align*}
 & \|u_{n_*}(t)\|_{L^\infty}^{q-\frac{r}{m_*}}|||u_{n_*}(t)|||_{r,\alpha_*; T^{\frac{1}{2}}}^{\frac{r}{m_*}}\\
 & \le CC_*^q\delta^q\left(t^{-\frac{N}{2}\frac{p+1}{q+1}}\Phi(t^{-1})^{-p\beta_B}\right)^{q-\frac{r}{m_*}}
 \left(t^{-\frac{N}{2}\left(\frac{p+1}{q+1}-\frac{1}{r}\right)}\Phi(t^{-1})^{-p\beta_B+\frac{\alpha_*}{r}}\right)^{\frac{r}{m_*}}\\
 & =CC_*^q\delta^qt^{-\frac{N}{2}-1+\frac{N}{2m_*}}\Phi(t^{-1})^{-1-\beta_B+\frac{\alpha_*}{m_*}},
 \quad t\in(0,T).
\end{align*}
Here we also used relations 
\begin{align*}
 & -\frac{N}{2}\frac{p+1}{q+1}\left(q-\frac{r}{m_*}\right)-\frac{N}{2}\left(\frac{p+1}{q+1}-\frac{1}{r}\right)\frac{r}{m_*}\\
 & \qquad\quad
 =-\frac{N}{2}\frac{pq+q}{q+1}+\frac{N}{2m_*}
 =-\frac{N}{2}-\frac{N}{2}\frac{pq-1}{q+1}+\frac{N}{2m_*}
=-\frac{N}{2}-1+\frac{N}{2m_*},\\
 & -p\beta_B\left(q-\frac{r}{m_*}\right)+\left(-p\beta_B+\frac{\alpha_*}{r}\right)\frac{r}{m_*}
=-pq\beta_B+\frac{\alpha_*}{m_*}
=-\frac{pq}{pq-1}+\frac{\alpha_*}{m_*}
=-1-\beta_B+\frac{\alpha_*}{m_*}.
\end{align*}
Similarly to \eqref{eq:4.20}, the first relation (resp.~the second relation) follows from \eqref{eq:4.1} (resp.~\eqref{eq:1.13}). 
These together with Lemma~\ref{Lemma:3.6}~(1) imply that 
\begin{align*}
 & \biggr|\biggr|\biggr|\,\int_0^t S(D_2(t-s))u_{n_*}(s)^q\,ds\,\biggr|\biggr|\biggr|_{m,\eta;T^{\frac{1}{2}}}\\
 & \le CC_*^q\delta^qt^{-\frac{N}{2}\left(1-\frac{1}{m}\right)}\Phi(t^{-1})^{\frac{\eta}{m}-\alpha_*}
 \cdot\Phi(t^{-1})^{-\beta_B+\alpha_*}\\
 & \qquad
 +CC_*^q\delta^qt^{-\frac{N}{2}-1+\frac{N}{2m_*}}\Phi(t^{-1})^{-1-\beta_B+\frac{\alpha_*}{m_*}}
 \cdot t^{-\frac{N}{2}\left(\frac{1}{m_*}-\frac{1}{m}\right)+1}\Phi(t^{-1})^{\frac{\eta}{m}-\frac{\alpha_*}{m_*}}\\
 & \le CC_*^q\delta^qt^{-\frac{N}{2}\left(1-\frac{1}{m}\right)}\Phi(t^{-1})^{\frac{\eta}{m}-\beta_B},\quad t\in(0,T).
\end{align*}
This implies \eqref{eq:4.19}. 
Furthermore, applying \eqref{eq:4.19} with $m=1$ and $\eta=\beta\in[\alpha_*,\beta_B]$, 
we obtain \eqref{eq:4.15} with $n=n_*$. 
\vspace{5pt}
\newline
\underline{Step 4.} 
We prove that \eqref{eq:4.12} and \eqref{eq:4.13} hold for $n=n_*+1$.
Taking small enough $\delta>0$ if necessary, 
by \eqref{eq:4.6}, \eqref{eq:4.8}, and \eqref{eq:4.19} with $m=1$ and $\eta=\beta_B$ 
we have 
\begin{align*}
|||v_{n_*+1}(t)|||_{1,\beta_B; T^{\frac{1}{2}}}
 & \le |||S(D_2t)\nu|||_{1,\beta_B; T^{\frac{1}{2}}}+\biggr|\biggr|\biggr|\,\int_0^t S(D_2(t-s))u_{n_*}(s)^q\,ds\,\biggr|\biggr|\biggr|_{1,\beta_B; T^{\frac{1}{2}}}\\
 & \le C_*\delta^{q-\epsilon}+CC_*^q\delta^q\le 2C_*\delta^{q-\epsilon},\quad t\in(0,T).
\end{align*}
Similarly, by \eqref{eq:4.6}, \eqref{eq:4.8}, and \eqref{eq:4.19} with $m=\infty$
we obtain 
\begin{align*}
 & \|v_{n_*+1}(t)\|_{L^\infty}
\le\|S(D_2t)\nu\|_{L^\infty}+\left\|\,\int_0^t S(D_2(t-s))u_{n_*}(s)^q\,ds\,\right\|_{L^\infty}\\
 & \qquad\quad 
 \le (C_*\delta^{q-\epsilon}+CC_*^q\delta^q)t^{-\frac{N}{2}}\Phi(t^{-1})^{-\beta_B}
 \le 2C_*\delta^{q-\epsilon} t^{-\frac{N}{2}}\Phi(t^{-1})^{-\beta_B},\quad t\in(0,T).
\end{align*}
These imply that \eqref{eq:4.12} and \eqref{eq:4.13} hold for $n=n_*+1$.
Therefore \eqref{eq:4.9}--\eqref{eq:4.15} hold for $n=0,1,2,\dots$, 
and the proof of Lemma~\ref{Lemma:4.1} is complete.
$\Box$
\vspace{5pt}
\newline
{\bf Proof of Theorem~\ref{Theorem:1.3}.}
Let $T_*>0$, $\epsilon>0$, and $\delta>0$ be as in Lemma~\ref{Lemma:4.1}. 
Let $\delta_B>0$ be such that $\delta_B\le\min\{\delta,\delta^{q-\epsilon}\}$, 
and assume \eqref{eq:1.14}. 
Let $\{(u_n,v_n)\}$ be as in \eqref{eq:2.1}, and define the limit function $(u,v)$ of $\{(u_n,v_n)\}$ by \eqref{eq:2.3}. 
Then we apply the arguments in Section~2.1 together with Lemma~\ref{Lemma:4.1} to see 
that $(u,v)$ is a solution to problem~\eqref{eq:P} in ${\mathbb R}^N\times(0,T)$ satisfying \eqref{eq:4.9}--\eqref{eq:4.15} with $(u_n,v_n)$ replaced by $(u,v)$. 
Furthermore, we deduce from \eqref{eq:3.8} 
that $(u,v)$ satisfies \eqref{eq:1.15} and \eqref{eq:1.16}. 
Thus Theorem~\ref{Theorem:1.3} follows.
$\Box$
\section{Proof of Theorem~\ref{Theorem:1.4}}
In this section we consider case~(C), that is, 
$$
 p=q=1+\frac{2}{N}.
$$
Similarly to Section~4, throughout this section, we set $\Phi(\tau)\coloneqq\log(e+\tau)$ for $\tau\ge 0$. 

Let $0\le \gamma_*<N/2$ and $T_*>0$. 
For any $T\in(0,T_*]$, by Proposition~\ref{Proposition:3.2} we find $C_*>0$ such that 
\begin{equation}
\label{eq:5.1}
\begin{split}
 & \sup_{0<t<T}\,\left\{|||S(D_1t)\mu|||_{1,\frac{N}{2}; T^{\frac{1}{2}}}+|||S(D_2t)\nu|||_{1,\frac{N}{2}; T^{\frac{1}{2}}}\right\}\le C_*\Lambda,\\
 & \sup_{0<t<T}\,\left\{t^{\frac{N}{2}\left(1-\frac{1}{p}\right)}\Phi(t^{-1})^{-\frac{\gamma_*}{p}+\frac{N}{2}}
 \left(|||S(D_1t)\mu|||_{p,\gamma_*; T^{\frac{1}{2}}}+|||S(D_2t)\nu|||_{p,\gamma_*; T^{\frac{1}{2}}}\right)\right\} \le C_*\Lambda,\\
 & \sup_{0<t<T}\,\left\{t^{\frac{N}{2}}\Phi(t^{-1})^{\frac{N}{2}}
  \left(\|S(D_1t)\mu\|_{L^\infty}+\|S(D_2t)\nu\|_{L^\infty}\right)\right\}
  \le C_*\Lambda,
\end{split}
\end{equation}
where $\Lambda\coloneqq|||\mu|||_{1,\frac{N}{2};T^{\frac{1}{2}}}+|||\nu|||_{1,\frac{N}{2};T^{\frac{1}{2}}}$.
Then we have:
\begin{lemma}
\label{Lemma:5.1}
Consider case~{\rm (C)}. 
Let $\{(u_n,v_n)\}$ be as in Section~{\rm 2.1}. 
Let $T_*$ and $\gamma_*$ be as in the above. 
Then there exists $\delta>0$ with the following properties: 
if $(\mu,\nu)$ satisfies 
\begin{equation}
\label{eq:5.2}
|||\mu|||_{1,\frac{N}{2};T^{\frac{1}{2}}}+|||\nu|||_{1,\frac{N}{2};T^{\frac{1}{2}}}\le\delta
\end{equation}
for some $T\in(0,T_*]$, then
\begin{align}
\label{eq:5.3}
 & \sup_{0<t<T}\,\left\{|||u_n(t)|||_{1,\frac{N}{2};T^{\frac{1}{2}}}+|||v_n(t)|||_{1,\frac{N}{2};T^{\frac{1}{2}}}\right\}\le 2C_*\delta,\\
\label{eq:5.4}
 & \sup_{0<t<T}\,\left\{t^{\frac{N}{2}\left(1-\frac{1}{p}\right)}\Phi(t^{-1})^{-\frac{\gamma_*}{p}+\frac{N}{2}}
 \left(|||u_n(t)|||_{p,\gamma_*;T^{\frac{1}{2}}}+|||v_n(t)|||_{p,\gamma_*;T^{\frac{1}{2}}}\right)\right\}\le 2C_*\delta,\\
\label{eq:5.5}
 & \sup_{0<t<T}\,\left\{t^{\frac{N}{2}}\Phi(t^{-1})^{\frac{N}{2}}
 \left(\|u_n(t)\|{L^\infty}+\|v_n(t)\|_{L^\infty}\right)\right\}\le 2C_*\delta,
\end{align}
for $n=0,1,2,\dots$, where $C_*$ is as in \eqref{eq:5.1}. 
Furthermore, for any $\eta\in[\gamma_*,N/2]$, 
there exists $C>0$ such that 
\begin{equation}
\label{eq:5.6}
\begin{split}
 & \biggr|\biggr|\biggr|\int_0^t S(D_2(t-s))u_n(s)^p\,ds\,\biggr|\biggr|\biggr|_{1,\eta;T^{\frac{1}{2}}}
 +\biggr|\biggr|\biggr|\int_0^t S(D_1(t-s))v_n(s)^p\,ds\,\biggr|\biggr|\biggr|_{1,\eta;T^{\frac{1}{2}}}\\
 & \le C\Phi(t^{-1})^{\eta-\frac{N}{2}}
 \end{split}
 \end{equation}
 for $t\in(0,T)$ and $n=0,1,2,\dots$.
\end{lemma}
{\bf Proof.}
Let $T_*\in(0,\infty)$, and assume \eqref{eq:5.2} for some $T\in(0,T_*]$. 
By induction we prove \eqref{eq:5.3}--\eqref{eq:5.6}. 
It follows from \eqref{eq:5.1} that \eqref{eq:5.3}--\eqref{eq:5.5} hold for $n=0$. 
We assume that \eqref{eq:5.3}--\eqref{eq:5.5} hold for some $n=n_*\in\{0,1,2,\dots\}$. 

Let $\ell\in[1,\infty]$ and $\eta\in[\gamma_*,N/2]$. Set $\ell_*\coloneqq1$ if $\ell=1$. 
If $\ell>1$, let $\ell_*\in(1,\ell)$ be such that  
\begin{equation}
\label{eq:5.7}
\frac{N}{2}\left(\frac{1}{\ell_*}-\frac{1}{\ell}\right)<1.
\end{equation}
By Proposition~\ref{Proposition:3.2} 
we obtain 
\begin{equation*}
\begin{split}
 & \biggr|\biggr|\biggr|\int_0^t S(D_2(t-s))u_{n_*}(s)^p\,ds\,\biggr|\biggr|\biggr|_{\ell,\eta;T^{\frac{1}{2}}}
  +\,\biggr|\biggr|\biggr|\int_0^t S(D_1(t-s))v_n(s)^p\,ds\,\biggr|\biggr|\biggr|_{\ell,\eta;T^{\frac{1}{2}}}\\
 & \le C\int_0^{t/2} (t-s)^{-\frac{N}{2}\left(1-\frac{1}{\ell}\right)}\Phi((t-s)^{-1})^{\frac{\eta}{\ell}-\gamma_*}
 \left(|||u_{n_*}(s)^p|||_{1,\gamma_*;T^{\frac{1}{2}}}+|||v_{n_*}(s)^p|||_{1,\gamma_*;T^{\frac{1}{2}}}\right)\,ds\\
 &+C\int_{t/2}^t (t-s)^{-\frac{N}{2}\left(\frac{1}{\ell_*}-\frac{1}{\ell}\right)}\Phi((t-s)^{-1})^{\frac{\eta}{\ell}-\frac{\gamma_*}{\ell_*}}
  \left(|||u_{n_*}(s)^p|||_{\ell_*,\gamma_*;T^{\frac{1}{2}}}+|||v_{n_*}(s)^p|||_{\ell_*,\gamma_*;T^{\frac{1}{2}}}\right)\,ds
\end{split}
\end{equation*}
for $t\in(0,T)$. 
On the other hand, 
by Lemma~\ref{Lemma:3.1}, \eqref{eq:5.4} with $n=n_*$, and \eqref{eq:5.5} with $n=n_*$ we have 
\begin{equation*}
\begin{split}
 &  |||u_{n_*}(t)^p|||_{1,\gamma_*;T^{\frac{1}{2}}}+|||v_{n_*}(t)^p|||_{1,\gamma_*;T^{\frac{1}{2}}}
=|||u_{n_*}(s)|||^p_{p,\gamma_*;T^{\frac{1}{2}}}+|||v_{n_*}(s)|||^p_{p,\gamma_*;T^{\frac{1}{2}}}\\
 & \le CC_*^p\delta^p\, t^{-\frac{N}{2}(p-1)}\Phi(t^{-1})^{\gamma_*-\frac{N}{2}p}
 =CC_*^p\delta^p\, t^{-1}\Phi(t^{-1})^{\gamma_*-\frac{N}{2}p}
\end{split}
\end{equation*}
and
\begin{equation*}
\begin{split}
 &  |||u_{n_*}(t)^p|||_{\ell_*,\gamma_*;T^{\frac{1}{2}}}+|||v_{n_*}(t)^p|||_{\ell_*,\gamma_*;T^{\frac{1}{2}}}
=|||u_{n_*}(s)^{p\ell_*}|||^{\frac{1}{\ell_*}}_{1,\gamma_*;T^{\frac{1}{2}}}+|||v_{n_*}(s)^{p\ell_*}|||^{\frac{1}{\ell_*}}_{1,\gamma_*;T^{\frac{1}{2}}}\\
 & \le\|u_{n_*}(s)\|_{L^\infty}^{p\left(1-\frac{1}{\ell_*}\right)}|||u_{n_*}(s)^p|||^{\frac{1}{\ell_*}}_{1,\gamma_*;T^{\frac{1}{2}}}
+\|v_{n_*}(s)\|_{L^\infty}^{p\left(1-\frac{1}{\ell_*}\right)}|||v_{n_*}(s)^p|||^{\frac{1}{\ell_*}}_{1,\gamma_*;T^{\frac{1}{2}}}\\\
 & \le CC_*^p\delta^p\biggr\{t^{-\frac{N}{2}}\Phi(t^{-1})^{-\frac{N}{2}}\biggr\}^{p\left(1-\frac{1}{\ell_*}\right)}
 \biggr\{t^{-1}\Phi(t^{-1})^{\gamma_*-\frac{N}{2}p}\biggr\}^{\frac{1}{\ell_*}}\\
 & \le CC_*^p\delta^p\, t^{-\frac{N}{2}\left(p-\frac{1}{\ell_*}\right)}\Phi(t^{-1})^{-\frac{N}{2}p+\frac{\gamma_*}{\ell_*}}
\end{split}
\end{equation*}
for $t\in(0,T)$. 
Since $T_*<\infty$ and 
$$
\gamma_*-\frac{N}{2}p=\gamma_*-\frac{N}{2}\left(1+\frac{2}{N}\right)=\gamma_*-\frac{N}{2}-1<-1,
$$
we deduce that
\begin{equation}
\label{eq:5.8}
\begin{split}
 & \biggr|\biggr|\biggr|\int_0^t S(D_2(t-s))u_{n_*}(s)^p\,ds\,\biggr|\biggr|\biggr|_{\ell,\eta;T^{\frac{1}{2}}}
  +\biggr|\biggr|\biggr|\int_0^t S(D_1(t-s))v_n(s)^p\,ds\,\biggr|\biggr|\biggr|_{\ell,\eta;T^{\frac{1}{2}}}\\
 & \le CC_*^p\delta^p\,t^{-\frac{N}{2}\left(1-\frac{1}{\ell}\right)}\Phi(t^{-1})^{\frac{\eta}{\ell}-\gamma_*}
 \int_0^{t/2}s^{-1}\Phi(s^{-1})^{\gamma_*-\frac{N}{2}p}\,ds\\
 & \qquad
 +CC_*^p\delta^p\,t^{-\frac{N}{2}\left(p-\frac{1}{\ell_*}\right)}\Phi(t^{-1})^{-\frac{N}{2}p+\frac{\gamma_*}{\ell_*}}
 \int_{t/2}^t (t-s)^{-\frac{N}{2}\left(\frac{1}{\ell_*}-\frac{1}{\ell}\right)}\Phi((t-s)^{-1})^{\frac{\eta}{\ell}-\frac{\gamma_*}{\ell_*}}\,ds\\
 & \le CC_*^p\delta^p\,t^{-\frac{N}{2}\left(1-\frac{1}{\ell}\right)}\Phi(t^{-1})^{-\frac{N}{2}+\frac{\eta}{\ell}}
 +CC_*^p\delta^pt^{-\frac{N}{2}\left(1-\frac{1}{\ell}\right)}\Phi(t^{-1})^{-\frac{N}{2}+\frac{\eta}{\ell}-1}\\
 & \le CC_*^p\delta^p\,t^{-\frac{N}{2}\left(1-\frac{1}{\ell}\right)}\Phi(t^{-1})^{-\frac{N}{2}+\frac{\eta}{\ell}}, \quad t\in(0,T).
\end{split}
\end{equation}
Here we used $\Phi(\tau)=\log(e+\tau)$ (resp.~Lemma~\ref{Lemma:3.6} and \eqref{eq:5.7}) 
in the estimate of the above integral on the interval $(0,t/2)$ (resp.~$(t/2,t)$). 
Then, by \eqref{eq:5.8} with $\ell=1$ we obtain \eqref{eq:5.6}. 
Furthermore, taking small enough $\delta>0$ if necessary, 
by \eqref{eq:5.1} and \eqref{eq:5.8} with $\ell=1$ and $\eta=N/2$ 
we see that
$$
\sup_{0<t<T}\,\left\{|||u_{n_*+1}(t)|||_{1,\frac{N}{2};T^{\frac{1}{2}}}+|||v_{n_*+1}(t)|||_{1,\frac{N}{2};T^{\frac{1}{2}}}\right\}\le C_*\delta+CC_*^p\delta^p\le 2C_*\delta.
$$
Thus \eqref{eq:5.3} holds with $n=n_*+1$. 
Similarly, taking small enough $\delta>0$ if necessary, 
by \eqref{eq:5.1} and \eqref{eq:5.8} with $\ell=p$ and $\eta=\gamma_*$ (resp.~with $\ell=\infty$) 
we obtain \eqref{eq:5.4} (resp.~\eqref{eq:5.5}) with $n=n_*+1$. 
Therefore we see that \eqref{eq:5.3}--\eqref{eq:5.6} hold for $n=0,1,2,\dots$, 
and the proof of Lemma~\ref{Lemma:5.1} is complete.
$\Box$\vspace{5pt}

\noindent
{\bf Proof of Theorem~\ref{Theorem:1.4}.}
Let $T_*\in(0,\infty)$. 
Let $\delta_C>0$ be small enough, and assume \eqref{eq:1.17} for some $T\in(0,T_*]$. 
Then, similarly to the proof of Theorem~\ref{Theorem:1.3}, 
by Lemma~\ref{Lemma:5.1} we find a solution~$(u,v)$ to problem~\eqref{eq:P} in ${\mathbb R}^N\times(0,T)$. 
Furthermore, the solution~$(u,v)$ satisfies \eqref{eq:5.3}--\eqref{eq:5.6} with $(u_n,v_n)$ replaced by $(u,v)$. 
Then, thanks to \eqref{eq:3.8}, $(u,v)$ is the desired solution. 
The proof of Theorem~\ref{Theorem:1.6} is complete. 
$\Box$
\section{Proofs of Theorems~\ref{Theorem:1.5} and \ref{Theorem:1.6}}
In this section we consider cases (D) and (E), that is, 
$$
\frac{q+1}{pq-1}>\frac{N}{2}\quad\mbox{and}\quad q\ge 1+\frac{2}{N}.
$$
Then 
\begin{equation}
\label{eq:6.1}
pq<1+\frac{2}{N}(q+1)\le q\left(1+\frac{2}{N}\right),\quad
\mbox{that is,}
\quad
p<1+\frac{2}{N}.
\end{equation}
Furthermore, it follows that 
\begin{equation}
\label{eq:6.2}
\delta\coloneqq-\frac{N}{2}\max\left\{p-\frac{N+2}{Nq},0\right\}+1>0,
\end{equation}
since 
$$
-\frac{N}{2}\left(p-\frac{N+2}{Nq}\right)+1
=-\frac{1}{q}\left(\frac{N}{2}(pq-1)-1\right)+1
>-\frac{1}{q}\left((q+1)-1\right)+1=0.
$$
Set 
$$
r_*\coloneqq\max\left\{\frac{Nq}{N+2},\frac{1}{p}\right\}\ge\frac{Nq}{N+2}\ge 1.
$$
Let $r^*\in(r_*,q)$. 
Assume that $(\mu,\nu)$ satisfies \eqref{eq:1.19} in case (D) 
(resp.~\eqref{eq:1.20} in case (E)).
By Proposition~\ref{Proposition:3.2}, \eqref{eq:1.12}, and \eqref{eq:2.7} we find $C_*>0$ such that 
\begin{equation}
\label{eq:6.3}
\begin{split}
 & \|S(D_1t)\mu\|_{L^{r^*}_{{\rm ul}}}
 \le C_*t^{-\frac{N}{2}\left(\frac{N+2}{Nq}-\frac{1}{r^*}\right)}\Phi(t^{-1})^{-1},\\
 & \|S(D_1t)\mu\|_{L^\infty}
 \le C_*t^{-\frac{N+2}{2q}}\Phi(t^{-1})^{-1},\\
 & \|S(D_2t)\nu\|_{L^1_{{\rm ul}}}\le C_*,\\
 & \|S(D_2t)\nu\|_{L^\infty}\le C_*t^{-\frac{N}{2}},
\end{split}
\end{equation}
for $t\in(0,1)$. 
Then we have:
\begin{lemma}
\label{Lemma:6.1}
Consider case~{\rm (D)} {\rm (}resp.~case~{\rm (E)}{\rm )}. 
Let $\{(u_n,v_n)\}$ be as in Section~{\rm 2.1}. 
Let $\Phi$ be a non-decreasing function in $[0,\infty)$ with properties~{\rm ($\Phi$1)}--{\rm ($\Phi$3)} and satisfy \eqref{eq:1.18}. 
Let $r^*$ and $C_*$ be as in the above. 
Assume that $(\mu,\nu)$ satisfies \eqref{eq:1.19} {\rm (}resp.~\eqref{eq:1.20}{\rm )}.
Then there exists $T\in(0,1)$ such that
\begin{align}
\label{eq:6.4}
 & \sup_{0<t<T}\left\{t^{\frac{N}{2}\left(\frac{N+2}{Nq}-\frac{1}{r^*}\right)}\Phi(t^{-1})\|u_n(t)\|_{L^{r^*}_{{\rm ul}}}\right\}\le 2C_*,\\
\label{eq:6.5}
 & \sup_{0<t<T}\left\{t^{\frac{N+2}{2q}}\Phi(t^{-1})\|u_n(t)\|_{L^\infty}\right\}\le 2C_*,\\
\label{eq:6.6}
 & \sup_{0<t<T}\|v_n(t)\|_{L^1_{{\rm ul}}}\le 2C_*,\\
 \label{eq:6.7} 
 & \sup_{0<t<T}\left\{t^{\frac{N}{2}}\|v_n(t)\|_{L^\infty}\right\}\le 2C_*,
\end{align}
for $n=0,1,2,\dots$. 
Furthermore, there exists $C>0$ such that 
\begin{align}
\label{eq:6.8}
 & \biggr|\biggr|\biggr|\int_0^t S(D_1(t-s))v_n(s)^p\,ds\,\biggr|\biggr|\biggr|_{\Phi,r_*,r_*}\le Ct^\delta\Phi(t^{-1}),\\
\label{eq:6.9}
 & \left\|\,\int_0^t S(D_2(t-s))u_n(s)^q\,ds\,\right\|_{L^1_{{\rm ul}}}\le C\int_0^t s^{-1}\Phi(s^{-1})^{-q}\,ds,
 \end{align}
 for $t\in(0,T)$ and $n=0,1,2,\dots$, where $\delta$ is as in \eqref{eq:6.2}. 
\end{lemma}
{\bf Proof of Lemma~\ref{Lemma:6.1}.}
By induction we obtain \eqref{eq:6.4}--\eqref{eq:6.9}. 
Let $T\in (0,1)$ be a constant to be chosen later. 
It follows from \eqref{eq:6.3} that \eqref{eq:6.4}--\eqref{eq:6.7} hold for $n=0$. 
We assume that \eqref{eq:6.4}--\eqref{eq:6.7} hold for some $n=n_*\in\{0,1,2,\dots\}$. 
Then, for any $\ell\in[1,\infty]$ with $\ell p\ge 1$, 
by Lemma~\ref{Lemma:3.1} we have 
\begin{equation}
\label{eq:6.10}
\|v_{n_*}(t)^p\|_{L^\ell_{{\rm ul}}}\le \|v_{n_*}(t)\|_{L^\infty}^{p-\frac{1}{\ell}}\|v_{n_*}(t)\|_{L^1_{{\rm ul}}}^{\frac{1}{\ell}}
\le (2C_*)^pt^{-\frac{N}{2}\left(p-\frac{1}{\ell}\right)},\quad t\in (0,T).
\end{equation}
%
\underline{Step 1.} 
We prove that \eqref{eq:6.8} holds for $n=n_*$ in the case of $r_*=Nq/(N+2)>1/p$. 
Let $r>1$ be such that 
$$
\frac{1}{p}<r<r_*=\frac{Nq}{N+2},\quad \frac{N}{2}\left(\frac{1}{r}-\frac{N+2}{Nq}\right)<1. 
$$
By Proposition~\ref{Proposition:3.2}, Lemma~\ref{Lemma:3.5}, Lemma~\ref{Lemma:3.6}, and \eqref{eq:6.10} we have
\begin{equation}
\label{eq:6.11}
\begin{split}
 & \biggr|\biggr|\biggr|\int_{t/2}^t S(D_1(t-s))v_{n_*}(s)^p\,ds\, \biggr|\biggr|\biggr|_{\Phi,r_*,r_*}\\
 & \le C\int_{t/2}^t(t-s)^{-\frac{N}{2}\left(\frac{1}{r}-\frac{N+2}{Nq}\right)}\Phi((t-s)^{-1})\|v_{n_*}(s)^p\|_{L^r_{{\rm ul}}}\,ds\\
 & \le C\int_{t/2}^t(t-s)^{-\frac{N}{2}\left(\frac{1}{r}-\frac{N+2}{Nq}\right)}\Phi((t-s)^{-1})s^{-\frac{N}{2}\left(p-\frac{1}{r}\right)}\,ds\\
 & \le Ct^{-\frac{N}{2}\left(p-\frac{N+2}{Nq}\right)+1}\Phi((t/2)^{-1})
 \le Ct^{-\delta}\Phi(t^{-1}),\quad t\in(0,T).
\end{split}
\end{equation}
On the other hand, 
if $p>1$, by Proposition~\ref{Proposition:3.2}, Lemma~\ref{Lemma:3.5}, \eqref{eq:6.1}, and \eqref{eq:6.10} we have 
\begin{equation}
\label{eq:6.12}
\begin{split}
 & \biggr|\biggr|\biggr|\int_0^{t/2} S(D_1(t-s))v_{n_*}(s)^p\,ds\,\biggr|\biggr|\biggr|_{\Phi,r_*,r_*}\\
 & \le C\int_0^{t/2}(t-s)^{-\frac{N}{2}\left(1-\frac{N+2}{Nq}\right)}\Phi((t-s)^{-1})\|v_{n_*}(s)^p\|_{L^1_{{\rm ul}}}\,ds\\
 & \le Ct^{-\frac{N}{2}\left(1-\frac{N+2}{Nq}\right)}\Phi(t^{-1})\int_0^{t/2}s^{-\frac{N}{2}(p-1)}\,ds\\
 & \le Ct^{-\frac{N}{2}\left(p-\frac{N+2}{Nq}\right)+1}\Phi((t/2)^{-1})
 \le Ct^\delta\Phi(t^{-1}),\quad t\in(0,T).
\end{split}
\end{equation}
Similarly, 
if $0<p\le 1$, then 
\begin{equation}
\label{eq:6.13}
\begin{split}
 & \biggr|\biggr|\biggr|\int_0^{t/2} S(D_1(t-s))v_{n_*}(s)^p\,ds\,\biggr|\biggr|\biggr|_{\Phi,r_*,r_*}\\
 & \le C\int_0^{t/2}(t-s)^{-\frac{N}{2}\left(p-\frac{N+2}{Nq}\right)}\Phi((t-s)^{-1})\|v_{n_*}(s)^p\|_{L^{\frac{1}{p}}_{{\rm ul}}}\,ds\\
 & \le C\int_0^{t/2}(t-s)^{-\frac{N}{2}\left(p-\frac{N+2}{Nq}\right)}\Phi((t-s)^{-1})\,ds\\
 & \le Ct^{-\frac{N}{2}\left(p-\frac{N+2}{Nq}\right)+1}\Phi((t/2)^{-1})\le Ct^\delta\Phi(t^{-1}),\quad t\in(0,T).
\end{split}
\end{equation}
By \eqref{eq:6.11}, \eqref{eq:6.12}, and \eqref{eq:6.13} we obtain 
$$
\biggr|\biggr|\biggr|\int_0^t S(D_1(t-s))v_{n_*}(s)^p\,ds\,\biggr|\biggr|\biggr|_{\Phi,r_*,1}\le Ct^\delta\Phi(t^{-1}),\quad t\in(0,T).
$$
Thus \eqref{eq:6.8} holds for $n=n_*$ in the case of $r_*>1/p$. 

On the other hand, if $r_*=1/p$, 
then $0<p\le 1$ and Proposition~\ref{Proposition:3.2} together with Lemma~\ref{Lemma:3.6} and \eqref{eq:6.10} implies that 
\begin{equation*}
\begin{split}
 & \biggr|\biggr|\biggr|\int_0^t S(D_1(t-s))v_{n_*}(s)^p\,ds\,\biggr|\biggr|\biggr|_{\Phi,r_*,r_*}
 =\,\biggr|\biggr|\biggr|\int_0^t S(D_1(t-s))v_{n_*}(s)^p\,ds\,\biggr|\biggr|\biggr|_{\Phi,\frac{1}{p},\frac{1}{p}}\\
 & \qquad\quad
 \le C\int_0^t\Phi((t-s)^{-1})\|v_{n_*}(s)^p\|_{L^{\frac{1}{p}}_{{\rm ul}}}\,ds
 \le C\int_0^t\Phi((t-s)^{-1})\,ds\\
 & \qquad\quad 
 \le Ct\Phi(t^{-1})
 =Ct^\delta\Phi(t^{-1}),\quad t\in(0,T).
\end{split}
\end{equation*}
This implies \eqref{eq:6.8} with $n=n_*$ in the case of $r_*=1/p$. 
Thus \eqref{eq:6.8} holds for $n=n_*$. 
\vspace{3pt}
\newline
\underline{Step 2.} 
We prove that \eqref{eq:6.4} and \eqref{eq:6.5} hold with $n=n_*+1$. 
It follows from \eqref{eq:6.10} that
\begin{equation}
\label{eq:6.14}
\left\|\,\int_{t/2}^t S(D_1(t-s))v_{n_*}(s)^p\,ds\,\right\|_{L^\infty}
\le C\int_{t/2}^t \|v_{n_*}(s)^p\|_{L^\infty}\,ds\le Ct^{-\frac{N}{2}p+1}
\end{equation}
for $t\in(0,T)$. 
Furthermore, if $p\ge 1$, by \eqref{eq:6.1} and \eqref{eq:6.10} we have 
\begin{equation}
\label{eq:6.15}
\left\|\,\int_0^{t/2} S(D_1(t-s))v_{n_*}(s)^p\,ds\,\right\|_{L^\infty}
\le C\int_0^{t/2}(t-s)^{-\frac{N}{2}}\|v_{n_*}(s)^p\|_{L^1_{{\rm ul}}}\,ds
\le Ct^{-\frac{N}{2}p+1}
\end{equation}
for $t\in(0,T)$. 
Similarly, if $0<p<1$, then 
\begin{equation}
\label{eq:6.16}
\left\|\,\int_0^{t/2} S(D_1(t-s))v_{n_*}(s)^p\,ds\,\right\|_{L^\infty} \le C\int_0^{t/2}(t-s)^{-\frac{N}{2}p}\|v_{n_*}(s)^p\|_{L^{\frac{1}{p}}_{{\rm ul}}}\,ds=Ct^{-\frac{N}{2}p+1}
\end{equation}
for $t\in(0,T)$.
By \eqref{eq:6.2}, \eqref{eq:6.14}, \eqref{eq:6.15}, and \eqref{eq:6.16} we see that 
$$
t^{\frac{N+2}{2q}}\Phi(t^{-1})\left\|\,\int_0^t S(D_1(t-s))v_{n_*}(s)^p\,ds\,\right\|_{L^\infty}
\le Ct^{\frac{N+2}{2q}-\frac{N}{2}p+1}\Phi(t^{-1})
\le Ct^\delta\Phi(t^{-1})
$$
for $t\in(0,T)$. 
Then, taking small enough $T>0$ if necessary, 
by Lemma~\ref{Lemma:3.5}~(2) and \eqref{eq:6.3} we obtain
$$
t^{\frac{N+2}{2q}}\Phi(t^{-1})\|u_{n_*}(t)\|_{L^\infty}\le C_*+Ct^\delta\Phi(t^{-1})
\le C_*+Ct^{\frac{\delta}{2}}
\le 2C_*,\quad t\in(0,T).
$$
Thus \eqref{eq:6.5} holds for $n=n_*+1$. 

Similarly, since $r^*>r_*\ge 1/p$, we find $m>1$ such that 
$$
\frac{1}{p}<m<r^*,\quad \frac{N}{2}\left(\frac{1}{m}-\frac{1}{r^*}\right)<1.
$$
It follows from a similar argument to those of \eqref{eq:6.15} and \eqref{eq:6.16} that
$$
\left\|\,\int_0^{t/2} S(D_1(t-s))v_{n_*}(s)^p\,ds\,\right\|_{L^{r^*}_{{\rm ul}}}
\le Ct^{-\frac{N}{2}\left(p-\frac{1}{r_*}\right)+1}
$$
for $t\in(0,T)$. 
Then, by \eqref{eq:6.10} we have 
\begin{equation}
  \notag 
  \begin{aligned}
    & 
    \left\|\,\int_{0}^t S(D_1(t-s))v_{n_*}(s)^p\,ds\,\right\|_{L^{r^*}_{{\rm ul}}}
    \\
    &
    \le 
    Ct^{-\frac{N}{2}\left(p-\frac{1}{r^*}\right)+1}
    + 
    C\int_{t/2}^t (t-s)^{-\frac{N}{2}\left(\frac{1}{m}-\frac{1}{r^*}\right)}\|v_{n_*}(s)^p\|_{L^m_{{\rm ul}}}\,ds
    \le Ct^{-\frac{N}{2}\left(p-\frac{1}{r^*}\right)+1}
  \end{aligned}
\end{equation}
for $t\in(0,T)$. 
Taking small enough $T>0$ if necessary, by Lemma~\ref{Lemma:3.5}~(2), \eqref{eq:6.2},
and \eqref{eq:6.3} we obtain 
\begin{align*}
t^{\frac{N}{2}\left(\frac{N+2}{Nq}-\frac{1}{r^*}\right)}\Phi(t^{-1}) \|u_{n_*+1}(t)\|_{L^{r^*}_{{\rm ul}}}
 & \le C_*+Ct^{-\frac{N}{2}\left(p-\frac{N+2}{Nq}\right)+1}\Phi(t^{-1})\\
 & \le C_*+Ct^{\delta}\Phi(t^{-1})\le C_*+Ct^{\frac{\delta}{2}}\le 2C_*,
\quad t\in(0,T).
 \end{align*}
Thus \eqref{eq:6.4} holds for $n=n_*+1$. 
\vspace{3pt}
\newline
\underline{Step 3.} 
We prove \eqref{eq:6.9} with $n=n_*$. 
Since $q>r^*$, it follows from \eqref{eq:6.4} and \eqref{eq:6.5} with $n=n_*$ that
\begin{align*}
 & \left\|\,\int_0^t S(D_2(t-s))u_{n_*}(s)^q\,ds\,\right\|_{L^1_{{\rm ul}}}\\
 & \le C\int_0^t \|u_{n_*}(s)^q\|_{L^1_{{\rm ul}}}\,ds
\le  C\int_0^t \|u_{n_*}(s)\|_{L^\infty}^{q-r^*}\|u_{n_*}(s)\|_{L^{r^*}_{{\rm ul}}}^{r^*}\,ds\\
 & \le C\int_0^t
 \left(s^{-\frac{N+2}{2q}}\Phi(s^{-1})^{-1}\right)^{q-r^*}
 \left(s^{-\frac{N}{2}\left(\frac{N+2}{Nq}-\frac{1}{r^*}\right)}\Phi(s^{-1})^{-1}\right)^{r^*}\,ds\\
 & \le C\int_0^t s^{-1}\Phi(s^{-1})^{-q}\,ds,\quad t\in(0,T).
\end{align*}
This implies \eqref{eq:6.9} with $n=n_*$. 
Furthermore, taking small enough $T$ if necessary, by \eqref{eq:1.18} we obtain 
$$
\|v_{n_*+1}(t)\|_{L^1_{{\rm ul}}}\le C_*+C\int_0^t s^{-1}\Phi(s^{-1})^{-q}\,ds\le 2C_*,\quad t\in(0,T).
$$
Thus \eqref{eq:6.6} holds for $n=n_*+1$. 
Similarly, taking small enough $T$ if necessary, we see that
\begin{align*}
t^{\frac{N}{2}}\|v_{n_*+1}(t)\|_{L^\infty} & \le C_*+Ct^{\frac{N}{2}}\int_0^{t/2} (t-s)^{-\frac{N}{2}}\|u_{n_*}(s)^q\|_{L^1_{{\rm ul}}}\,ds
 +Ct^{\frac{N}{2}}\int_{t/2}^t \|u_{n_*}(s)^q\|_{L^\infty}\,ds\\
  & \le C_*+C\int_0^{t/2} s^{-1}\Phi(s^{-1})^{-q}\,ds+Ct^{\frac{N}{2}}\int_{t/2}^t\left(s^{-\frac{N+2}{2q}}\Phi(s^{-1})^{-1}\right)^q\,ds\\
  & \le C_*+\int_0^t s^{-1}\Phi(s^{-1})^{-q}\,ds\le 2C_*,\quad t\in(0,T).
\end{align*}
Thus \eqref{eq:6.7} holds for $n=n_*+1$. 
The proof of Lemma~\ref{Lemma:6.1} is complete.
$\Box$\vspace{5pt}

\noindent
{\bf Proofs of Theorems~\ref{Theorem:1.5} and \ref{Theorem:1.6}.}
Similarly to the proof of Theorem~\ref{Theorem:1.3}, 
by Lemma~\ref{Lemma:6.1} we find a solution~$(u,v)$ to problem~\eqref{eq:P} in ${\mathbb R}^N\times(0,T)$
for some $T>0$. Furthermore, the solution~$(u,v)$ satisfies \eqref{eq:6.4}--\eqref{eq:6.9} with $(u_n,v_n)$ replaced by $(u,v)$. 
Then we deduce from~\eqref{eq:3.8} that $(u,v)$ is the desired solution. 
Thus Theorems~\ref{Theorem:1.5} and \ref{Theorem:1.6} follows. 
$\Box$
\section{Discussions}
Taking into the account of Proposition~\ref{Proposition:1.1}, 
we discuss the optimality of assumptions in our theorems. 
We remark that, in cases~(B)--(F), problem~\eqref{eq:P} possesses no global-in-time positive solutions 
(see assertion~(3) in Section~1). 
\vspace{3pt}
\newline
\underline{Case (A)}:
Consider case~(A). 
Set 
\begin{align*}
      & \mu(x)=c_{a,1} |x|^{-\frac{2(p+1)}{pq-1}}\chi_{B(0,1)}(x)\quad\mbox{in}\quad{\mathbb R}^N,\\
      & \nu(x)=c_{a,2} |x|^{-\frac{2(q+1)}{pq-1}}\chi_{B(0,1)}(x)\quad\mbox{in}\quad{\mathbb R}^N,
\end{align*}
where $c_{a,1}$, $c_{a,2}>0$. 
Let $\alpha_A$ and $\beta_B$ be as in \eqref{eq:1.5}. 
Then 
$$
\|\mu\|_{M(r_1^*,\alpha_A;\infty)}=C_1c_{a,1},
\quad
\|\nu\|_{M(r_2^*,\beta_A;\infty)}=C_1'c_{a,2},
$$
where $C_1$ and $C_1'$ are independent of $c_{a,1}$ and $c_{a,2}$.
Then, if $c_{a,1}$ and $c_{a,2}$ are small enough, 
Theorem~\ref{Theorem:1.1} implies that problem~\eqref{eq:P} possesses a global-in-time solution. 
On the other hand, if either $c_{a,1}$ or $c_{a,2}$ is large enough, then 
Proposition~\ref{Proposition:1.1}~(a) implies that problem~\eqref{eq:P} possesses no local-in-time solutions. 
This means that, if the constant $\delta_A$ in Theorem~\ref{Theorem:1.1} is large enough, 
then problem~\eqref{eq:P} does not necessarily possess local-in-time solutions.
\vspace{3pt}
\newline
\underline{Case (B)}:
Consider case~(B). Set 
\begin{align*}
      & \mu(x)=c_{b,1} |x|^{-\frac{2(p+1)}{pq-1}}\left[\log\left(e+\frac{1}{|x|}\right)\right]^{-\frac{p}{pq-1}}\chi_{B(0,1)}(x)
      \quad\mbox{in}\quad{\mathbb R}^N,\\
      & \nu(x)=c_{b,2} |x|^{-N}\left[\log\left(e+\frac{1}{|x|}\right)\right]^{-\frac{1}{pq-1}-1}\chi_{B(0,1)}(x)
     \quad\,\,\,\,\mbox{in}\quad{\mathbb R}^N,
\end{align*}
where $c_{b,1}$, $c_{b,2}>0$. 
Then 
\begin{align*}
 & \mu^*(s)\asymp c_{b,1}s^{-\frac{1}{N}\frac{2(p+1)}{pq-1}}\left[\log\left(e+\frac{1}{s}\right)\right]^{-\frac{p}{pq-1}}\chi_{(0,\omega_N)}(s)\quad\mbox{for $s>0$},\\
 & \nu^{*}(s)\asymp c_{b,2}s^{-1}\left[\log\left(e+\frac{1}{s}\right)\right]^{-\frac{1}{pq-1}-1}\chi_{(0,\omega_N)}(s)\quad\mbox{for $s>0$},\\
 & \nu^{**}(s)\asymp 
 \left\{
 \begin{array}{ll}
 \displaystyle{c_{b,2}s^{-1}\left[\log\left(e+\frac{1}{s}\right)\right]^{-\frac{1}{pq-1}}} & \quad\mbox{for $s\in(0,\omega_N)$},\vspace{5pt}\\
 c_{b,2}s^{-1}& \quad\mbox{for $s\in[\omega_N,\infty)$}.\vspace{3pt}
 \end{array}\right. 
\end{align*}
These imply that 
$$
\|\mu\|_{\frac{q+1}{p+1},\alpha_B;1}+|||\nu|||_{1,\beta_B;1}=C_2(c_{b,1}+c_{b,2}),
$$
where $C_2$ is independent of $c_{a,1}$ and $c_{a,2}$. 
Then, if $c_{b,1}$ and $c_{b,2}$ are small enough, then 
Theorem~\ref{Theorem:1.3} implies that problem~\eqref{eq:P} possesses a local-in-time solution. 
On the other hand, if either $c_{a,1}$ or $c_{a,2}$ is large enough, then 
Proposition~\ref{Proposition:1.1}~(b) implies that problem~\eqref{eq:P} possesses no local-in-time solutions. 
This means that, if the constant $\delta_B$ in Theorem~\ref{Theorem:1.3} is large enough, 
then problem~\eqref{eq:P} does not necessarily possess local-in-time solutions.
\vspace{3pt}
\newline
\underline{Case (C)}:
Consider case~(C). 
Set 
\begin{align*}
      & \mu(x)=c_{c,1}|x|^{-N}\left[\log\left(e+\frac{1}{|x|}\right)\right]^{-\frac{N}{2}-1}\chi_{B(0,1)}(x)\quad\mbox{in}\quad{\mathbb R}^N,\\
      & \nu(x)=c_{c,2}|x|^{-N}\left[\log\left(e+\frac{1}{|x|}\right)\right]^{-\frac{N}{2}-1}\chi_{B(0,1)}(x)\quad\mbox{in}\quad{\mathbb R}^N,
\end{align*}
where $c_{c,1}$, $c_{c,2}>0$. 
Then 
\begin{align*}
 & c_{c,1}^{-1}\mu^*(s)
 =c_{c,2}^{-1}\nu^*(s)
 \asymp 
s^{-1} \left[\log\left(e+\frac{1}{s}\right)\right]^{-\frac{N}{2}-1}\chi_{(0,\omega_N)}(s)
 \quad\mbox{for}\quad s>0,\\
 & c_{c,1}^{-1}\mu^{**}(s)=c_{c,2}^{-1}\nu^{**}(s)\asymp 
 \left\{
 \begin{array}{ll}
s^{-1}\displaystyle{\left[\log\left(e+\frac{1}{s}\right)\right]^{-\frac{N}{2}}} & \mbox{for}\quad s\in(0,\omega_N),\vspace{5pt}\\
s^{-1} & \mbox{for}\quad s\in[\omega_N,\infty).\vspace{3pt}
 \end{array}\right.
\end{align*}
These imply that 
$$
|||\mu|||_{1,\frac{N}{2};1}+|||\nu|||_{1,\frac{N}{2};1}=C_3(c_{c,1}+c_{c,2}),
$$
where $C_3$ is independent of $c_{a,1}$ and $c_{a,2}$. 
Then, if $c_{c,1}$ and $c_{c,2}$ are small enough, then 
Theorem~\ref{Theorem:1.4} implies that problem~\eqref{eq:P} possesses a local-in-time solution. 
On the other hand, if either $c_{c,1}$ or $c_{c,2}$ is large enough, then 
Proposition~\ref{Proposition:1.1}~(c) implies that problem~\eqref{eq:P} possesses no local-in-time solutions. 
This means that, if the constant $\delta_C$ in Theorem~\ref{Theorem:1.4} is large enough, 
then problem~\eqref{eq:P} does not necessarily possess local-in-time solutions.
\vspace{3pt}
\newline
\underline{Case (D)}:
Consider case~(D). 
Let $\Phi$ be a non-decreasing function in $[0,\infty)$ with properties $(\Phi 1)$--$(\Phi3)$. 
Let $\nu\in{\mathcal M}$ and set 
$$
\mu(x)=|x|^{-\frac{N+2}{q}}\Phi(|x|^{-1})^{-1}\chi_{B(0,1)}(x)\quad\mbox{in}\quad{\mathbb R}^N.
$$
It follows from Lemma~\ref{Lemma:3.5} that
$$
\mu^*(s)\asymp s^{-\frac{N+2}{Nq}}\Phi(s^{-1})^{-1}\chi_{(0,\omega_N)}(s),\quad s>0.
$$
This implies that 
$$
\mu\in L_{{\rm ul}}^{\frac{Nq}{N+2},\infty}\Phi(L)^{\frac{Nq}{N+2}}. 
$$
Then Theorem~\ref{Theorem:1.5} implies that problem~\eqref{eq:P} possesses a local-in-time solution if 
\begin{equation}
\label{eq:7.1}
\int_0^1 s^{-1}\Phi(s^{-1})^{-q}\,ds<\infty
\quad\mbox{and}\quad
\nu\in{\mathcal M}_{{\rm ul}}.
\end{equation}
Next, we assume that $r^{-\epsilon}\Phi(r^{-1})^{-1}$ is decreasing in $(0,1)$ for some $\epsilon>0$. 
Then Proposition~\ref{Proposition:1.1}~(d) implies that, 
if either
$$
\int_0^1 s^{-1}\Phi(s^{-1})^{-q}\,ds=\infty
\quad\mbox{or}\quad
\nu\not\in{\mathcal M}_{{\rm ul}},
$$
then problem~\eqref{eq:P} does not possess no local-in-time solutions. 
Thus problem~\eqref{eq:P} does not necessarily possess a local-in-time solution without \eqref{eq:7.1}.
\vspace{3pt}
\newline
\underline{Case (E)}:
Consider case~(E). 
Let $\Psi$ be a non-decreasing function in $[0,\infty)$ with $(\Phi 1)$--$(\Phi3)$ such that 
\begin{equation}
\label{eq:7.2}
\int_0^1 \tau^{-1}\Psi(\tau^{-1})^{-1}\,d\tau=1. 
\end{equation}
Let $\nu\in{\mathcal M}$ and set 
$$
\mu(x)=|x|^{-N}\Psi(|x|^{-1})^{-1}\chi_{B(0,1)}(x)\quad\mbox{in}\quad{\mathbb R}^N.
$$
It follows from Lemma~\ref{Lemma:3.5} that
$$
\mu^*(s)\asymp s^{-1}\Psi(s^{-1})^{-1}\chi_{(0,\omega_N)}(s),\quad 
\mu^{**}(s)\asymp s^{-1}\int_0^s \tau^{-1}\Psi(\tau^{-1})^{-1}\chi_{(0,\omega_N)}(\tau)\,d\tau,
$$
for $s>0$. 
Set 
$$
\Phi(s)\coloneqq
\left(\int_0^{s^{-1}} \tau^{-1}\Psi(\tau^{-1})^{-1}\chi_{(0,1)}(\tau)
\,d\tau\right)^{-1},\quad s>0.
$$ 
Then $\Phi$ is a non-decreasing function in $[0,\infty)$ and $\Phi(0)=1$ by \eqref{eq:7.2}. 
Furthermore, by assumption~$(\Phi 2)$ for $\Psi$ we have 
$$
    \Phi(a^2) 
    = 
    \left( 2 \int_0^{a^{-1}} \tau^{-1}\Psi(\tau^{-2})^{-1}\chi_{(0,1)}(\tau) \,d\tau\right)^{-1}
    \le C\Phi(a),\quad a>0,
$$
which implies that $\Phi$ satisfies $(\Phi 2)$. 
In addition, for any $\delta>0$, by assumption~$(\Phi 3)$ for $\Psi$ 
we find $C_\delta>0$ such that 
$$
\left(\frac{\tau_2}{\tau_1}\tau^{-1}\right)^{\delta}\Psi\left(\frac{\tau_2}{\tau_1}\tau^{-1}\right)^{-1}
\ge C_\delta^{-1} \tau^{-\delta}\Psi(\tau^{-1})^{-1}
$$
for small enough $\tau>0$ and all $\tau_1$, $\tau_2>0$ with $\tau_1\le\tau_2$. 
This implies that 
$$
\tau^{-1}\Psi\left(\frac{\tau_2}{\tau_1}\tau^{-1}\right)^{-1}
\ge C_\delta^{-1}\frac{\tau_1^\delta}{\tau_2^\delta}\tau^{-1}\Psi(\tau^{-1})^{-1}
$$
for small enough $\tau>0$ and all $\tau_1$, $\tau_2>0$ with $\tau_1\le\tau_2$. 
Then 
\begin{equation}
  \notag 
  \begin{aligned}
    \tau_2^{-\delta} \Phi(\tau_2) 
    &=
    \tau_2^{-\delta} 
    \left(\int_0^{\tau_1^{-1}} \tau^{-1}\Psi\left(\frac{\tau_2}{\tau_1}\tau^{-1}\right)^{-1}
    \chi_{(0,\tau_1^{-1}\tau_2)}(\tau)
    \,d\tau\right)^{-1}
    \\
    &\le 
    C_\delta \tau_1^{-\delta}\left(\int_0^{\tau_1^{-1}} \tau^{-1}\Psi(\tau^{-1})^{-1}
    \chi_{(0,\tau_1^{-1}\tau_2)}(\tau)
    \,d\tau\right)^{-1}
    \le 
    C_\delta \tau_1^{-\delta}\Phi(\tau_1)
  \end{aligned}
\end{equation}
for large enough $\tau_1$, $\tau_2$ with $\tau_1\le\tau_2$. 
Thus $\Phi$ satisfies $(\Phi 3)$. 
Since $\mu\in {\mathfrak L}_{{\rm ul}}^{1,\infty}\Phi({\mathfrak L})$, 
we observe from Theorem~\ref{Theorem:1.6} that problem~\eqref{eq:P} possesses a local-in-time solution if 
\begin{equation}
\label{eq:7.3}
\int_0^1 s^{-1}\Phi(s^{-1})^{-q}\,ds<\infty
\quad\mbox{and}\quad
\nu\in{\mathcal M}_{{\rm ul}}.
\end{equation}
On the other hand, setting 
$$
h_2(|x|)\coloneqq\Psi(|x|^{-1})^{-1}, 
$$
by Proposition~\ref{Proposition:1.1}~(e) we see that, 
if either
$$
\int_0^1 s^{-1}\Phi(s^{-1})^{-q}\,ds=\infty
\quad\mbox{or}\quad
\nu\not\in{\mathcal M}_{{\rm ul}},
$$
then problem~\eqref{eq:P} possesses no local-in-time solutions.  
Thus problem~\eqref{eq:P} does not necessarily possess a local-in-time solution without \eqref{eq:7.3}.
\vspace{3pt}
\newline
\underline{Case (F)}:
Consider case~(F). 
In Theorem~\ref{Theorem:1.2} we obtain a local-in-time solution if $\mu$, $\nu\in{\mathcal M}_{{\rm ul}}$. 
On the other hand, 
we observe from Proposition~\ref{Proposition:1.1}~(f) that 
$\mu$, $\nu\in{\mathcal M}_{{\rm ul}}$ is a necessary and sufficient condition 
for problem~\eqref{eq:P} to possess a local-in-time solution.
\medskip

\noindent{\bf Acknowledgments.}
YF was supported in part  by JSPS KAKENHI Grant Number 23K03179.
KI and TK  were supported in part by JSPS KAKENHI Grant Number JP19H05599.
YF and  TK were also supported in part by JSPS KAKENHI Grant Number JP22KK0035. 
\medskip

\noindent
{\bf  Conflict of Interest.}
The authors state no conflict of interest. 

\end{document}